\documentclass[12pt]{article}

\usepackage{indentfirst}
\usepackage{amsfonts}
\usepackage{amsmath}
\usepackage{amssymb}
\usepackage{amsbsy, amsthm}
\usepackage[margin=1in]{geometry}

\newtheorem{dfn}{Definition}[section]
\newtheorem{theorem}[dfn]{Theorem}
\newtheorem{lemma}[dfn]{Lemma}

\newenvironment{pf}{\noindent{\bf Proof.}}
{\enspace\vrule height5pt depth0pt width5pt}
\newenvironment{pf1}[1]{\noindent{\bf Proof of #1.}}
{\enspace\vrule height5pt depth0pt width5pt}
\def\deg {{\rm deg}}

\def\cl {{\rm cl}}
\def\X {{\mathcal X}}

\def\T {{\mathcal T}}
\def\E {{\mathcal E}}
\def\Eb {\overline{{\mathcal E}}}

\def\Se {{\mathcal S}}
\def\L {{\mathcal L}}
\def\F {{\mathcal F}}
\def\R {{\mathcal R}}
\def\P {{\mathcal P}}
\def\D {{\mathcal D}}
\def\Q {{\mathcal Q}}
\def\M {{\mathcal M}}
\def\TM {{\mathcal T}{\mathcal M}}
\def\I {{\mathcal I}}
\def\W {{\mathcal W}}

\begin{document}
\title{Packing and Covering Immersions in 4-Edge-Connected Graphs}
\author{Chun-Hung Liu\thanks{chliu@math.tamu.edu. Current address: Department of Mathematics, Texas A\&M University, College Station, TX 77843--3368. Partially supported by NSF under Grant No.~DMS-1664593, DMS-1929851 and DMS-1954054.} \\
\small Department of Mathematics, \\
\small Princeton University,\\
\small Princeton, NJ 08544-1000, USA}


\maketitle

\begin{abstract}
\noindent
A graph $G$ contains another graph $H$ as an immersion if $H$ can be obtained from a subgraph of $G$ by splitting off edges and removing isolated vertices.
In this paper, we prove an edge-variant of the Erd\H{o}s-P\'{o}sa property with respect to the immersion containment in 4-edge-connected graphs.
More precisely, we prove that for every graph $H$, there exists a function $f$ such that for every 4-edge-connected graph $G$, either $G$ contains $k$ pairwise edge-disjoint subgraphs each containing $H$ as an immersion, or there exists a set of at most $f(k)$ edges of $G$ intersecting all such subgraphs.
This theorem is best possible in the sense that the 4-edge-connectivity cannot be replaced by the 3-edge-connectivity.
\end{abstract}

\section{Introduction}
In this paper {\em graphs} are finite and are permitted to have loops and parallel edges.
Many questions in graph theory or combinatorial optimization can be formulated as follows.
Given a set of graphs $\F$ and a graph $G$, what is the maximum number of disjoint subgraphs of $G$ where each isomorphic to a member of $\F$ or what is the minimum number of vertices that are required to meet all such subgraphs?
We call the former problem the {\it packing problem} and the maximum number the {\it packing number}, and we call the latter problem the {\it covering problem} and the minimum the {\it covering number}.
For example, if $\F$ consists of the graph $K_2$, then the packing number is the maximum size of a matching and the covering number is the minimum size of a vertex cover; if $\F$ is the set of cycles, the covering number is the minimum size of a feedback vertex set.

In view of combinatorial optimization, the packing problem and the covering problem can be formulated as integer programming problems.
And the covering problem is the dual of the packing problem.
Furthermore, it is easy to see that the packing number is at most the covering number.
On the other hand, it is natural to ask when the covering number can be bounded by a function of the packing number from above.
In other words, we hope that the optimal solutions of the packing problem and covering problem are bounded by functions of each other.

Formally, a set of graphs $\F$ has the {\em Erd\H{o}s-P\'{o}sa property} if for every integer $k$, there exists a number $f(k)$ such that for every graph $G$, either $G$ contains $k$ disjoint subgraphs each isomorphic to a member of $\F$, or there exists $Z \subseteq V(G)$ with $\lvert Z \rvert \leq f(k)$ such that $G-Z$ does not contain a subgraph isomorphic to a member of $\F$.
A classical result of Erd\H{o}s and P\'{o}sa \cite{ep} states that the set of cycles has the Erd\H{o}s-P\'{o}sa property.
Hence the packing number and the covering number for the set of cycles are tied together.

This theorem was later generalized by Robertson and Seymour in terms of graph minors.
A graph is a {\em minor} of another graph if the former can be obtained from a subgraph of the latter by contracting edges.
For every graph $H$, define $\M(H)$ to be the set of graphs containing $H$ as a minor.
Robertson and Seymour \cite{rs V} proved that $\M(H)$ has the Erd\H{o}s-P\'{o}sa property if and only if $H$ is planar.
So the aforementioned result of Erd\H{o}s and P\'{o}sa follows from the case that $H$ is the one-vertex graph with one loop.

A variant of the minor containment is the topological minor containment.
We say that a graph is a {\em topological minor} of another graph if the former can be obtained from a subgraph of the latter by repeatedly contracting edges incident with at least one vertex of degree two.
Minor containment and topological minor containment are closely related.
For example, they are equivalent for characterizing planar graphs.

However, minors and topological minors behave much differently with respect to the Erd\H{o}s-P\'{o}sa property.
For every graph $H$, define $\TM(H)$ to be the set of graphs containing $H$ as a topological minor.
Unlike graph minors, the Erd\H{o}s-P\'{o}sa property for $\TM(H)$ is not equivalent with the planarity of $H$.
The author, Postle and Wollan \cite{lpw} provided a characterization of graphs $H$ in which $\TM(H)$ has the Erd\H{o}s-P\'{o}sa property and proved that it is NP-hard to decide whether $\TM(H)$ has the Erd\H{o}s-P\'{o}sa property for the input graph $H$.

The topological minor relation can be equivalently defined as follows.
A graph $H$ with no isolated vertices is a topological minor of another graph $G$ if there exist an injection $\pi_V$ from $V(H)$ to $V(G)$ and a function $\pi_E$ that maps the edges $e$, say with ends $u,v$, of $H$ to paths in $G$ from $\pi_V(u)$ to $\pi_V(v)$ (if $e$ is a loop with the end $v$, then $\pi_E(e)$ is a cycle in $G$ containing $\pi_V(v)$) such that $\pi_E(e_1)$ and $\pi_E(e_2)$ are internally disjoint for distinct edges $e_1,e_2$ of $H$.
Note that we consider a loop as a cycle as well.

We say that a graph $H$ (allowing isolated vertices) is an {\em immersion} of a graph $G$ if the mentioned internally disjoint property is replaced by the edge-disjoint property.
Formally, an {\em $H$-immersion} in $G$ is a pair of functions $\Pi=(\pi_V,\pi_E)$ such that the following hold.
	\begin{itemize}
		\item $\pi_V$ is an injection from $V(H)$ to $V(G)$.
		\item $\pi_E$ maps $E(H)$ to the set of subgraphs of $G$ such that for every edge $e$ of $H$, if $e$ has distinct ends $x,y$, then $\pi_E(e)$ is a path with ends $\pi_V(x)$ and $\pi_V(y)$, and if $e$ is the loop with end $v$, then $\pi_E(e)$ is a cycle containing $\pi_V(v)$.
		\item If $e_1,e_2$ are distinct edges of $H$, then $\pi_E(e_1)$ and $\pi_E(e_2)$ are edge-disjoint.
	\end{itemize}
We denote the subgraph $\bigcup_{e \in E(H)}\pi_E(e) \cup \bigcup_{v \in V(H)}\pi_V(v)$ of $G$ by $\Pi(H)$.
We say that two $H$-immersions $\Pi=(\pi_V,\pi_E)$ and $\Pi'=(\pi'_V,\pi'_E)$ are {\it edge-disjoint} if $\bigcup_{W \in \pi_E(E(H))}E(W)$ is disjoint from $\bigcup_{W \in \pi_E'(E(H))}E(W)$. 
(In this paper, for any function $f$ and any subset $X$ of its domain, we define $f(X)$ to be the set $\{f(x): x \in X\}$.)
Equivalently, $\Pi$ and $\Pi'$ are edge-disjoint if and only if $\Pi(H)$ and $\Pi'(H)$ are edge-disjoint subgraphs of $G$.

As immersions consist of edge-disjoint paths, it is reasonable to ask for packing edge-disjoint copies of immersions instead of disjoint copies.
Furthermore, one vertex can meet more than one edge-disjoint immersion, so it is more natural to cover these edge-disjoint subgraphs by edges instead of by vertices.
This motivates an edge-variant of the Erd\H{o}s-P\'{o}sa property.

We say that a set $\F$ of graphs has the {\em edge-variant of the Erd\H{o}s-P\'{o}sa property} if for every integer $k$, there exists an integer $f(k)$ such that for every graph $G$, either $G$ contains $k$ edge-disjoint subgraphs each isomorphic to a member of $\F$, or there exists $Z \subseteq E(G)$ with $\lvert Z \rvert \leq f(k)$ such that $G-Z$ has no subgraph isomorphic to a member of $\F$.
Raymond, Sau and Thilikos \cite{rst-ep} proved that $\M(\theta_r)$ has the edge-variant of the Erd\H{o}s-P\'{o}sa property, where $\theta_r$ is the loopless graph on two vertices with $r$ edges.

For every graph $H$, define $\I(H)$ to be the set of graphs containing $H$ as an immersion.
$\I(H)$ does not have the edge-variant of the Erd\H{o}s-P\'{o}sa property for every graph $H$.
The necessary conditions for graphs $H$ for which $\TM(H)$ has the Erd\H{o}s-P\'{o}sa property mentioned in \cite{lpw} are necessary for graphs $H$ for which $\I(H)$ has the edge-variant of the Erd\H{o}s-P\'{o}sa property.
On the other hand, even though a family of graphs does not have the edge-variant of the Erd\H{o}s-P\'{o}sa property, this family possibly has this property if we restrict the host graphs to be members of a smaller class of graphs.
For example, the set of odd cycles does not have the edge-variant of the Erd\H{o}s-P\'{o}sa property, but Kawarabayashi and Kobayashi \cite{kk} proved that it has the edge-variant of the Erd\H{o}s-P\'{o}sa property in 4-edge-connected graphs.
We address the same direction in this paper and prove that for every graph $H$, $\I(H)$ has the edge-variant of the Erd\H{o}s-P\'{o}sa property in 4-edge-connected graphs.
In fact, we prove the following theorem that is slightly stronger than the previous statement.

\begin{theorem} \label{main intro}
For every graph $H$, there exists a function $f: {\mathbb N} \rightarrow {\mathbb N}$ such that for every graph $G$ whose every component is 4-edge-connected and for every positive integer $k$, either $G$ contains $k$ edge-disjoint subgraphs each containing $H$ as an immersion, or there exists $Z \subseteq E(G)$ with $\lvert Z \rvert \leq f(k)$ such that $G-Z$ has no $H$-immersion.
\end{theorem}

Recall that the necessary conditions for which $\TM(H)$ has the Erd\H{o}s-P\'{o}sa property mentioned in \cite{lpw} provide necessary conditions for which $\I(H)$ has the edge-variant of the Erd\H{o}s-P\'{o}sa property.
The constructions in \cite{lpw} show that the 4-edge-connectivity cannot be replaced by the 3-edge-connectivity.
In addition, Kakimura and Kawarabayashi \cite{kk1} proved that Theorem \ref{main intro} is true if $G$ is 4-edge-connected and $H$ is a complete graph; they also provided some example showing that 3-edge-connectivity is not enough.
On the other hand, Giannopoulou, Kwon, Raymond and Thilikos \cite{gkrt} proved that the requirement of being 4-edge-connected can be dropped if $H$ is a loopless connected planar subcubic graph.

We remark that statements analogous to Theorem \ref{main intro} with respect to minors and topological minors are not true.
In other words, there does not exist a constant $c$ such that for every graph $H$, $\M(H)$ and $\TM(H)$ have the Erd\H{o}s-P\'{o}sa property even if the host graphs are $c$-connected.
Let $H$ be a graph such that $\M(H)$ (or $\TM(H)$, respectively) does not have the Erd\H{o}s-P\'{o}sa property.
So there exists a positive integer $k$ such that for every positive integer $N$, there exists a graph $G_N$ that does not contain $k$ disjoint subgraphs where each of them is a member of $\M(H)$ (or $\TM(H)$, respectively) such that there exists no $Z \subseteq V(G)$ with $\lvert Z \rvert \leq N$ hitting all such subgraphs.
For positive integers $N,c$, let $G_{c,N}$ be the graph obtained from $G_N$ by adding $c$ new vertices and adding edges from these $c$ vertices to all vertices in $G_N$.
Then for any positive integers $c,N$, the graph $G_{c,N}$ is $c$-connected but does not contain $k+c$ disjoint subgraphs where each of them is a member of $\M(H)$ (or $\TM(H)$, respectively), and there exists no hitting set in $G_{c,N}$ with size at most $N$.
So there is no absolute constant $c$ that would ensure that $\M(H)$ and $\TM(H)$ have the Erd\H{o}s-P\'{o}sa property in $c$-connected graphs.

In fact, we prove a stronger version of Theorem \ref{main intro} (see Theorem \ref{EP comp near 4-edge-conn}).
Theorem \ref{EP comp near 4-edge-conn} states that the 4-edge-connectivity can be replaced by the condition of having no edge-cut of order exactly three. (In fact, Theorem \ref{EP comp near 4-edge-conn} is even slightly stronger than this.)
Note that every Eulerian graph has no edge-cut of order three, so Theorem \ref{EP comp near 4-edge-conn} implies that the edge-variant of the Erd\H{o}s-P\'{o}sa property holds if the host graphs are Eulerian graphs.
Recall that Kakimura and Kawarabayashi \cite{kk1} proved that Theorem \ref{main intro} is true if $G$ is 4-edge-connected and $H$ is a complete graph; we remark that what they actually proved is stronger: the 4-edge-connectivity of $G$ can be replaced by the condition that no ``minimal edge-cut'' has size three.
Our Theorem \ref{EP comp near 4-edge-conn} also implies the stronger setting in \cite{kk1}.
See the remark after the proof of Theorem \ref{EP comp near 4-edge-conn} for the details.

We also consider the following version of the half-integral packing problem in this paper.
For every loopless graph $H$, an {\em $H$-half-integral immersion} in $G$ is a pair of functions $(\pi_V,\pi_E)$ such that the following hold.
	\begin{itemize}
		\item $\pi_V$ is an injection from $V(H)$ to $V(G)$.
		\item $\pi_E$ maps every edge $e$ with ends $u,v$ of $G$ to a path in $G$ from $\pi_V(u)$ to $\pi_V(v)$.
		\item For every edge $e$ of $G$, there exist at most two edges $e_1,e_2$ of $H$ such that $e \in \pi_E(e_1)$ and $e \in \pi_E(e_2)$.
	\end{itemize}
The following theorem shows that the 4-edge-connectivity can be dropped if we consider the following version of half-integral packing of half-integral immersions.

\begin{theorem} \label{half-integral intro}
For every loopless graph $H$, there exists a function $f: {\mathbb N} \rightarrow {\mathbb N}$ such that for every graph $G$ and for every positive integer $k$, either 
	\begin{enumerate}
		\item $G$ contains $k$ $H$-half-integral immersions $(\pi_V^{(1)},\pi_E^{(1)}), ..., (\pi_V^{(k)},\pi_E^{(k)})$ such that for every edge $e$ of $G$, there exist at most two pairs $(i,e')$ with $1 \leq i \leq k$ and $e' \in E(H)$ such that $e \in \pi_E^{(i)}(e')$, or 
		\item there exists $Z \subseteq E(G)$ with $\lvert Z \rvert \leq f(k)$ such that $G-Z$ has no $H$-half-integral immersion.
	\end{enumerate}
\end{theorem}

We remark that a result about half-integral packing of topological minors was proved by the author \cite{l}, but the notion of the half-integral packing in \cite{l} is different from the one in above theorem.

More recent developments about the Erd\H{o}s-P\'{o}sa property can be found in a survey of Raymond and Thilikos \cite{rt}.

\subsection{Overview of this paper}

Now we roughly sketch the proof of Theorem \ref{main intro} and describe the organization of this paper.
More detailed sketches of proofs will be included in later sections when we are ready to prove them.
We need following ingredients for the proof of Theorem \ref{main intro}.

\subsubsection{Ingredient 1: edge-tangles}

Tangle is one important notion in Robertson and Seymour's Graph Minors series.
It defines a ``consistent orientation'' for each separation of small order in a graph and has been proven to be useful in dealing with problems related to the Erd\H{o}s-P\'{o}sa property. 
If one can delete a small number of vertices from a graph such that the packing number of every component of the remaining graph is smaller than packing number of the original graph, then one can obtain a hitting set of small size by an easy inductive argument.
So we can assume that there must exist a component whose packing number is not smaller than the packing number of the original graph.
Note that such a component is unique as we cannot have two components whose packing number are not smaller than the packing number of the original graph.
Hence, as long as we delete a small number of vertices, we know which component of the graph obtained by vertex-deletion is the ``most important''.
Given a separation of a graph, by simply seeing which side has this ``important'' component, we obtain an orientation for separations to define a tangle.
As we address the edge-variant of the Erd\H{o}s-P\'{o}sa property in this paper, we develop a similar machinery called ``edge-tangles'' which gives an orientation for each edge-cut of small order in Section \ref{sec:edge-tangles}. 
This section is a preparation of many results of this paper and includes formal definitions of tangles and edge-tangles.

The notion of edge-tangles is natural but its explicit form seems unnoticed by the community until the first version of this paper was submitted for publication.
We remark that Diestel and Oum \cite{do2,do} extended the ideas of ``consistent orientations'' for separations of graphs to a much more general setting for ``abstract separation systems'' to prove a general strong duality theorem.
After a version of this paper was submitted for publication, it was pointed out by Reinhard Diestel (via private communication with the author) that the concept of edge-tangles defined in this paper coincides with a special case of their abstract separation systems when graphs are loopless.
In particular, applying their duality theorem for abstract separation systems in \cite{do2} to edge-tangles, they \cite[Section 5.2]{do} noted that edge-tangles are dual to low ``carving width'' as defined by Seymour and Thomas \cite{st}.
The existence of such a duality theorem for edge-tangles was independently asked by Robin Thomas via private communication with the author around 2013 when the author introduced the notion of edge-tangles to develop a structure theorem for excluding immersions.
In addition, Diestel, Hundertmark and Lemanczyk \cite{dhl} applied their general work to edge-tangles to derive the classical Gomory-Hu tree theorem.
Besides those developments for abstract separation systems, in this paper we consider different aspects for edge-tangles, mainly on developing structure theorems with respect to immersions and edge-tangles and its application to Erd\H{o}s-P\'{o}sa type problems.
We omit the details and formal definitions of the terms mentioned in this paragraph, as this paper does not rely on them.

In Section \ref{sec:edge-tangles}, we develop basic theory related to edge-tangles.
In particular, in Sections \ref{subsec:edge-cut_line_graph} and \ref{subsec:basic_prop_edge-tangles}, we show some basic properties for edge-tangles and build a relationship between edge-tangles in a graph $G$ and tangles in the line graph of $G$.
Those results will be used in later sections of this paper.

It is known that if a graph $G$ contains a graph $H$ as a minor, then tangles in $H$ ``induce'' tangles in $G$; and every tangle of large order has a ``subtangle'' induced by a large grid minor.
We show analogous results in Section \ref{subsec:imm_edge-tangles}: if a graph $G$ contains a graph $H$ as an immersion, then edge-tangles in $H$ ``induce'' edge-tangles in $G$ (Lemma \ref{immersion induces edge-tangle}); every edge-tangle of large order has a ``sub-edge-tangle'' determined by a large degree vertex or induced by a large wall immersion (Lemma \ref{truncation large deg or wall}).

We prove other lemmas about edge-tangles in Section \ref{subsec:edge-tangles_others}.
Those lemmas will be used in later sections.

\subsubsection{Ingredient 2: a structure theorem for excluding immersions}

The next step to prove Theorem \ref{main intro} is to study the structure of minimum counterexamples to Theorem \ref{main intro}.
If $G$ is a graph that does not contain $k$ edge-disjoint $H$-immersions, then $G$ does not contain an $H'$-immersion for some larger graph $H'$.
So we shall prove a structure theorem for excluding a fixed graph as an immersion.
Such a structure theorem was developed by the author in an unpublished manuscript in 2013.
In this paper we include a proof of part of this theorem which is sufficient for proving Theorem \ref{main intro}.
More specifically, we will prove in this paper a structure theorem (Theorem \ref{excluding immersion}) for graphs that forbids a fixed graph as an immersion with respect to an edge-tangle that ``grasps'' a large set of pairwise edge-disjoint but pairwise intersecting subgraphs (called ``thorns'').
The formal definition of thorns is included in Section \ref{sec:excluding immersions structure}.

Roughly speaking, Theorem \ref{excluding immersion} states that if an $H$-immersion free graph $G$ has an edge-tangle of large order grasping a large thorns, then we can sweep all except few vertices into the non-important side of edge-cuts in this edge-tangle so that the ``resulting graph'' is ``simpler'' than $H$ in terms of the supply of large degree vertices.
To prove Theorem \ref{excluding immersion}, we first prove strengthenings of some Menger-type results of Robertson and Seymour \cite{rs XXIII} and Marx and Wollan \cite{mw} in Section \ref{sec:spider theorems}, and then we complete the proof of this structure theorem in Section \ref{sec:excluding immersions structure} by using results proved in Sections \ref{sec:edge-tangles} and \ref{sec:spider theorems}.

We remark that a structure theorem about excluding a fixed graph as an immersion was proved by Wollan \cite{w}.
But it seems to us that Wollan's structure theorem is not sufficiently informative to be applied in this paper.
In addition, the structure theorem proved in this paper is a ``local version'' and is a critical component of the proof of a ``global version'' of an excluding immersion structure theorem proved in a later paper of the author \cite{l_global}.
The global version has other applications, see \cite{l_global} for more details.

\subsubsection{Ingredient 3: 4-edge-connectivity and thorns}

Then, in Section \ref{sec:4-edge-connected edge-tangle} we will show that the 4-edge-connectivity of ``sufficiently large'' graphs can ensure that every edge-tangle ``grasps'' such a large thorns and hence the aforementioned structure theorem can be applied to minimum counterexamples to Theorem \ref{main intro}.
This is achieved by Theorem \ref{edge-tangle in 4-edge-connected control edge-minor}.
Recall that Lemma \ref{truncation large deg or wall} shows that every edge-tangle of large order has a sub-edge-tangle determined by a large degree vertex or induced by a large wall immersion.
Both a large degree vertex and a large grid immersion define a large thorns.
So the remaining key step in proving Theorem \ref{edge-tangle in 4-edge-connected control edge-minor} is to show that in any 4-edge-connected graph, one can obtain a large grid immersion from a large wall immersion.

\subsubsection{Ingredient 4: preserving edge-connectivity}

Finally, we prove Theorems \ref{main intro} and \ref{half-integral intro} in Section \ref{sec:EP main}.

Recall that if one can find an edge-cut of small order such that the subgraph induced on each side has smaller packing number for $H$-immersions, then one can apply induction on the packing number to obtain a small hitting set of each of these two subgraphs, and one can obtain a hitting set of the whole graph by further collecting the edges between the two sides of the edge-cut.
Hence in the minimum counterexample to Theorem \ref{main intro}, no edge-cut of small order has this property.
In particular, at most one side can contain an $H$-immersion.
Furthermore, when $H$ is connected, at least one side must contain an $H$-immersion, for otherwise the edges between the two sides is a small hitting set.
Hence we obtain an orientation of each edge-cut of small order by indicating the side having an $H$-immersion is important, so we obtain an edge-tangle.
Our structure theorem (Theorem \ref{excluding immersion}) ensures that for any $H$-immersion, one can find an edge-cut of small order such that the non-important side contains part of this $H$-immersion.
By repeatedly applying machinery developed in Section \ref{sec:spider theorems}, we can repeatedly enlarge the portion of this $H$-immersion contained in the non-important side until the entire $H$-immersion is contained in the non-important side, which is a contradiction by the definition of the edge-tangle.
This is the purpose of Section \ref{sec:isolating an immersion} and is achieved by Lemma \ref{isolating immersion}. 

In fact, careful readers might notice that there are issues with the above strategy connected to induction on the packing number and defining an edge-tangle. 
One concern is that we require $H$ to be connected.
This concern can be solved relatively easily by induction on the number of components of $H$, and Lemma \ref{isolating immersion} actually already takes care of it.
The other concern is more substantial: we tried to apply induction hypothesis to the subgraphs induced by each side of an edge-cut.
Note that such a subgraph is not necessarily 4-edge-connected, so we cannot apply induction to it.
The key strategy is to apply induction to the graph obtained from contracting one side of the edge-cut.
Such a graph preserves the 4-edge-connectivity, but the packing number is not necessarily smaller than the original graph even though the other side of the edge-cut contains an $H$-immersion.
To solve this issue, we actually prove a stronger version of Theorem \ref{main intro} by allowing some vertices having labels, where the label of each vertex can be roughly considered the number of $H$-immersions that this vertex represents.
More details are included in Section \ref{sec:EP main}.

Section \ref{sec:EP main} formally proves this stronger setting of Theorem \ref{main intro} and solves the aforementioned concern about preserving edge-connectivity.
Theorem \ref{EP comp near 4-edge-conn} is the strongest version in this paper and it implies Theorem \ref{main intro}.
We remark that we do not require the edge-connectivity in Section \ref{sec:isolating an immersion} and Lemma \ref{isolating immersion}, as they only require the edge-tangles to grasp a thorns.
Moreover, Theorem \ref{half-integral intro} is a simple corollary of this stronger version of Theorem \ref{main intro}, and its proof is included in Section \ref{sec:EP main}.

\subsection{Notations}

We define some notations to conclude this section.
Given a subset $X$ of the vertex-set $V(G)$ of a graph $G$, the subgraph of $G$ induced by $X$ is denoted by $G[X]$, and the set of vertices that are not in $X$ but adjacent to some vertices in $X$ is denoted by $N_G(X)$.
When $X = \{v\}$, we write $N_G(\{v\})$ as $N_G(v)$ for simplicity.
We define $N_G[X] = N_G(X) \cup X$ and $N_G[v] = N_G(v) \cup \{v\}$.
A graph is {\it simple} if it does not contain parallel edges and loops.
The {\it line graph} of a graph $G$, denoted by $L(G)$, is the simple graph with $V(L(G)) = E(G)$, and every pair of vertices $x,y \in V(L(G))$ are adjacent in $L(G)$ if and only if $x,y$ are two edges having a common end in $G$.
For every $v \in V(G)$, define $\cl(v)$ to be the clique in $L(G)$ consisting of the edges of $G$ incident with $v$.
The {\it degree} of a vertex $v$ in a graph $G$, denoted by $\deg_G(v)$, is the number of edges of $G$ incident with $v$, where each loop is counted twice.
A vertex of $G$ is an {\it isolated vertex} in $G$ if it is not incident with any edge.
Note that a vertex is non-isolated even if all edges incident with it are loops.
If $G$ is a graph and $Y \subseteq V(G)$, then $G-Y$ is the graph $G[V(G)-Y]$; if $Y \subseteq E(G)$, then $G-Y$ is the graph with $V(G-Y)=V(G)$ and $E(G-Y)=E(G)-Y$.
For a positive integer $k$, a graph $G$ is {\it $k$-edge-connected} if $G$ contains at least two vertices and $G-F$ is connected for every $F \subseteq E(G)$ with $\lvert F \rvert < k$.
For every positive integer $n$, we denote the set $\{1,2,...,n\}$ by $[n]$ for short.

\section{Tangles and edge-tangles}
\label{sec:edge-tangles}

\subsection{Edge-cuts and separations of line graphs} \label{subsec:edge-cut_line_graph}

A {\it separation} of a graph $G$ is an ordered pair $(A,B)$ of edge-disjoint subgraphs of $G$ with $A \cup B=G$, and the {\it order} of $(A,B)$ is $\lvert V(A) \cap V(B) \rvert$.

A separation $(A,B)$ of $G$ is {\it normalized} if every vertex $v \in V(A) \cap V(B)$ is adjacent to a vertex of $A-V(B)$ and adjacent to a vertex in $B-V(A)$.
The {\it normalization} of a separation $(A,B)$ of a graph $G$ is the separation $(A^*,B^*)$ of $G$ defined as follows.
	\begin{itemize}
		\item Let $S_1$ be the set of all non-isolated vertices $v$ in $G$ contained in $V(A) \cap V(B)$ with $N_A(v) \subseteq V(B)$.
		Let $A'$ be the graph $A-S_1$ and let $B'$ be the subgraph of $G$ such that $(A',B')$ is a separation of $G$ with $V(A') \cap V(B') = V(A) \cap V(B) - S_1$.
		In other words, $(A',B')$ is obtained from $(A,B)$ by removing all non-isolated vertices $v \in V(A) \cap V(B)$ of $G$ with $N_A(v) \subseteq V(B)$ from $A$ and putting all edges of $G$ incident with $v$ into $B$.
		Note that $V(B')=V(B)$.

		\item Let $S_2$ be the set of all edges in $B'$ whose every end is in $V(A') \cap V(B')$.
		Let $A''=A' \cup S_2$ and $B''=B'-S_2$.
		Note that $(A'',B'')$ is a separation of $G$, and every loop of $G$ incident with some vertex in $V(A'') \cap V(B'')$ belongs to $A''$.

		\item Let $S_3$ be the set of all isolated vertices of $B''$.
		Note that $S_3$ consists of the isolated vertices of $G$ contained in $V(B'')$ and some vertices contained in $V(A'') \cap V(B'')$ that are not adjacent to any vertex in $V(B'')-V(A'')$ by the definition of $A''$.
		Define $A^*$ to be the graph obtained from $A''$ by adding $S_3-V(A'')$, and define $B^*=B''-S_3$.
		That is, we remove all isolated vertices of $B''$ from $B''$ and put them into $A''$.
	\end{itemize} 
	
The following lemma shows some basic properties of the normalization of a separation and will be used in the rest of the section.

\begin{lemma} \label{normalization basic}
Let $G$ be a graph and $(A,B)$ a separation of $G$.
If $(A^*,B^*)$ is the normalization of $(A,B)$, then the following hold.
	\begin{enumerate}
		\item $(A^*,B^*)$ is normalized.
		\item The order of $(A^*,B^*)$ is at most the order of $(A,B)$.
		\item If $e \in E(B)-E(B^*)$, then every end of $e$ belongs to $V(A) \cap V(B)$.
		\item If $e \in E(B^*)-E(B)$, then $e$ is incident with some vertex $v$ in $V(A) \cap V(B)$ with $N_A(v) \subseteq V(B)$.
	\end{enumerate}
\end{lemma}

\begin{pf}
Let $G$ be a graph, $(A,B)$ a separation of $G$, and $(A^*,B^*)$ the normalization of $(A,B)$.
Let $S_1,S_2,S_3$ be the sets and let $(A',B'),(A'',B'')$ be the separations mentioned in the definition of $(A^*,B^*)$, respectively.
It is clear that $V(A^*) \cap V(B^*) \subseteq V(A'') \cap V(B'') = V(A') \cap V(B') = V(A) \cap V(B) - S_1$, and $V(A)-V(B) \subseteq V(A^*)-V(B^*)$.

We first prove Statement 1.
Let $v \in V(A^*) \cap V(B^*)$.
So $v$ is a non-isolated vertex of $G$ and $v \in V(A) \cap V(B)-S_1$.
Hence $N_A(v) \not \subseteq V(B)$.
That is, $v$ is adjacent to some vertex in $V(A)-V(B) \subseteq V(A^*)-V(B^*)$.
Since $v \in V(B^*)$, $v$ is not an isolated vertex of $B''$.
So $v$ is adjacent in $B''$ to some vertex $u \in V(B'')$.
Note that $u$ is not an isolated vertex of $B''$.
Hence $u \in V(B'')-S_3 =V(B^*)$.
Since $v \in V(A^*) \cap V(B^*) \subseteq V(A') \cap V(B')$, $u \not \in V(A') \cap V(B')$, for otherwise every edge incident with both $u,v$ belongs to $S_2$ and $u$ is not adjacent to $v$ in $B''$.
Since $V(A'') \cap V(B'') = V(A') \cap V(B')$, $u \in V(B'')-V(A'')$.
Since $u \not \in S_3$, $u \in V(B^*)-V(A^*)$.
So $v$ is adjacent to a vertex in $V(B^*)-V(A^*)$.
This shows that $(A^*,B^*)$ is normalized.

Statement 2 immediately follows from the fact that $V(A^*) \cap V(B^*) \subseteq V(A) \cap V(B)-S_1$.

Now we prove Statement 3.
Let $e \in E(B)-E(B^*)$.
Since $e \in E(B)$, $e \not \in E(A)$.
So $e \not \in E(A')$ and hence $e \in E(B')$.
Since every vertex in $S_3$ is an isolated vertex in $B''$ and $e \not \in E(B^*)$, $e \not \in E(B'')$.
So $e \in E(B')-E(B'') \subseteq S_2$.
Hence every end of $e$ belongs to $V(A') \cap V(B') \subseteq V(A) \cap V(B)$.
This shows Statement 3.

Finally, we prove Statement 4.
Let $e \in E(B^*)-E(B)$.
Since every vertex in $S_3$ is an isolated vertex of $B''$, $e \in E(B^*)=E(B'') \subseteq E(B')$.
So $e \not \in E(A')$.
Since $e \not \in E(B)$, $e \in E(A)$.
Hence $e$ is incident with some vertex $w$ in $S_1$.
But every vertex in $S_1$ satisfies that $w \in V(A) \cap V(B)$ and $N_A(w) \subseteq V(B)$.
This completes the proof.
\end{pf}

\bigskip

An {\it edge-cut} of a graph $G$ is an ordered partition $[A,B]$ of $V(G)$, where some of $A$ and $B$ is allowed to be empty.
The {\it order} of an edge-cut $[A,B]$, denoted by $\lvert [A,B] \rvert$, is the number of edges with one end in $A$ and one end in $B$.
For an edge $e$ of $G$, we write $e \in [A,B]$ if $e$ has one end in $A$ and one end in $B$.

The {\it partner} of a normalized separation $(A,B)$ of the line graph $L(G)$ of $G$ is the edge-cut $[A',B']$ of $G$ satisfying that $A'$ is the union of the set of isolated vertices of $G$ and the set $\{v \in V(G): \cl(v) \subseteq V(A)\}$, and $B'=\{v \in V(G): \cl(v) \subseteq V(B), \cl(v) \neq \emptyset\}$.

\begin{lemma} \label{partner of sep}
Let $G$ be a graph, and let $(A,B)$ be a separation of $L(G)$.
If $(A,B)$ is normalized, then the partner $[A',B']$ of $(A,B)$ is a well-defined edge-cut of $G$, and the order of $(A,B)$ equals the order of $[A',B']$.
\end{lemma}

\begin{pf}
We first show that the partner $[A',B']$ of $(A,B)$ is a well-defined edge-cut of $G$.
That is, $A' \cup B' = V(G)$ and $A' \cap B'=\emptyset$.
Let $v \in V(G)$.
If $v$ is an isolated vertex of $G$, then $v \in A'$; otherwise, $\cl(v)$ is a non-empty clique, so $\cl(v) \subseteq V(A)$ or $\cl(v) \subseteq V(B)$, and hence $v \in A' \cup B'$.
So $A' \cup B' = V(G)$.
Suppose to the contrary that $v \in A' \cap B'$.
Then $\cl(v) \subseteq V(A) \cap V(B)$ and $\cl(v) \neq \emptyset$.
So there exist $e_0 \in \cl(v)$ and a set $X \subseteq V(G)-\{v\}$ with $\lvert X \rvert \leq 1$ such that $N_{L(G)}(e_0) \subseteq \cl(v) \cup \bigcup_{u \in X}\cl(u)$.
But since $(A,B)$ is normalized and $e_0 \in \cl(v) \subseteq V(A) \cap V(B)$, $e_0$ is adjacent in $L(G)$ to a vertex in $A-V(B)$ and a vertex in $B-V(A)$.
Since $\cl(v) \subseteq V(A) \cap V(B)$, $\bigcup_{u \in X}\cl(u)$ intersects both $V(A)-V(B)$ and $V(B)-V(A)$.
But it is impossible since $\lvert X \rvert \leq 1$, a contradiction.
This shows that $[A',B']$ is an edge-cut of $G$.

Now we show that the order of $(A,B)$ equals the order of $[A',B']$.
Let $e \in [A',B']$ with ends $u,v$, where $u \in A'$ and $v \in B'$.
So $\cl(u) \subseteq V(A)$ and $\cl(v) \subseteq V(B)$.
Hence, $e \in \cl(u) \cap \cl(v) \subseteq V(A) \cap V(B)$.
This implies that the order of $(A,B)$ is at least the order of $[A',B']$.

On the other hand, let $e \in V(A) \cap V(B)$.
Since $(A,B)$ is normalized, $e$ is adjacent to a vertex $e_A$ of $L(G)$ in $V(A)-V(B)$ and a vertex $e_B$ of $L(G)$ in $V(B)-V(A)$.
So $e$ and $e_A$ have a common end $x$ in $G$, and $e$ and $e_B$ have a common end $y$ of $G$.
Since $e_A \not \in V(B)$, $\cl(x) \subseteq V(A)$ and hence $x \in A'$.
Similarly, $\cl(y) \subseteq V(B)$ and $y \in B'$.
So $x \neq y$ and they are the ends of $e$.
This proves that $e \in [A',B']$ and the order of $(A,B)$ is at most the order of $[A',B']$.
\end{pf}

\subsection{Basic properties of edge-tangles} \label{subsec:basic_prop_edge-tangles}

Let $\theta$ be an integer.
A {\it tangle} $\T$ in a graph $G$ of {\it order} $\theta$ is a set of separations of $G$ of order less than $\theta$ such that
\begin{enumerate}
	\item[(T1)] for every separation $(A,B)$ of $G$ of order less than $\theta$, either $(A,B) \in \T$ or $(B,A) \in \T$;
	\item[(T2)] if $(A_1, B_1), (A_2,B_2), (A_3,B_3) \in \T$, then $A_1 \cup A_2 \cup A_3 \neq G$;
	\item[(T3)] if $(A,B) \in \T$, then $V(A) \neq V(G)$.
\end{enumerate}
The notion of tangles was first defined by Robertson and Seymour in \cite{rs X}.
We call (T1), (T2) and (T3) the {\it first}, {\it second} and {\it third tangle axioms}, respectively.
Note that (T2) implies that if $(A,B) \in \T$, then $(B,A) \not \in \T$.

An {\it edge-tangle} $\E$ in a graph $G$ of {\it order} $\theta$ is a set of edge-cuts of $G$ of order less than $\theta$ such that the following hold.
\begin{enumerate}
	\item[(E1)] For every edge-cut $[A,B]$ of $G$ of order less than $\theta$, either $[A,B] \in \E$ or $[B,A] \in \E$;
	\item[(E2)] If $[A_1,B_1],[A_2,B_2],[A_3,B_3] \in \E$, then $B_1 \cap B_2 \cap B_3 \neq \emptyset$.
	\item[(E3)] If $[A,B] \in \E$, then $G$ has at least $\theta$ edges incident with vertices in $B$.
\end{enumerate}
We call (E1), (E2) and (E3) the {\it first}, {\it second} and {\it third edge-tangle axioms}, respectively.
Note that if an edge-tangle $\E$ of order $\theta \geq 1$ in $G$ exists, then $[\emptyset, V(G)] \in \E$ by (E1) and (E2), so $\lvert E(G) \rvert \geq \theta$ by (E3).
Furthermore, for every $[A,B] \in \E$, there exists an edge of $G$ whose every end is in $B$ by (E3).

The following lemma is simple but useful.
It shows that the ``orientation'' of the edge-cuts given by an edge-tangle is ``consistent''.

\begin{lemma} \label{easy edge-tangle}
Let $\theta$ be a positive integer.
Let $G$ be a graph and $\E$ an edge-tangle of order $\theta$ in $G$.
If $[A,B],[C,D] \in \E$, then the following hold.
	\begin{enumerate}
		\item If the order of $[A \cup C, B \cap D]$ is less than $\theta$, then $[A \cup C, B \cap D] \in \E$.
		\item If $A' \subseteq A$ and $[A',V(G)-A']$ is an edge-cut of $G$ of order less than $\theta$, then $[A',V(G)-A'] \in \E$.
	\end{enumerate}
\end{lemma}

\begin{pf}
We first prove Statement 1.
Assume that $[A \cup C, B \cap D]$ has order less than $\theta$.
By (E1), either $[A \cup C, B \cap D] \in \E$ or $[B \cap D, A \cup C] \in \E$.
Since $[A,B],[C,D] \in \E$ and $B \cap D \cap (A \cup C)=\emptyset$, $[B \cap D, A \cup C] \not \in \E$ by (E2).
So $[A \cup C, B \cap D] \in \E$.
This shows Statement 1.

Now we prove Statement 2.
Let $A' \subseteq A$, and assume that $[A',V(G)-A']$ is an edge-cut of $G$ of order less than $\theta$.
By (E1), either $[A',V(G)-A'] \in \E$ or $[V(G)-A',A'] \in \E$.
Since $[A,B] \in \E$ and $B \cap A'=\emptyset$, $[V(G)-A',A'] \not \in \E$ by (E2).
So $[A',V(G)-A'] \in \E$.
This shows Statement 2.
\end{pf}

\bigskip

The following lemma shows that the vertex-set of any component with at most one edge can be moved to either side of an edge-cut without flipping the ``orientation'' given by an edge-tangle.

\begin{lemma} \label{moving isolated edges}
Let $\theta$ be an integer with $\theta \geq 2$ and let $G$ be a graph.
Then the following hold.
	\begin{enumerate}
		\item Let $S$ be the vertex-set of a component of $G$ with at most one edge.
		If $[A,B] \in \E$, then $[A \cup S, B-S] \in \E$.
		\item Let $D$ be the union of the vertex-sets of the components of $G$ with at most one edge.
		If $[A,B] \in \E$, then $[A \cup D, B-D] \in \E$.
	\end{enumerate}
\end{lemma}

\begin{pf}
We first prove Statement 1.
Suppose to the contrary that $[A \cup S,B-S] \not \in \E$.
Note that every edge with one end in $A \cup S$ and one end in $B-S$ is an edge with one end in $A$ and one end in $B$.
Hence the order of $[A \cup S,B-S]$ is at most the order of $[A,B]$.
By (E1), $[B-S,A \cup S] \in \E$.
The edge-cut $[S,V(G)-S]$ has order 0, so either $[S,V(G)-S] \in \E$ or $[V(G)-S,S] \in \E$ by (E1).
Since $\theta \geq 2$ and there exists at most one edge incident with $S$, (E3) implies that $[V(G)-S,S] \not \in \E$.
So $[S, V(G)-S] \in \E$.
Hence $[A,B], [B-S,A \cup S], [S,V(G)-S]$ are edge-cuts in $\E$ such that $B \cap (A \cup S) \cap V(G)-S = \emptyset$, contradicting (E2).
This proves Statement 1.

Now we prove Statement 2.
Let $S_1,S_2,...,S_k$ be the subsets of $V(G)$ such that each $S_i$ is the vertex-set of some component of $G$ with at most one edge.
So $D = \bigcup_{i=1}^k S_i$.
For every $j \in [k]$, let $D_j = \bigcup_{i=1}^j S_i$.
We shall prove that $[A \cup D_j, B-D_j] \in \E$ for $j \in [k]$ by induction on $j$.
The case $j=1$ immediately follows from Statement 1 of this lemma.
So we may assume that $j \geq 2$ and $[A \cup D_{j-1},B-D_{j-1}] \in \E$.
Applying Statement 1 of this lemma by taking $[A,B]=[A \cup D_{j-1},B-D_{j-1}]$ and $S=S_j$, we know that $[A \cup D_j, B-D_j] = [(A \cup D_{j-1}) \cup S_j, (B-D_{j-1})-S_j] \in \E$.
This shows that $[A \cup D_j, B-D_j] \in \E$ for every $j \in [k]$.
Hence $[A \cup D, B-D] = [A \cup D_k, B-D_k] \in \E$.
This proves the lemma.
\end{pf}

\bigskip

The next step is to build a relationship between edge-tangles in $G$ and tangles in its line graph $L(G)$.

Given an edge-tangle $\E$ of order $\theta$ in $G$, the {\it conjugate} $\Eb$ of $\E$ is the set of separations of $L(G)$ of order less than $\lceil \theta/3 \rceil$ such that $(A,B) \in \Eb$ if and only if the partner of the normalization of $(A,B)$ is in $\E$.

One reason for considering separations of $L(G)$ of order less than $\lceil \theta/3 \rceil$ only instead of considering separations of $L(G)$ of order less than $\theta$ is due to a technicality in the proof of the following lemma which shows a relationship between $\E$ and $\Eb$.

\begin{lemma} \label{conjugate of edge-tangle}
Let $\theta$ be an integer with $\theta \geq 2$ and $G$ a graph. 
If $\E$ is an edge-tangle of order $3\theta-2$ of $G$, then $\Eb$ is a tangle of order $\theta$ in $L(G)$.
\end{lemma}

\begin{pf}
Observe that every member of $\Eb$ has order less than $\lceil \frac{3\theta-2}{3} \rceil = \theta$.
We shall prove that $\Eb$ satisfies tangle axioms (T1), (T2) and (T3).
Note that $\lvert E(G) \rvert \geq 3\theta-2$ since $G$ has an edge-tangle of order $3\theta-2$.

\noindent{\bf Claim 1:} $\Eb$ satisfies (T1).

\noindent{\bf Proof of Claim 1:}
Let $(A,B)$ be a separation of $L(G)$ of order less than $\theta$.
Let $(A_1,B_1)$ and $(B_2,A_2)$ be the normalizations of $(A,B)$ and $(B,A)$, respectively.
And let $[A'_1,B'_1]$ and $[B_2',A_2']$ be the partners of $(A_1,B_1)$ and $(B_2,A_2)$, respectively.
If any of $[A'_1,B'_1]$ and $[B'_2,A'_2]$ is in $\E$, then $(A,B)$ or $(B,A)$ is in $\Eb$, and we are done.
So we may assume that none of $[A'_1,B'_1]$ and $[B'_2,A'_2]$ is in $\E$.
By Lemmas \ref{normalization basic} and \ref{partner of sep}, the order of $[A'_1,B'_1]$ and $[B'_2,A'_2]$ are less than $\theta$, so $[B'_1,A'_1]$ and $[A'_2,B'_2]$ are in $\E$ by (E1).
Let $D$ be the union of the vertex-sets of the components of $G$ with at most one edge.
Let $[B_1'',A_1'']=[B_1' \cup D, A_1' -D]$ and let $[A_2'', B_2''] = [A_2' \cup D, B_2'-D]$.
By Statement 2 of Lemma \ref{moving isolated edges}, $[B_1'',A_1'']$ and $[A_2'',B_2'']$ are in $\E$.

Let $v \in A''_1 \cap B''_2$.
So $v \not \in D$ and hence $v$ does not belong to any component of $G$ with at most one edge.
In particular, $v$ is not an isolated vertex in $G$ and $\cl(v) \neq \emptyset$.
Note that $\cl(v)$ is not a set consisting of one isolated vertex in $L(G)$, for otherwise the vertex in $\cl(v)$ is the only edge of some component of $G$ and $v \in D$.
Since $v \in A_1'' \subseteq A_1'$, $\cl(v) \subseteq V(A_1)$ as $v$ is not an isolated vertex in $G$.
Similarly, $\cl(v) \subseteq V(B_2)$.
So $\cl(v) \subseteq V(A_1) \cap V(B_2)$.

Suppose that $\cl(v) \not \subseteq V(A)$.
Then some vertex in $\cl(v)$ is in $V(B)-V(A)$ and is not an isolated vertex of $B$, so this vertex is in $V(B_1)-V(A_1)$.
Hence $\cl(v)-V(A_1) \neq \emptyset$, a contradiction.

So $\cl(v) \subseteq V(A)$.
Similarly, $\cl(v) \subseteq V(B)$, for otherwise $\cl(v)-V(B_2) \neq \emptyset$.
Hence $\cl(v) \subseteq V(A) \cap V(B)$.

This shows that $\bigcup_{u \in A''_1 \cap B''_2} \cl(u) \subseteq V(A) \cap V(B)$.
Therefore, the number of edges of $G$ incident with some vertex in $A_1'' \cap B_2''$ is at most $\lvert \bigcup_{u \in A''_1 \cap B''_2} \cl(u) \rvert \leq \lvert V(A) \cap V(B) \rvert < \theta$.
In particular, the number of edges with one end in $A_1'' \cap B_2''$ and one end in $B_1'' \cup A_2''$ is less than $\theta$.
Hence, (E1) implies that either $[A_1'' \cap B_2'', B_1'' \cup A_2''] \in \E$ or $[B_1'' \cup A_2'', A_1'' \cap B_2''] \in \E$.
But (E3) excludes the latter case, so $[A_1'' \cap B_2'', B_1'' \cup A_2''] \in \E$.
However, $[B_1'',A_1'']$, $[A_2'',B_2'']$ and $[A_1'' \cap B_2'',B_1''\cup A_2'']$ belong to $\E$, but $A_1'' \cap B_2'' \cap (B_1'' \cup A_2'') = \emptyset$, contradicting (E2).
This proves that $\Eb$ satisfies (T1).
$\Box$

\smallskip

Next, we show that $\Eb$ satisfies (T3).
Suppose that $(A,B) \in \Eb$ with $V(A)=V(L(G))$.
So the partner $[A',B']$ of the normalization of $(A,B)$ is in $\E$.
Note that $\bigcup_{v \in B'} \cl(v) \subseteq V(B) = V(A) \cap V(B)$, so the number of edges incident with vertices in $B'$ is at most $\lvert V(A) \cap V(B) \rvert < \theta \leq 3\theta-2$, contradicting (E3).
Consequently, $\Eb$ satisfies (T3).

Now we suppose that $\Eb$ does not satisfy (T2).
So there exist separations $(A_1,B_1)$, $(A_2,B_2)$, $(A_3,B_3)$ in $\Eb$ such that $A_1 \cup A_2 \cup A_3 = L(G)$.
For each $i \in [3]$, let $(A_i^*,B_i^*)$ be the normalization of $(A_i, B_i)$, and let $[A_i',B_i']$ be the partner of $(A_i^*,B_i^*)$.
By the definition of $\Eb$, $[A_i',B_i'] \in \E$ for $i \in [3]$.
The number of edges of $G$ incident with vertices in $\bigcap_{i=1}^3 B_i'$ is at most $\lvert \bigcup_{v \in \bigcap_{i=1}^3 B_i'} \cl(v) \rvert \leq \lvert \bigcap_{i=1}^3 V(B_i^*) \rvert \leq \lvert \bigcap_{i=1}^3 V(B_i) \rvert$.
However, $\bigcap_{i=1}^3 V(B_i) \subseteq \bigcup_{i=1}^3 V(A_i \cap B_i)$, as $A_1 \cup A_2 \cup A_3 = L(G)$.
So the number of edges of $G$ incident with vertices in $\bigcap_{i=1}^3 B_i'$ is at most $\lvert \bigcup_{i=1}^3 V(A_i \cap B_i)\rvert \leq 3(\theta-1)$.
In addition, $\lvert [A_1' \cup A_2', B_1' \cap B_2'] \rvert \leq \sum_{i=1}^2 \lvert [A_i',B_i'] \rvert \leq 2(\theta-1)<3\theta-2$, so $[A_1' \cup A_2', B_1' \cap B_2'] \in \E$ by Lemma \ref{easy edge-tangle}.
Similarly, $\lvert [A_1' \cup A_2' \cup A_3', B_1' \cap B_2' \cap B_3'] \rvert \leq \sum_{i=1}^3 \lvert [A_i',B_i'] \rvert \leq 3(\theta-1)<3\theta-2$, so $[A_1' \cup A_2' \cup A_3', B_1' \cap B_2' \cap B_3'] \in \E$ by Lemma \ref{easy edge-tangle}.
Hence by (E3), the number of edges of $G$ incident with vertices in $\bigcap_{i=1}^3 B_i'$ is at least $3\theta-2$, a contradiction.
This proves that $\Eb$ satisfies (T3).
Consequently, $\Eb$ is a tangle of order $\theta$ in $L(G)$.
\end{pf}

\bigskip

Let $G$ be a graph and $\E$ a collection of edge-cuts of $G$ of order less than a positive number $\theta$, and let $X \subseteq E(G)$.
Define $\E-X$ to be the set of edge-cuts of $G-X$ of order less than $\theta-\lvert X \rvert$ such that $[A,B] \in \E-X$ if and only if $[A,B] \in \E$.

\begin{lemma} \label{edge-tangle deleting edges}
Let $G$ be a graph and $\theta$ be a positive integer.
If $\E$ is an edge-tangle in $G$ of order $\theta$ and $X$ is a subset of $E(G)$ with $\lvert X \rvert < \theta$, then $\E-X$ is an edge-tangle in $G-X$ of order $\theta-\lvert X \rvert$.
\end{lemma}

\begin{pf}
If $[A,B]$ is an edge-cut of order less than $\theta-\lvert X \rvert$ in $G-X$, then $[A,B]$ is an edge-cut in $G$ of order less than $\theta$.
So for every edge-cut $[A,B]$ of $G-X$ of order less than $\theta-\lvert X \rvert$, since $\E$ is an edge-tangle in $G$ of order $\theta$, either $[A,B] \in \E$ or $[B,A] \in \E$, and hence either $[A,B] \in \E-X$ or $[B,A] \in \E-X$ by the definition of $\E-X$.
This shows that $\E-X$ satisfies (E1).

Since $\E$ satisfies (E2), $B_1 \cap B_2 \cap B_3 \neq \emptyset$, for any edge-cuts $[A_1,B_1],[A_2,B_2],[A_3,B_3] \in \E$.
So for any members $[A_1,B_1],[A_2,B_2],[A_3,B_3]$ of $\E-X$, we have $[A_1,B_1],[A_2,B_2],[A_3.B_3] \in \E$ by the definition of $\E-X$, so $B_1 \cap B_2 \cap B_3 \neq \emptyset$.
Hence $\E-X$ satisfies (E2).

If $[A,B] \in \E-X$, then $[A,B] \in \E$, so $G$ contains at least $\theta$ edges incident with vertices in $B$.
Hence, $G-X$ contains at least $\theta-\lvert X \rvert$ edges incident with vertices in $B$.
So $\E-X$ satisfies (E3).
This proves that $\E-X$ is an edge-tangle of order $\theta-\lvert X \rvert$.
\end{pf}

\bigskip

Let $\E$ be an edge-tangle in a graph $G$.
We say that a subset $Y$ of $E(G)$ is {\it free} with respect to $\E$ if there exist no $Z \subseteq Y$ and $[A,B] \in \E-Z$ of order less than $\lvert Y-Z \rvert$ such that every edge in $Y-Z$ has every end in $A$.

\begin{lemma} \label{subset of free set is free}
Let $G$ be a graph and $\E$ an edge-tangle in $G$ of order $\theta \geq 1$.
Let $Z$ be a subset of $E(G)$ with $\lvert Z \rvert < \theta$ and let $X$ be a subset of $E(G)-Z$ such that $X$ is free with respect to $\E-Z$.
If $\lvert X \rvert \leq \theta-\lvert Z \rvert$, then for every $X' \subseteq X$ and $Z' \subseteq Z$, $X'$ is free with respect to $\E-Z'$.
\end{lemma}

\begin{pf}
Suppose to the contrary that $X'$ is not free with respect to $\E-Z'$.
Then there exist $W' \subseteq X'$ and $[A,B] \in (\E-Z')-W'=\E-(Z' \cup W')$ of order less than $\lvert X'-W' \rvert$ such that every edge in $X'-W'$ has every end in $A$.
Let $W=W' \cup (X-X')$.
So $X'-W' = X-W$.
Since $Z' \subseteq Z$ and $W' \subseteq X'$, $[A,B]$ is an edge-cut of $G-Z$ of order less than $\lvert X'-W' \rvert + \lvert W' \rvert \leq \lvert X' \rvert \leq \lvert X \rvert \leq \theta-\lvert Z \rvert$.
Since $\E-Z$ has order $\theta-\lvert Z \rvert$ by Lemma \ref{edge-tangle deleting edges}, by (E1), either $[A,B] \in \E-Z$ or $[B,A] \in \E-Z$.

Since $W \subseteq X$ and $X$ is free with respect to $\E-Z$, either $[A,B] \not \in (\E-Z)-W$ or the order of $[A,B]$ in $G-(Z \cup W)$ is at least $\lvert X-W \rvert = \lvert X'-W' \rvert$. 
Since $Z \cup W \supseteq Z' \cup W'$, the order of $[A,B]$ in $G-(Z \cup W)$ is at most the order of $[A,B]$ in $G-(Z' \cup W')$, which is less than $\lvert X'-W' \rvert$.
So $[A,B] \not \in (\E-Z)-W$.
Hence $[A,B] \not \in \E-Z$ by the definition of $(\E-Z)-W$.
So $[B,A] \in \E-Z$.
By the definition of $\E-Z$, $[B,A] \in \E$.
Since $[A,B] \in \E-(Z' \cup W')$, the order of $[B,A]$ in $G-(Z' \cup W')$ is less than $\theta-\lvert Z' \cup W' \rvert$.
Since $[B,A] \in \E$, $[B,A] \in \E-(Z' \cup W')$ by the definition of $\E-(Z' \cup W')$.
So $[A,B],[B,A] \in \E-(Z' \cup W')$, contradicting (E2).
\end{pf}

\bigskip

Let $\T$ be a tangle in a graph $G$.
We say that a subset $X$ of $V(G)$ is {\it free} with respect to $\T$ if there does not exist $(A,B) \in \T$ of order less than $\lvert X \rvert$ such that $X \subseteq V(A)$.

Note that for every graph $G$, $V(L(G))=E(G)$.
So for every subset of $E(G)$, we can consider whether it is free with respect to an edge-tangle $\E$ in $G$ and whether it is free with respect to the conjugate $\Eb$ of $\E$.

\begin{lemma} \label{free wrt edge free wrt vertex}
Let $\E$ be an edge-tangle in a graph $G$, and let $\Eb$ be the conjugate of $\E$.
Let $X,Z$ be disjoint subsets of $E(G)$.
If $X$ is free with respect to $\E-Z$, then $X$ is free with respect to $\Eb-Z$.
\end{lemma}

\begin{pf}
Suppose that $X$ is not free with respect to $\Eb-Z$.
Then there exists a separation $(A,B) \in \Eb-Z$ of $L(G)-Z$ of order less than $\lvert X \rvert$ such that $X \subseteq V(A)$.
We may assume that the order of $(A,B)$ is as small as possible, and subject to that, $V(B)$ is inclusion-wise minimal.
So every vertex in $V(A) \cap V(B) - X$ is adjacent to a vertex in $V(A)-V(B)$ and adjacent to a vertex in $V(B)-V(A)$; every vertex in $V(A) \cap V(B) \cap X$ is adjacent to a vertex in $V(B)-V(A)$.
Furthermore, $B$ has no isolated vertices by the minimality of $(A,B)$.
Let $(A',B')$ be a separation of $L(G)$ with $V(A') = V(A) \cup Z$ and $V(B')=V(B) \cup Z$.
Since $(A,B) \in \Eb-Z$, $(A',B') \in \Eb$.
Let $(A^*,B^*)$ be the normalization of $(A',B')$.
So $V(B^*)-Z=V(B')-Z$, $V(A^*) \cap V(B^*)-(X \cup Z) = V(A') \cap V(B')-(X \cup Z)$ and $\lvert V(A^*) \cap V(B^*)-(X \cup Z) \rvert = \lvert V(A') \cap V(B') \rvert - \lvert V(A') \cap V(B') \cap (X \cup Z) \rvert$.

Suppose that $(A^*,B^*) \not \in \Eb$.
Then $(B^*,A^*) \in \Eb$ by (T1).
By (T1) and (T3), $(G[V(A') \cap V(B')], G-E(G[V(A') \cap V(B')])) \in \Eb$.
Since $V(B')-Z=V(B^*)-Z$, we know $(A',B'), (B^*,A^*), \allowbreak (G[V(A') \cap V(B')], G-E(G[V(A') \cap V(B')]))$ are members of $\Eb$ such that $A' \cup B^* \cup G[V(A') \cap V(B')] = L(G)$, contradicting (T2).
So $(A^*,B^*) \in \Eb$.

Let $[C,D]$ be the partner of $(A^*,B^*)$.
So $[C,D] \in \E$.
Let $W = (V(A') \cap V(B')) \cap (X \cup Z)$.
Note that $W$ is a subset of $V(L(G))=E(G)$.
Every edge $e$ in $X-W$ has every end in $C$ since it is a vertex in $V(A')-V(B')$.
And the order of $[C,D]$ in $G-W$ equals $\lvert V(A^*) \cap V(B^*)-(X \cup Z) \rvert  = \lvert V(A') \cap V(B') \rvert - \lvert W \rvert < \lvert X \rvert - \lvert W-Z \rvert \leq \lvert X-(W \cap X) \rvert$.

Let $\theta$ be the order of $\E$.
So the order of $\Eb$ is at most $\lceil \theta/3 \rceil$.
Since $(A',B') \in \Eb$, $\lvert V(A') \cap V(B') \rvert < \theta/3$.
So the order of $[C,D]$ in $G-W$ is at most $\lvert V(A') \cap V(B') \rvert - \lvert W \rvert < \theta/3-\lvert W \rvert <\theta-\lvert W \rvert$.
Hence $[C,D] \in \E-W$.

Therefore, $[C,D] \in (\E-Z)-(W \cap X)$ is an edge-cut of $(G-Z)-(W \cap X)$ of order less than $\lvert X - (X \cap W) \rvert$ such that every edge in $X-(X \cap W)$ has every end in $C$.
So $X$ is not free with respect to $\E-Z$, a contradiction.
\end{pf}

\bigskip

The converse of Lemma \ref{free wrt edge free wrt vertex} is also true when $Z=\emptyset$, subject to a requirement on the size of the set, as shown in the following lemma.

\begin{lemma} \label{free wrt vertex free wrt edge}
Let $\E$ be an edge-tangle in a graph $G$ of order at least two, and let $\Eb$ be the conjugate of $\E$.
Let $X$ be a subset of $E(G)$.
Denote the order of $\Eb$ by $\theta$.
If $X$ is free with respect to $\Eb$ and $\lvert X \rvert \leq \theta$, then $X$ is free with respect to $\E$.
\end{lemma}

\begin{pf}
Suppose to the contrary that $X$ is not free with respect to $\E$.
So there exist $W \subseteq X$ and $[A,B] \in \E-W$ of order less than $\lvert X - W \rvert$ such that every edge in $X-W$ has every end in $A$.
We further assume that the order of $[A,B]$ is as small as possible, and subject to that, $\lvert A \rvert$ is as large as possible.

\noindent{\bf Claim 1:} The following hold.
	\begin{itemize}
		\item Every non-isolated vertex of $G$ is incident with an edge of $G-W$ whose every end is contained in $A$ or an edge of $G-W$ whose every end is contained in $B$.
		\item $A$ contains the vertex-set of every component of $G$ with at most one edge.
	\end{itemize}

\noindent{\bf Proof of Claim 1:}
The first statement of this claim immediately follows from the minimality of the order of $[A,B]$.
Furthermore, the minimality of $\lvert [A,B] \rvert$, the maximality of $A$, and Statement 2 of Lemma \ref{moving isolated edges} imply the second statement of this claim.
$\Box$

\smallskip

Define $B'$ to be the induced subgraph of $L(G-W)$ such that $V(B')$ is the set of edges of $G-W$ incident with vertices of $B$.
Define $(A',B')$ to be a separation of $L(G-W)$ such that $V(A') \cap V(B')$ is the set of edges of $G-W$ with one end in $A$ and one end in $B$. 
Note that the order of $(A',B')$ is at most the order of $[A,B]$ in $G-W$.
So $\lvert V(A') \cap V(B') \rvert < \lvert X-W \rvert$.

Let $(A^*,B^*)$ be a separation of $L(G)$ with $V(A^*)=V(A') \cup W$ and $V(B^*)=V(B') \cup W$.
So $X \subseteq V(A^*)$ and $\lvert V(A^*) \cap V(B^*) \rvert < \lvert X \rvert \leq \theta$.
By (T1), $(A^*,B^*) \in \Eb$ or $(B^*,A^*) \in \Eb$.

Let $(A'',B'')$ be the normalization of $(A^*,B^*)$.
Let $[C,D]$ be the partner of $(A'',B'')$ in $G$.

\noindent{\bf Claim 2:} $[C,D]=[A,B]$.

\noindent{\bf Proof of Claim 2:}
Suppose that there exists a vertex $v \in C-A$.
Since $A$ contains all isolated vertices in $G$ by Claim 1, $v$ is not an isolated vertex in $G$.
Since $v \in C$, $\cl(v) \subseteq V(A'')$.
Since $v \not \in A$, $v \in B$.
By Claim 1, $v$ is incident with an edge $e$ of $G-W$ whose every end is contained in $B$.
So $e \in V(B')-V(A') = V(B^*)-V(A^*)$.
Since $v \in B$, $e$ is not the only edge of some component of $G$ by Claim 1.
So $e$ is not an isolated vertex in $B^*$.
Hence $e \in V(B'')-V(A'')$.
But $e \in \cl(v)-V(A'')$, a contradiction.
This shows $C \subseteq A$.

Suppose that there exists a vertex $u \in A-C$.
Since $C$ contains all isolated vertices of $G$, $u$ is not an isolated vertex of $G$.
Since $u \in A$, $u$ is incident with an edge $f$ of $G-W$ whose every end in $A$ by Claim 1.
Hence $f \in V(A')-V(B')=V(A^*)-V(B^*)$.
So $f \in V(A'')-V(B'')$ and hence $\cl(u) \subseteq V(A'')$.
This implies that $u \in C$, a contradiction.

Therefore, $A=C$.
Since $\{A,B\}$ and $\{C,D\}$ are partitions of $V(G)$, $[A,B]=[C,D]$.
$\Box$

\smallskip

By Claim 2, since $[C,D] = [A,B] \in \E$, $(A^*,B^*) \in \Eb$.
But $X \subseteq V(A^*)$ and $\lvert V(A^*) \cap V(B^*) \rvert < \lvert X \rvert$, so $X$ is not free with respect to $\Eb$, a contradiction.
This proves that $X$ is free with respect to $\E$.
\end{pf}

\subsection{Immersions and edge-tangles} \label{subsec:imm_edge-tangles}

The following lemma provides a way to obtain an edge-tangle from an immersion.

\begin{lemma} \label{immersion induces edge-tangle}
Let $H$ be a graph and $\E'$ an edge-tangle of order $\theta$ in $H$.
Let $G$ be a graph that contains an $H$-immersion $(\pi_V, \pi_E)$.
If $\E$ is the set of all edge-cuts $[A,B]$ of $G$ of order less than $\theta$ such that there exists $[A',B'] \in \E'$ with $\pi_V(A')=A \cap \pi_V(V(H))$, then $\E$ is an edge-tangle of order $\theta$ in $G$.
\end{lemma}

\begin{pf}
We shall show that $\E$ satisfies the edge-tangle axioms (E1), (E2) and (E3).
Note that for every edge-cut $[A,B]$ of $G$, $A \cap \pi_V(V(H))$ and $B \cap \pi_V(V(H))$ are two disjoint subsets of $\pi_V(V(H))$ such that their union is $\pi_V(V(H))$, so there exists an edge-cut $[A',B']$ of $H$ such that $\pi_V(A') = A \cap \pi_V(V(H))$ and $\pi_V(B') = B \cap \pi_V(V(H))$.

We first prove that $\E$ satisfies (E1).
Let $[A,B]$ be an edge-cut of $G$ of order less than $\theta$.
Let $[A',B']$ an edge-cut of $H$ such that $\pi_V(A') = A \cap \pi_V(V(H))$ and $\pi_V(B') = B \cap \pi_V(V(H))$.
Since $(\pi_V,\pi_E)$ is an $H$-immersion, there are at least $\lvert [A', B'] \rvert$ edge-disjoint paths in $G$ from $A \cap \pi_V(V(H))$ to $B \cap \pi_V(V(H))$.
So $\lvert [A',B'] \rvert \leq \lvert [A,B] \rvert < \theta$.
Hence, one of $[A',B']$ and $[B',A']$ is in $\E'$, and hence one of $[A,B]$ and $[B,A]$ is in $\E$.
So $\E$ satisfies (E1).

Now we prove that $\E$ satisfies (E2).
For each $i \in [3]$, let $[A_i,B_i] \in \E$ be an edge-cut of $G$.
By the definition of $\E$, for each $i \in [3]$, there exists $[A'_i, B'_i] \in \E'$ such that $\pi_V(A_i')=A_i \cap \pi_V(V(H))$ and hence $\pi_V(B'_i) = B_i \cap \pi_V(V(H))$.
Since $\E'$ is an edge-tangle in $H$, $B_1' \cap B_2' \cap B_3'$ contains a vertex $v$ of $H$.
So $\pi_V(v) \in B_1 \cap B_2 \cap B_3$.
This proves that $\E$ satisfies (E2).

Finally, we prove that $\E$ satisfies (E3).
Let $[A,B] \in \E$.
By the definition of $\E$, there exists $[A',B'] \in \E'$ such that $\pi_V(A') = A \cap \pi_V(V(H))$ and hence $\pi_V(B') = B \cap \pi_V(V(H))$.
Since $\E'$ satisfies (E3), $H$ contains at least $\theta$ edges incident with vertices in $B'$.
So there are at least $\theta$ edge-disjoint subgraphs of $G$ each containing a vertex in $\pi_V(B') \subseteq B$ and containing an edge of $G$.
Therefore, $G$ contains at least $\theta$ edges incident with vertices in $B$.
Consequently, $\E$ is an edge-tangle in $G$.
\end{pf}

\bigskip

We call the edge-tangle $\E$ defined in Lemma \ref{immersion induces edge-tangle} the {\it edge-tangle induced by the $H$-immersion $(\pi_V,\pi_E)$ and the edge-tangle $\E'$ in $H$}.

The {\it $m \times n$ wall} is the simple graph with vertex-set $\{(i,j): 1 \leq i \leq n, 1 \leq j \leq m\}$ and edge-set $\{(i,j)(i+1,j): 1 \leq i \leq n-1, 1 \leq j \leq m\} \cup \{(2a-1,2b-1)(2a-1,2b): 1 \leq a \leq \lceil n/2 \rceil, 1 \leq b \leq \lfloor m/2 \rfloor\} \cup \{(2a,2b)(2a,2b+1): 1 \leq a \leq \lfloor n/2 \rfloor, 1 \leq b \leq \lfloor (m-1)/2 \rfloor\}$.
The {\it $i$-th row} of the $m \times n$ wall is the subgraph induced by $\{(x,i): 1 \leq x \leq n\}$.
The {\it $k$-th column} of the $m \times n$ wall is the subgraph induced by $\{(x,y): 2k-1 \leq x \leq \min\{2k,n\},1 \leq y \leq m\}$.
Hence, the $m \times n$ wall contains $m$ rows and $\lceil n/2 \rceil$ columns.

It was proved by Robertson, Seymour and Thomas \cite{rst} (see Theorem \ref{truncation tangle} below) that every tangle is ``equivalent'' to a subdivision of a wall.
The next objective of this section is to prove an analogous result (Lemma \ref{truncation large deg or wall}) about edge-tangles and immersions, which will be used in Section \ref{sec:4-edge-connected edge-tangle}.

For a graph $H$, an $H$-immersion $(\pi_v, \pi_e)$ is an {\it $H$-subdivision} if
\begin{itemize}
	\item for every pair of distinct edges $e_1,e_2$ of $H$, $V(\pi_E(e_1) \cap \pi_E(e_2)) \subseteq \pi_V(S)$, where $S$ is the set of the common ends of $e_1,e_2$, and 
	\item for every $e \in E(H)$, $\pi_V(V(H)) \cap V(\pi_E(e)) = \pi_V(S)$, where $S$ is the set of ends of $e$.
\end{itemize}
We say that a {\it tangle $\T$ is induced by a $r \times r$ wall-subdivision} $(\pi_V,\pi_E)$ if for every $(A,B) \in \T$, there exists a row of the wall such that $E(B)$ intersects $\pi_E(e)$ for every edge $e$ of this row.

One corollary of the following restatement of \cite[(2.3)]{rst} is that every graph with a tangle of large order contains a subdivision of a large wall.

\begin{theorem} \cite[(2.3)]{rst} \label{truncation tangle}
	Let $\theta \geq 2$, and let $\T$ be a tangle in $G$ of order at least $20^{\theta^4(2\theta-1)}$.
	If $\T' \subseteq \T$ is a tangle of order $\theta$, then $\T'$ is induced by a $\theta \times \theta$ wall-subdivision.
\end{theorem}

\begin{lemma} \label{small edge-cut wall}
Let $r$ and $\theta$ be positive integers with $\theta \leq r$.
Let $G$ be the $r \times 2r$ wall.
If $[A,B]$ is an edge-cut of order less than $\theta$ of $G$, then 
	\begin{enumerate}
		\item exactly one of $A$ and $B$ contains all vertices of a column of $G$, 
		\item exactly one of $A$ and $B$ contains all vertices of a row of $G$,
		\item $A$ contains all vertices of a column if and only if $A$ contains all vertices of a row, 
		\item exactly one of $A$ and $B$ intersects vertices in at least $\theta$ columns of $G$,
		\item exactly one of $A$ and $B$ intersects vertices in at least $\theta$ rows of $G$,  
		\item $A$ intersects vertices in at least $\theta$ columns of $G$ if and only if $A$ contains all vertices of a column of $G$, and 
		\item $A$ intersects vertices in at least $\theta$ rows of $G$ if and only if $A$ contains all vertices of a column of $G$.
	\end{enumerate}
\end{lemma}

\begin{pf}
Note that $G$ has $r$ rows and $r$ columns.
Suppose that $A$ contains all vertices of a column and $B$ contains all vertices of another column.
Then every row must contain an edge in $[A,B]$, so the order of $[A,B]$ is at least $r$, a contradiction.
Suppose that none of $A$ and $B$ contains all vertices of a column.
Then every column must contain an edge in $[A,B]$, so the order of $[A,B]$ is at least $r$, a contradiction.
So exactly one of $A$ and $B$ contains all vertices of a column.
Similarly, exactly one of $A$ and $B$ contains all vertices of a row.
Furthermore, if $A$ contains all vertices of a column, then $B$ cannot contain all vertices of a row, so $A$ contains all vertices of a row as well.
Similarly, if $B$ contains all vertices of a column, then $B$ contains all vertices of a row.
This shows Statements 1-3.

Since $r \geq \theta$, every column intersects vertices in at least $\theta$ rows.
By Statement 1, at least one of $A$ and $B$ intersects vertices in at least $\theta$ rows.
If one of $A$ and $B$ contains all vertices of a column and the other intersects vertices in at least $\theta$ rows, then there are at least $\theta$ rows containing both vertices in $A$ and in $B$, so there are at least $\theta$ edges between $A$ and $B$, a contradiction.
So $A$ contains all vertices of a column if and only if $A$ intersects vertices in at least $\theta$ rows.
This shows Statement 7.
Then Statements 1 and 7 imply Statement 5.

Similarly, $A$ contains all vertices of a row if and only if $A$ intersects vertices in at least $\theta$ columns.
So Statement 6 follows from Statement 3, and Statement 4 follows from Statements 5-7.
\end{pf}

\begin{lemma} \label{wall has edge-tangle}
Let $r$ and $\theta$ be positive integers.
Let $G$ be the $r \times 2r$ wall.
Let $\E$ be the set of all edge-cuts $[A,B]$ of order less than $\theta$ of $G$ satisfying that $B$ intersects vertices in at least $\theta$ columns of $G$.
If $r \geq 2\theta$, then $\E$ is an edge-tangle of $G$ of order $\theta$.
\end{lemma}

\begin{pf}
Let $[A,B]$ be an edge-cut of $G$ of order less than $\theta$.
By Lemma \ref{small edge-cut wall}, $[A,B]$ or $[B,A]$ is in $\E$.
Hence, $\E$ satisfies (E1).

Let $[A_i,B_i] \in \E$ be edge-cuts of $G$ for $i \in [3]$.
For each $i \in [3]$, $B_i$ intersects vertices in at least $\theta$ columns of $G$, so $B_i$ contains all vertices of a column of $G$ by Lemma \ref{small edge-cut wall}. 
For each $i \in [3]$, let $c_i$ be a column of $G$ contained in $B_i$.
Suppose that $B_1 \cap B_2 \cap B_3 = \emptyset$.
If $c_1=c_2$, then $B_3$ is disjoint from $c_1$, so $A_3$ contains $c_1$, a contradiction.
So by symmetry, we may assume that $c_1,c_2,c_3$ are pairwise distinct.
Since $c_1$ is contained in $A_2 \cup A_3$, we know that $A_2$ or $A_3$, say $A_2$, contains at least one half vertices of $c_1$, so $A_2$ intersects at least $r/2$ rows.
But $B_2$ contains $c_2$, so there are at least $r/2$ edges with one end in $A_2$ and one end in $B_2$.
Therefore, $[A_2,B_2]$ has order at least $r/2 \geq \theta$, a contradiction.
Hence, $\E$ satisfies (E2).

For every $[A,B] \in \E$, $B$ contains a column in $G$ by Lemma \ref{small edge-cut wall}, so there are at least $r \geq \theta$ edges in $G$ incident with some vertices in $B$.
This proves that $\E$ is an edge-tangle.
\end{pf}

\bigskip

We call the edge-tangle $\E$ mentioned in Lemma \ref{wall has edge-tangle} the {\it natural edge-tangle} in the $r \times 2r$ wall of order $\theta$.
Here is a short summary about natural edge-tangles in a wall.

\begin{lemma} \label{summary edge-tanlge induced by wall immersion}
Let $r$ and $\theta$ be positive integers with $r \geq 2\theta$.
Let $W$ be the $r \times 2r$ wall.
Let $G$ be a graph and $\Pi=(\pi_V,\pi_E)$ a $W$-immersion in $G$.
If $\E$ is the edge-tangle in $G$ induced by $\Pi$ and the natural edge-tangle of order $\theta$ in $W$, then for every edge-cut $[A,B]$ of $G$ of order less than $\theta$, the following are equivalent.
	\begin{enumerate}
		\item $[A,B] \in \E$.
		\item $B$ intersects the image of $\pi_V$ of vertices in at least $\theta$ columns of $W$.
		\item $B$ contains the image of $\pi_V$ of all vertices of some column of $W$.
		\item $B$ contains the image of $\pi_V$ of all vertices of some row of $W$.
	\end{enumerate}
\end{lemma}

\begin{pf}
Statements 1 and 2 are equivalent by the definition of $\E$.
Statements 2-4 are equivalent by Lemma \ref{small edge-cut wall}.
\end{pf}

\bigskip

We need the following lemma to prove Lemma \ref{truncation large deg or wall}. 
It states that whenever each ``vertex'' of a large ``grid'' is labelled with a bounded number of labels in a way that every label is used by a bounded number of times, one can find a large ``subgrid'' such that the sets of labels used in this ``subgrid'' are pairwise disjoint.

\begin{lemma} \label{grid pigeonhole}
	For any nonnegative integers $s,t,p,q$, there exist integers $s^*=s^*(s,t,p,q),t^*=t^*(s,t,p,q)$ such that the following holds.
	Let $I,J$ be sets with $\lvert I \rvert=t^*$ and $\lvert J \rvert=s^*$.
	Let $U$ be a set, and let $f$ be a function that maps each pair $(i,j) \in I \times J$ to a subset of $U$ of size at most $p$.
	If for every $u \in U$, $\lvert \{(i,j) \in I \times J: u \in f((i,j))\} \rvert \leq q$, then there exist $I' \subseteq I$ with $\lvert I' \rvert=t$ and $J' \subseteq J$ with $\lvert J' \rvert = s$ such that $f((x,y)) \cap f((x',y')) =\emptyset$ for distinct $(x,y),(x',y') \in I' \times J'$.
\end{lemma}

\begin{pf}
	We shall prove this lemma by induction on $s$.
	When $s=0$, the lemma holds obviously.
	So we may assume that $s \geq 1$ and this lemma holds for every smaller $s$.
	
	Let $s_1=s^*(s-1,t,p,q)$ and $t_1=t^*(s-1,t,p,q)$.
	Define $s^*(s,t,p,q)=s_1+pqt_1+1$ and $t^*(s,t,p,q) = pqt_1$.

	Without loss of generality, we may assume that $I=[t^*]$ and $J=[s^*]$.
	Note that for every $x \in [t^*]$, there are at most $pq$ elements $x' \in [t^*]$ such that $f((x,1)) \cap f((x',1)) \neq \emptyset$.
	So there exists a subset $I_1$ of $I$ with $\lvert I_1 \rvert = \frac{\lvert I \rvert}{pq} = t_1$ such that $f((x,1)) \cap f((x',1)) = \emptyset$ for distinct $x,x' \in I_1$.
	Note that $\lvert \bigcup_{x \in I_1}f((x,1)) \rvert \leq pt_1$.
	So there exists a subset $J_1$ of $J-\{1\}$ with $\lvert J_1 \rvert \geq \lvert J \rvert-1 - pt_1q = s_1$ such that $f((x,y)) \cap f((i,1)) = \emptyset$ for every $(x,y) \in I_1 \times J_1$ and $i \in I_1$.
	By the induction hypothesis, there exist $I_2 \subseteq I_1$ with $\lvert I_2 \rvert = t$ and $J_2 \subseteq J_1$ with $\lvert J_2 \rvert =s-1$ such that $f((x,y)) \cap f((x',y')) = \emptyset$ for distinct $(x,y),(x',y') \in I_2 \times J_2$.
	Define $I'=I_2$ and $J'=\{1\} \cup J_2$.
	Then $f((x,y)) \cap f((x',y')) = \emptyset$ for distinct $(x,y),(x',y') \in I' \times J'$.
\end{pf}

\begin{lemma} \label{truncation large deg or wall}
For every positive integers $\theta$ and $d$ with $\theta \geq 2$, there exists an integer $w=w(\theta,d)$ such that if $\E$ is an edge-tangle in a graph $G$ of order at least $w$, then either there exists $v \in V(G)$ incident with at least $d$ edges in $G$ such that $v \in B$ for every $[A,B] \in \E_\theta$, or $\E_\theta$ is induced by a $2\theta \times 4\theta$ wall-immersion and the natural edge-tangle of order $\theta$ in the $2\theta \times 4\theta$ wall, where $\E_\theta$ is the edge-tangle in $G$ of order $\theta$ such that $\E_\theta \subseteq \E$.
\end{lemma}

Now we sketch the proof of Lemma \ref{truncation large deg or wall}.
Since $\E$ is an edge-tangle in $G$ of large order, $\Eb$ is a tangle in $L(G)$ of large order, so it is induced by a very large wall subdivision $\Pi$ in $L(G)$ by Theorem \ref{truncation tangle}.
If there exists a vertex $v$ of $G$ such that $\cl(v)$ contains many branch vertices of $\Pi$, then it is not hard to show that this vertex $v$ satisfies the conclusion of Lemma \ref{truncation large deg or wall}.
If there exists no such vertex $v$ exists, then there exists a smaller (but still sufficiently large) wall subdivision $\Pi^*$ such that every branch vertex of $\Pi^*$ is a branch vertex of $\Pi$, and $\cl(u)$ contains at most one branch vertex of $\Pi^*$ for every $u \in V(G)$.
Such a wall subdivision $\Pi^*$ in $L(G)$ defines a wall immersion $\Pi'$ in $G$ in an obvious way.
Then one can show that the wall immersion $\Pi'$ satisfies the conclusion of Lemma \ref{truncation large deg or wall} by using the relationship between $\E$ and $\Eb$.

\bigskip

\begin{pf1}{Lemma \ref{truncation large deg or wall}}
Let $\theta$ and $d$ be positive integers with $\theta \geq 2$.
	Let $\theta' = s_{\ref{grid pigeonhole}}(2\theta,4\theta,4, \allowbreak (2\theta d)^2) + t_{\ref{grid pigeonhole}}(2\theta,4\theta,4,(2\theta d)^2) + 5\theta$, where $s_{\ref{grid pigeonhole}}$ and $t_{\ref{grid pigeonhole}}$ are the integers $s^*$ and $t^*$ mentioned in Lemma \ref{grid pigeonhole}. 
	Define $w=20^{64\theta'^5}$.
	
	Denote the $2\theta' \times 2\theta'$ wall by $W$.
	Let $G$ be a graph, and let $\E$ be an edge-tangle in $G$ of order at least $w$.
	By Lemma \ref{conjugate of edge-tangle}, $\Eb$ is a tangle of order at least $w/3-1 \geq 20^{(2\theta')^4(4\theta'-1)}$ in $L(G)$.
	For every integer $t$, let $\Eb_t$ be the tangle in $L(G)$ of order $t$ with $\Eb_t \subseteq \Eb$.
	By Theorem \ref{truncation tangle}, $\Eb_{2\theta'}$ is induced by a $W$-subdivision $(\pi_V,\pi_E)$ in $L(G)$.
	
	\noindent{\bf Claim 1:} If there exists $v \in V(G)$ such that $\cl(v)$ contains at least $(2\theta d)^2$ vertices in $\pi_V(V(W))$, then $v$ is incident with at least $d$ edges in $G$, and $v \in B$ for every $[A,B] \in \E_\theta$.
	
	\noindent{\bf Proof of Claim 1:}
	Let $v$ be a vertex of $G$ such that $\cl(v)$ contains at least $(2\theta d)^2$ vertices in $\pi_V(V(W))$.
	So $\lvert \cl(v) \rvert \geq (2\theta d)^2$ and $v$ is incident with at least $(2\theta d)^2 \geq d$ edges in $G$.
	
	Suppose that there exists an edge-cut $[A,B] \in \E_\theta$ such that $v \in A$.
	We may assume that the order of $[A,B]$ is as small as possible, and subject to that, $A$ is maximal.
	By Lemma \ref{easy edge-tangle} and the maximality of $A$, $A$ contains all isolated vertices of $G$.
	
	Since the order of $[A,B]$ is less than $\theta$ and $v \in A$ is incident with at least $(2\theta d)^2 \geq \theta$ edges in $G$, some edge of $G$ incident with $v$ has every end in $A$.
	
	We define the following.
	\begin{itemize}
		\item[(i)] Define $B'$ to be a subgraph of $L(G)$ such that $V(B')$ consists of the edges of $G$ incident with vertices in $B$.
		\item[(ii)] Define $(A',B')$ to be a separation of $L(G)$ such that $V(A') \cap V(B')$ consists of the edges of $G$ with one end in $A$ and one end in $B$.
		\item[(iii)] Subject to (i) and (ii), $E(A')$ is maximal.
	\end{itemize}
	Since the order of $[A,B]$ is minimal and some edge of $G$ incident with $v$ has every end in $A$, $(A',B')$ is normalized.
	Hence $[A,B]$ is the partner of $(A',B')$.
	Therefore, $(A',B') \in \Eb$.

	Since $\lvert V(A') \cap V(B') \rvert < \theta  \leq 2\theta'$, $(A',B') \in \Eb_{2\theta'}$.
	Since $\Eb_{2\theta'}$ is induced by a $W$-subdivision $(\pi_V,\pi_E)$ in $L(G)$, there exists a row $r$ of $W$ such that $E(B')$ intersects $\pi_E(e)$ for every edge $e$ of the row $r$.
	Since $\lvert V(A') \cap V(B') \rvert < \theta$, there are at most $\theta$ vertices $x$ of $W$ in $r$ such that $\pi_V(x) \in V(A')-V(B')$.
	That is, there are at least $2\theta'-\theta$ vertices $x$ of $W$ in $r$ such that $\pi_V(x) \in V(B')$.
	
	Suppose that there exists a row $r'$ of $W$ other than $r$ such that there are at least $4\theta$ vertices $x$ of $r'$ such that $\pi_V(x) \in V(A')-V(B')$.
	Then there exist $2\theta$ columns $c_1,c_2,...,c_{2\theta}$ of $W$ such that $V(A')-V(B')$ intersects the image of $\pi_V$ of some vertices of each of $c_i$.
	Since there are at least $2\theta'-\theta$ vertices $x$ of $W$ in $r$ such that $\pi_V(x) \in V(B')$, there are at least $\theta$ columns $c$ in $\{c_i: 1 \leq i \leq 2\theta\}$ such that both $V(A')-V(B')$ and $V(B')$ intersect the image of $\pi_V$ of some vertices in $c$.
	Hence there are at least $\theta$ disjoint paths from $V(A')$ to $V(B')$, a contradiction.
	
	Therefore, for every row of $W$, there are at least $2\theta'-4\theta \geq 1$ vertices $x$ of this row such that $\pi_V(x) \in V(B')$.
	In particular, $V(B')$ intersects the image of $\pi_V$ of the vertices of each row.
	Since $\lvert V(A') \cap V(B') \rvert < \theta$, $V(A')$ intersects the image of $\pi_V$ of vertices in at most $\theta-1$ rows of $W$.
	Since $\theta' > \theta-1$, $V(B')$ intersects the image of $\pi_V$ of vertices in every column of $W$.
	
	Since $v \in A$, $V(A') \supseteq \cl(v)$ contains at least $(2\theta d)^2$ vertices in $\pi_V(V(W))$.
	So $V(A')$ contains the image of $\pi_V$ of some vertices in either at least $2\theta d$ rows of $W$ or at least $\theta d$ columns of $W$.
	Since $V(A')$ intersects the image of $\pi_V$ of vertices in at most $\theta-1$ rows of $W$, the former is impossible.
So $V(A')$ contains the image of $\pi_V$ of some vertices in at least $\theta d$ columns of $W$.
	But $V(B')$ intersects the image of $\pi_V$ of vertices in every column of $W$.
	So there exist $\theta d \geq \theta$ disjoint paths from $V(A')$ to $V(B')$, a contradiction.
	$\Box$

\smallskip

	By Claim 1, we may assume that for every $v \in V(G)$, $\cl(v)$ contains at most $(2\theta d)^2$ vertices in $\pi_V(V(W))$, for otherwise the lemma holds.
	
	Denote the $2\theta \times 4\theta$ wall by $W'$.
	
	\noindent{\bf Claim 2:} There exists a $W'$-subdivision $\Pi^* = (\pi_V^*,\pi_E^*)$ in $L(G)$ such that 
	\begin{itemize}
		\item $\pi_V^*(V(W')) \subseteq \pi_V(V(W))$ and 
		\begin{itemize}
			\item if $x \in V(W')$ and $x' \in V(W)$ with $\pi_V^*(x)=\pi_V(x')$, then $x' = (i,j)$ for some $i,j$ with $2\theta+1 \leq i \leq 2\theta'-2\theta$ and $\theta+1 \leq j \leq \theta'-\theta$, and
			\item if $x,y$ are two vertices in the same row of $W'$, then there exist two vertices $x',y'$ in the same row of $W$ such that $\pi_V^*(x)=\pi_V(x')$ and $\pi_V^*(y)=\pi_V(y')$, 
		\end{itemize}
		\item $\bigcup_{e \in E(W')}\pi_E^*(e) \subseteq \bigcup_{e \in E(W)}\pi_E(e)$, and 
		\item for every $v \in V(G)$, $\lvert \cl(v) \cap \pi_V^*(V(W')) \rvert \leq 1$.
	\end{itemize}
	
	\noindent{\bf Proof of Claim 2:}
	Let $I=\{i: \theta+1 \leq i \leq \theta'-\theta\}$ and $J = \{j: \theta+1 \leq j \leq \theta'-\theta\}$.
	For every $(i,j) \in I \times J$, define $f((i,j)) = \{v \in V(G): \{\pi_V((2i-1,j)),\pi_V((2i,j))\} \cap \cl(v) \neq \emptyset\}$.
	Note that $\lvert f((i,j)) \rvert \leq 4$ for each $(i,j) \in I \times J$.
	In addition, for each $v \in V(G)$, $\lvert \{(i,j) \in I \times J: v \in f((i,j))\} \leq \lvert \cl(v) \cap \pi_V(V(W)) \rvert \leq (2\theta d)^2$.
	By Lemma \ref{grid pigeonhole}, there exist $I' \subseteq I$ with $\lvert I' \rvert = 4\theta$ and $J' \subseteq J$ with $\lvert J' \rvert = 2\theta$ such that $f((i,j)) \cap f((i',j')) = \emptyset$ for distinct $(i,j),(i',j') \in I' \times J'$.
	
	Denote the elements of $I'$ by $x_1<x_2<...<x_{4\theta}$ and denote the elements of $J'$ by $y_1<y_2<...<y_{2\theta}$.
	For each $i \in [2\theta]$ and $j \in [2\theta]$, define $\pi_V^*((2i-1,j))=\pi_V((2x_{2i-1}-1,y_j))$ and $\pi_V^*((2i,j))=\pi_V((2x_{2i},y_j))$.
	So $\pi_V^*(V(W')) \subseteq \pi_V(V(W))$, and for every $v \in V(G)$, $\lvert \cl(v) \cap \pi_V^*(V(W')) \rvert \leq 1$.
	Furthermore, if $x,y$ are two vertices in the same row of $W'$, then there exist two vertices $x',y'$ in the same row of $W$ such that $\pi_V^*(x)=\pi_V(x')$ and $\pi_V^*(y)=\pi_V(y')$.
	Note $\theta+1 \leq x_1<x_{4\theta} \leq \theta'-\theta$ and $\theta+1 \leq y_1<y_{2\theta} \leq \theta'-\theta$, so if $x \in V(W')$ and $x' \in V(W)$ with $\pi_V^*(x)=\pi_V(x')$, then $x' = (i,j)$ for some $i,j$ with $2\theta+1 \leq i \leq 2\theta'-2\theta$ and $\theta+1 \leq j \leq \theta'-\theta$.
	It is obvious that one can define $\pi_E^*$ such that $(\pi_V^*,\pi_E^*)$ is a $W'$-subdivision in $L(G)$ such that $\bigcup_{e \in E(W')}\pi_E^*(e) \subseteq \bigcup_{e \in E(W)}\pi_E(e)$.
	$\Box$

\smallskip
	
	Now we define a $W'$-immersion $(\pi_V',\pi_E')$ in $G$.	
	\begin{itemize}
		\item Define $\pi'_V$ to be the function that maps each vertex of $x$ of $W'$ to a vertex $v$ of $G$ such that $\pi_V^*(x) \in \cl(v)$ and $\lvert N_{\Pi(W')}(\pi_V^*(x)) \cap \cl(v) \rvert = \max_{u \in V(G)} \lvert N_{\Pi(W')}(\pi_V^*(x)) \cap \cl(u) \rvert$.
		\item Define $\pi'_E$ to be the function that maps each edge $e$ of $W'$ to the path in $G$ whose edge-set is the union of the set of internal vertices of $\pi^*_E(e)$ and the set $U_x \cup U_y$, where $x,y$ are the ends of $e$ and for each $u \in \{x,y\}$, the set $U_u$ satisfies
		\begin{itemize}
			\item if the vertex in $N_{\pi_E(e)}(u)$ is not in $\cl(\pi'_V(u))$, then $U_u=\{\pi^*_V(u)\}$, and 
			\item if the vertex in $N_{\pi_E(e)}(u)$ is in $\cl(\pi'_V(u))$, then $U_u=\emptyset$.
		\end{itemize}
	\end{itemize}
	It is clear that $(\pi_V',\pi_E')$ is a $W'$-immersion in $G$.
	
	To prove this lemma, it is sufficient to show that $\E_\theta$ is induced by $(\pi_V',\pi_E')$ and the natural edge-tangle of order $\theta$ in $W'$.
	Let $\E''$ be the edge-tangle induced by $(\pi_V',\pi_E')$ and the natural edge-tangle of order $\theta$ in $W'$.
	Note that $\E''$ has order $\theta$.
	
	Suppose to the contrary that $\E_\theta \neq \E''$.
	So there exists $[A,B] \in \E_\theta-\E''$.
	We further assume that the order of $[A,B]$ is as small as possible, and subject to that, $A$ is maximal.
	
	Since $[\{v\}, V(G)-\{v\}] \in \E_\theta \cap \E''$ for every vertex $v$ of $G$ incident with less than $\theta$ edges of $G$, $A$ contains all isolated vertices of $G$ by Lemma \ref{easy edge-tangle} and the maximality of $A$.
	Similarly, by the minimality of the order of $[A,B]$ and Lemma \ref{easy edge-tangle}, every vertex in $A$ that is a non-isolated vertex of $G$ is incident with an edge whose every end is in $A$, and every vertex in $B$ that is a non-isolated vertex of $G$ is incident with an edge whose every end is in $B$.

Define $(A',B')$ to be the separation of $L(G)$ such that 
	\begin{itemize}
		\item[(iv)] $V(A') \cap V(B')$ consists of the vertices of $L(G)$ corresponding to the edges with one end in $A$ and one end in $B$, 
		\item[(v)] $V(B')$ consists of the vertices of $L(G)$ corresponding to the edges of $G$ incident with vertices in $B$, and 
		\item[(vi)] subject to (iv) and (v), $E(A')$ is maximal.
	\end{itemize}

	Since for every non-isolated vertex of $G$, it is not an isolated vertex in $G[A]$ or $G[B]$, so it has a neighbor in the same side of the edge-cut.
	Hence every vertex in $V(A') \cap V(B')$ is adjacent to a vertex in $V(A')-V(B')$ and a vertex in $V(B')-V(A')$. 
	So $(A',B')$ is normalized.
	Since $A$ contains all isolated vertices of $G$, $[A,B]$ is the partner of $(A',B')$.
	Since $[A,B] \in \E_\theta$, $(A',B') \in \Eb_\theta$.
	Since $\Eb_{2\theta'}$ is induced by a $W$-subdivision $(\pi_V,\pi_E)$ in $L(G)$, $E(B')$ intersects every path in the image of $\pi_E$ of all edges of a row of $W$.
	Since the order of $(A',B')$ is less than $\theta$ and $W$ is a $2\theta' \times 2\theta'$ wall, $V(B')-V(A')$ contains at least $2\theta'-2\theta$ vertices in the image of $\pi_V$ of the vertices of a row of $W$.
	Furthermore, for each row of $W$, $V(A')$ contains at most $6\theta$ vertices in the image of $\pi_V$ of the vertices of this row, for otherwise there are at least $\theta$ disjoint paths from $V(A')$ to $V(B')$.
	Hence for every row of $W$, $V(B')-V(A')$ contains the image of $\pi_V$ of some vertex of this row.
	
	Since $[A,B] \not \in \E''$, $[B,A] \in \E''$ by (E1).
	So $A$ contains the image of $\pi_V'$ of all vertices of a row of $W'$ by Lemma \ref{summary edge-tanlge induced by wall immersion}.
	Recall that by the definition of $\pi_V'$, for every $v \in \pi_V'(V(W'))$, $\cl(v) \cap \pi_V^*(V(W')) \neq \emptyset$.
	And by Claim 2, for distinct $v_1,v_2 \in \pi_V'(V(W'))$, $\cl(v_1) \cap \cl(v_2) \cap \pi_V^*(V(W')) = \emptyset$.
	In addition, by Claim 2, if $x,y$ are two vertices in the same row of $W'$, then there exist two vertices $x',y'$ in the same row of $W$ such that $\pi_V^*(x)=\pi_V(x')$ and $\pi_V^*(y)=\pi_V(y')$.
	Hence, there exist a row $r$ of $W$ and a set $S$ such that $S \subseteq V(A')$, and $S$ consist of $\theta$ vertices in the image of $\pi_V$ of vertices in the row $r$.
	By Claim 2, $r$ is the $j$-th row of $W$, for some $\theta+1 \leq j \leq \theta'-\theta$.
	
	Let $S=\{v_1,v_2,...,v_\theta\}$.
	For each $i \in [\theta]$, let $x_i$ be the vertex of $W$ such that $\pi_V(x_i)=v_i$.
	Since $\theta+1 \leq j \leq \theta'-\theta$, there exist distinct rows $r_1,r_2,...,r_\theta$ of $W$ other than $r$ such that for every permutation $\sigma: [\theta] \rightarrow [\theta]$, there exist disjoint paths $P_{\sigma,1},P_{\sigma,2},...,P_{\sigma,\theta}$ in $L(G)$ such that for every $i \in [\theta]$, $P_{\sigma,i}$ is from $v_i$ to $\pi_V(x_{\sigma,i}')$ for some vertex $x_{\sigma,i}'$ of $W$ in the row $r_{\sigma(i)}$ such that $P_{\sigma,i}$ is contained in the image of $\pi_E$ of the column of $W$ containing $x_i$.
	Since for each $i \in [\theta]$, $V(B')-V(A')$ contains the image of $\pi_V$ of some vertex of each $r_i$, we know that for every permutation $\sigma: [\theta] \rightarrow [\theta]$ and for each $i \in [\theta]$, there exists a path $Q_{\sigma,i}$ contained in the image of $\pi_E$ of $r_i$ such that $Q_{\sigma,i}$ is from $\pi_V(x_{\sigma,i}')$ to a vertex in $V(B')-V(A')$.
In addition, we say $Q_{\sigma,i}$ is {\it decreasing} if the index of the column containing $\pi_V(x_{\sigma,i}')$ is at least the index of the column containing the other end of $Q_{\sigma,i}$; otherwise, we say $Q_{\sigma,i}$ is {\it increasing}.

For each $i \in [\theta]$, assuming $\sigma^*(i')$ is defined for every $1 \leq i' \leq i-1$, define $\sigma^*(i)$ to be the element in $[\theta]-\{\sigma^*(i'): 1 \leq i' \leq i-1\}$ such that $Q_{\sigma^*,i}$ is decreasing if possible, and subject to this,
	\begin{itemize}
		\item if $Q_{\sigma^*,i}$ is decreasing, then $P_{\sigma^*,i}$ is as short as possible, and
		\item if $Q_{\sigma^*,i}$ is increasing, then $P_{\sigma^*,i}$ is as long as possible.
	\end{itemize}
Then $\sigma^*$ is a permutation on $[\theta]$.
Since $P_{\sigma^*,1} \cup Q_{\sigma^*,1}, P_{\sigma^*,2} \cup Q_{\sigma^*,2}, ..., P_{\sigma^*,\theta} \cup Q_{\sigma^*,\theta}$ cannot be $\theta$ disjoint paths in $L(G)$ from $V(A')$ to $V(B')-V(A')$, there exist $1 \leq a < b \leq \theta$ such that $P_{\sigma^*,a} \cup Q_{\sigma^*,a}$ intersects $P_{\sigma^*,b} \cup Q_{\sigma^*,b}$.

Suppose that $Q_{\sigma^*,a}$ is decreasing.
Since $P_{\sigma^*,a} \cup Q_{\sigma^*,a}$ intersects $P_{\sigma^*,b} \cup Q_{\sigma^*,b}$, $Q_{\sigma^*,b}$ is decreasing.
But the choice of $P_{\sigma^*,a}$ implies that $P_{\sigma^*,a} \cup Q_{\sigma^*,a}$ and $P_{\sigma^*,b} \cup Q_{\sigma^*,b}$ are disjoint, a contradiction.

So $Q_{\sigma^*,a}$ is increasing.
In particular, the ends of $Q_{\sigma^*,a}$ are not contained in the same column.
If $Q_{\sigma^*,b}$ is increasing, then the choice of $P_{\sigma^*,a}$ implies that $P_{\sigma^*,a} \cup Q_{\sigma^*,a}$ and $P_{\sigma^*,b} \cup Q_{\sigma^*,b}$ are disjoint, a contradiction.
So $Q_{\sigma^*,b}$ is decreasing.
Since $P_{\sigma^*,a} \cup Q_{\sigma^*,a}$ intersects $P_{\sigma^*,b} \cup Q_{\sigma^*,b}$, the choice of $P_{\sigma^*,a}$ implies that the index of the column containing the end of $Q_{\sigma^*,b}$ other than $\pi_V(x'_{\sigma^*,b})$ is at most the index of the column containing $\pi_V(x'_{\sigma^*,a})$, so $Q_{\sigma^*,a}$ can be chosen to be decreasing, a contradiction.
This proves the lemma. 
\end{pf1}

\subsection{Other useful lemmas} \label{subsec:edge-tangles_others}

The following two lemmas will be used in Section \ref{sec:isolating an immersion}.

\begin{lemma} \label{A disjoint}
Let $G$ be a graph and $\E$ an edge-tangle in $G$.
Let $p$ be a positive integer and let $[A_1,B_1], ..., [A_p,B_p] \in \E$.
For each $i$ with $1 \leq i \leq p$, let $X_i$ be the set of edges of $G$ between $A_i$ and $B_i$.
Assume that for every $i$ with $1 \leq i \leq p$ and for every $v \in A_i$, there exists a path in $G[A_i]$ from $v$ to an end of an edge in $X_i$.
Assume the order of $\E$ is greater than $\lvert \bigcup_{i=1}^p X_i \rvert$.
If $\bigcup_{i=1}^p X_i$ is free with respect to $\E$ and $X_i \cap X_j = \emptyset$ for every pair of distinct $i,j$, then $A_i \cap A_j = \emptyset$ for every pair of distinct $i,j$.
\end{lemma}

\begin{pf}
There is nothing to prove if $p=1$, so we may assume that $p \geq 2$.
First, we suppose that there exists an edge $e \in X_1$ such that both ends of $e$ are in $A_2$.
Let $Z = (\bigcup_{i=1}^p X_i)-e$.
Since the order of $\E$ is greater than $\lvert \bigcup_{i=1}^p X_i \rvert$, $[A_2, B_2] \in \E-Z$ has order zero in $G-Z$.
But both ends of the unique member $e$ of $\bigcup_{i=1}^p X_i-Z$ are in $A_2$.
So $\bigcup_{i=1}^p X_i$ is not free with respect to $\E$, a contradiction.
Hence, no edge in $X_1$ has both ends in $A_2$.

Similarly, for every pair of distinct $i,j$, no edge in $X_i$ has both ends in $A_j$.

Now we suppose that there exist distinct $i,j$ such that $A_i \cap A_j \neq \emptyset$.
Let $v \in A_i \cap A_j$.
Let $P_i$ be a path in $G[A_i]$ from $v$ to an end of an edge $e'$ in $X_i$.
Since $X_i$ is disjoint from $X_j$ and some end of $e'$ is in $B_j$, $e'$ has both ends in $B_j$.
So $P_i$ intersects $X_j$.
But $P_i$ is contained in $G[A_i]$, so every edge in $X_j \cap E(P_i)$ has both ends in $A_i$, a contradiction.
This proves the lemma.
\end{pf}

\begin{lemma} \label{no 2 edges}
	Let $\xi$ be a positive integer.
	Let $G$ be a graph and $\E$ an edge-tangle in $G$ of order at least $\xi+2$.
	If $Z$ is a subset of $E(G)$ with $\lvert Z \rvert \leq \xi$, then there exists a set $W$ consisting of two edges of $G-Z$ with at least one common end such that $W$ is free with respect to $\E-Z$.
\end{lemma}

\begin{pf}
	Suppose to the contrary that every set consisting of two edges of $G-Z$ with at least one common end is not free with respect to $\E-Z$.
	That is, for every edges $e_1,e_2 \in E(G)-Z$ sharing at least one common end, there exist $Y \subset \{e_1,e_2\}$ and $[A,B] \in \E-(Y \cup Z)$ of order at most $1 -\lvert Y \rvert$ such that every edge in $\{e_1,e_2\}-Y$ has every end in $A$.
	Let $G_1,...,G_c$ be the components of $G-Z$.
	
	\noindent{\bf Claim 1:} There exists a unique $i$ with $1 \leq i \leq c$ such that $[V(G)-V(G_i),V(G_i)] \in \E-Z$.
	
	\noindent{\bf Proof of Claim 1:}
	We first prove that there exists $i$ with $1 \leq i \leq c$ such that $[V(G)-V(G_i),V(G_i)] \in \E-Z$.
	Suppose to the contrary that $[V(G)-V(G_i),V(G_i)] \not\in \E-Z$ for every $i \in [c]$.
	Since $\E-Z$ has order at least two, by (E1), $[V(G_i),V(G)-V(G_i)] \in \E-Z$ for every $i \in [c]$.
	We prove that $[\bigcup_{j=1}^kV(G_j), V(G)-\bigcup_{j=1}^kV(G_j)] \in \E-Z$ for every $k \in [c]$ by induction on $k$.
	The case $k=1$ is obviously true, so we may assume that $k \geq 2$ and $[\bigcup_{j=1}^{k-1}V(G_j), V(G)-\bigcup_{j=1}^{k-1}V(G_j)] \in \E-Z$.
	Then $[\bigcup_{j=1}^kV(G_j), V(G)-\bigcup_{j=1}^kV(G_j)] \in \E-Z$ by Lemma \ref{easy edge-tangle}.
	Hence $[\bigcup_{j=1}^kV(G_j), V(G)-\bigcup_{j=1}^kV(G_j)] \in \E-Z$ for every $k \in [c]$.
	But when $k=c$, $[V(G),\emptyset] = [\bigcup_{j=1}^cV(G_j), V(G)-\bigcup_{j=1}^cV(G_j)] \in \E-Z$, contradicting (E3).
	This shows the existence of $i$.
	
	Now we show the uniqueness of $i$.
	Suppose there exist distinct $a,b \in [c]$ such that $[V(G)-V(G_a), V(G_a)]$ and $[V(G)-V(G_b),V(G_b)]$ belong to $\E-Z$.
	But $a,b$ are distinct, so $V(G_a) \cap V(G_b) = \emptyset$, contradicting (E2).
	This proves the claim.
	$\Box$

	\smallskip

	Without loss of generality, we may assume that $[V(G)-V(G_1),V(G_1)] \in \E-Z$.
	Define $\E'$ to be the set of edge-cuts of $G_1$ such that $[A,B] \in \E'$ if and only if $[A,B]$ has order less than two and $[A \cup \bigcup_{i=2}^c V(G_i), B] \in \E-Z$.
		
Suppose that $\E'$ is not an edge-tangle in $G_1$ of order two.
It is easy to see that $\E'$ satisfies (E2) and (E3).
So $\E'$ does not satisfy (E1).
Hence there exists an edge-cut $[A,B]$ of $G_1$ such that $[B, A \cup \bigcup_{i=2}^c V(G_i)]$ and $[A, B \cup \bigcup_{i=2}^c V(G_i)]$ belong to $\E-Z$, but these two edge-cuts together with $[\bigcup_{i=2}^c V(G_i), V(G_1)]$ are three edge-cuts in $\E-Z$ such that $(A \cup \bigcup_{i=2}^c V(G_i)) \cap (B \cup \bigcup_{i=2}^c V(G_i)) \cap V(G_1)=\emptyset$, contradicting (E2).

Hence $\E'$ is an edge-tangle in $G_1$ of order two.
	
By considering the cut-edges of $G_1$, it is well-known that there exist a tree $T$ and a partition $(X_t: t \in V(T))$ of $V(G_1)$ such that 
	\begin{itemize}
		\item $G_1[X_t]$ either has only one vertex or is 2-edge-connected for every $t \in V(T)$,
		\item for every adjacent vertices $t_1,t_2$ of $T$, there exists uniquely one edge between $X_{t_1}$ and $X_{t_2}$, and 
		\item every edge of $G_1$ either has every end in $X_t$ for some $t \in V(T)$, or has one end in $X_{t_1}$ and one end in $X_{t_2}$ for some adjacent vertices $t_1,t_2$ of $T$.
	\end{itemize}
For each edge $e=t_1t_2$ of $T$, let $T_{e,t_1}$ and $T_{e,t_2}$ be the components of $T-e$ containing $t_1$ and $t_2$, respectively, and define $Y_{e,t_1} = \bigcup_{t \in V(T_{e,t_1})}X_t$ and $Y_{e,t_2} = \bigcup_{t \in V(T_{e,t_2})} X_t$.
Since $\E'$ has order two, by (E1) and (E2), exactly one of $[Y_{e,t_1}, Y_{e,t_2}]$ and $[Y_{e,t_2}, Y_{e,t_1}] \in \E-Z$.
If the former happens, then we orientate the edge $e$ from $t_1$ to $t_2$; otherwise, we orientate the edge $e$ from $t_2$ to $t_1$.
So we obtain an orientation of $E(T)$ and hence $T$ has a vertex $t^*$ of out-degree zero.
	
Given two edges $e,f$ of $G-Z$ with at least one common end, by the assumption, there exist $Y \subset \{e,f\}$ and $[A_{e,f},B_{e,f}] \in \E-Z$ of order at most $1 - \lvert Y \rvert$ such that every end of the edges in $\{e,f\}-Y$ is in $A$.
Let $[A'_{e,f},B'_{e,f}]$ be the edge-cut $[A_{e,f} \cap V(G_1), B_{e,f} \cap V(G_1)]$ of $G_1$ of order at most  $1- \lvert Y \rvert$.
	If $[B'_{e,f},A'_{e,f}] \in \E'$, then $[(B_{e,f} \cap V(G_1)) \cup \bigcup_{i=2}^cV(G_i), A_{e,f} \cap V(G_1)] \in \E-Z$ by the definition of $\E'$, but $[A_{e,f},B_{e,f}]$ also belongs to $\E-Z$, contradicting (E2).
	So $[B'_{e,f}, A'_{e,f}] \not \in \E'$.
	By (E1), $[A'_{e,f},B'_{e,f}] \in \E'$.
	Since $[A'_{e,f},B'_{e,f}]$ has order at most one, $B'_{e,f}$ contains $X_{t^*}$.
	
	We first claim that $X_{t^*}$ is a single vertex.
	Suppose $X_{t^*}$ contains at least two vertices.
	Then $G_1[X_{t^*}]$ is 2-edge-connected.
	We choose $e,f$ to be two edges of $G_1[X_{t^*}]$ sharing at least one common end.
	By the definition, one of $e,f$ has every end in $A'_{e,f}$.
	But as proved in the previous paragraph, $B'_{e,f}$ contains $X_{t^*}$ and hence contains the ends of $e$ and $f$.
	So $A'_{e,f} \cap B'_{e,f} \neq \emptyset$, a contradiction.
	
	Hence $X_{t^*}$ contains exactly one vertex $v$.
	Since $\E'$ has order at least two, $v$ is incident with at least two edges of $G_1$.
	Let $e,f$ be two edges of $G_1$ incident with $v$.
	Since one of $e,f$ has every end in $A'_{e,f}$, $v \in A'_{e,f}$, a contradiction.
	This proves the lemma.
\end{pf}

\section{Spider theorems}
\label{sec:spider theorems}

The main result of this section is Lemma \ref{colorful edge spider} which is an edge-version of a result (see Lemma \ref{colorful vertex spider} below) that is slightly stronger than a theorem implicitly proved by Robertson and Seymour \cite{rs XXIII} and explicitly proved by Marx and Wollan \cite{mw}.
Lemma \ref{colorful edge spider} enables us to show that given collections of ``interesting sets'' of edges, either we can extend a set free with respect to an edge-tangle by adding many sets from those given collections, or we can delete a bounded number of edges to make some collection of ``interesting sets'' containing no free sets.
This lemma will be frequently used in this paper.

We need the following lemma, which is a slightly stronger form of \cite[Theorem 7.2]{rs XXIII}.
And it can be proved by simply modifying the proof in \cite{rs XXIII}.

\begin{lemma} \label{stronger rs spider}
	Let $h \geq 1$ and $w \geq 0$ be integers.
	Let $\T$ be a tangle in a graph $G$, and let $W \subseteq V(G)$ be free with respect to $\T$ with $\lvert W \rvert \leq w$.
	If $\T$ has order at least $(w+h)^{h+1}+h$, then there exists $W' \subseteq V(G)$ with $W \subseteq W'$ and $\lvert W' \rvert \leq (w+h)^{h+1}$ such that for every $(C,D) \in \T$ of order $\lvert W \rvert+h_C$ with $W \subseteq V(C)$, where $h_C$ is an integer with $h_C<h$, there exists $(A^*,B^*) \in \T$ with $W' \subseteq V(A^* \cap B^*)$, $\lvert V(A^* \cap B^*) - W' \rvert \leq h_C$ and $C \subseteq A^*$.
\end{lemma}

\begin{pf}
	For every $(A,B) \in \T$ and every $v \in V(A) \cap V(B)$, the {\it $\T$-successor of $(A,B)$ via $v$} is the separation $(A',B')$ of $G$ such that 
	\begin{itemize}
		\item[(i)] $v \not \in V(B')$, $A \subseteq A'$ and $B' \subseteq B$,
		\item[(ii)] subject to (i), the order of $(A',B')$ is as small as possible, and
		\item[(iii)] subject to (i) and (ii), $B'$ is minimal.
	\end{itemize} 
	Let $A_0$ be the graph such that $V(A_0)=W$ and $E(A_0)=\emptyset$.
	Let $\T_0 = \{(A_0,G)\}$, and for $i \geq 1$, let $\T_i$ be the set of all $\T$-successors $(A',B')$ of members $(A,B)$ of $\T_{i-1}$ via some vertex in $V(A) \cap V(B)$ with $\lvert V(A') \cap V(B') \rvert < \lvert W \rvert + h$.
	Let $W' = \bigcup_{0 \leq i \leq h} \bigcup_{(A,B) \in \T_i} (V(A) \cap V(B))$.
	It is proved in \cite[Theorem 7.2]{rs XXIII} that $W \subseteq W'$, $\lvert W' \rvert \leq (w+h)^{h+1}$, and every member of $\T_i$ has order at least $\lvert W \rvert + i-1$ for every $i \in [h+1]$.
	
	Let $(C,D) \in \T$ of order $\lvert W \rvert+h_C$ with $W \subseteq V(C)$ for some integer $h_C$ with $h_C<h$.
	Note that $h_C \geq 0$ since $W$ is free with respect to $\T$.
	Let $(C^*,D^*) \in \T$ be the separation of $G$ such that 
	\begin{itemize}
		\item[(iv)] the order of $(C^*,D^*)$ is at most $\lvert W \rvert + h_C$, $C \subseteq C^*$ and $D^* \subseteq D$,
		\item[(v)] subject to (iv), the order of $(C^*,D^*)$ is minimal, and 
		\item[(vi)] subject to (iv) and (v), $C^*$ is maximal.
	\end{itemize}
	Let $(A,B) \in \T_i$ for some $i$ with $0 \leq i \leq h$ such that $A \subseteq C^*$ and $D^* \subseteq B$.
	Note that such an $(A,B)$ exists as $(A_0,G)$ is a candidate.
	We assume that $i$ is as large as possible.
	
	Note that if $i \neq 0$, then $(A,B)$ is a $\T$-successor of a member of $\T_{i-1}$, so either $(A,B)=(C^*,D^*)$ or the order of $(A,B)$ is smaller than the order of $(C^*,D^*)$, for otherwise $(C^*,D^*)$ is a better candidate for being in $\T_i$ than $(A,B)$ by (i)-(iii).
	Furthermore, if $(A,B) \neq (C^*,D^*)$, then $\lvert W \rvert + i-1 \leq \lvert V(A) \cap V(B) \rvert < \lvert V(C^*) \cap V(D^*) \rvert \leq \lvert W \rvert + h_C$, so $i \leq h_C<h$.
	
	Suppose that $V(A) \cap V(B) \not \subseteq V(C^*) \cap V(D^*)$.
	Let $v \in (V(A) \cap V(B)) - (V(C^*) \cap V(D^*))$.
	Let $(A',B')$ be the $\T$-successor of $(A,B)$ via $v$.
	Note that the order of $(A',B')$ is at most the order of $(C^*,D^*)$ as $A \subseteq C^*$ and $D^* \subseteq B$.
	Since $V(A) \cap V(B) \not \subseteq V(C^*) \cap V(D^*)$, $(A,B) \neq (C^*,D^*)$, so $i \leq h-1$.
	By the maximality of $i$, either $A' \not \subseteq C^*$, or $D^* \not \subseteq B'$.
	By the maximality of $C^*$, the order of $(C^* \cup A',D^* \cap B')$ is greater than the order of $(C^*,D^*)$.
	So the order of $(C^* \cap A', D^* \cup B')$ is smaller than the order of $(A',B')$ by the submodularity.
	But $v \not \in V(D^* \cup B')$ and $A \subseteq C^* \cap A'$ and $D^* \cup B' \subseteq B$, so $(A',B')$ is not the $\T$-successor of $(A,B)$ via $v$ by (ii), a contradiction.
	Hence, $V(A) \cap V(B) \subseteq V(C^*) \cap V(D^*)$.
	
	Let $(A^*,B^*)$ be the separation of $G$ such that $V(A^*) = V(C^*) \cup W'$ and $V(B^*) = V(D^*) \cup W'$, and subject to that, $A^*$ is maximal.
	If $i=0$, then the order of $(A,B)$ is $\lvert W \rvert$; if $i \geq 1$, then the order of $(A,B)$ is at least $\lvert W \rvert +i-1$.
	So the order of $(A,B)$ is at least $\lvert W \rvert$.
	Since $V(A) \cap V(B) \subseteq V(C^*) \cap V(D^*) \cap W'$, the order of $(A^*,B^*)$ is at most $\lvert V(C^*) \cap V(D^*) \rvert - \lvert V(A) \cap V(B) \rvert + \lvert W' \rvert \leq (w+h)^{h+1} + h_C$.
	So $(A^*,B^*) \in \T$.
	In addition, $W' \subseteq V(A^*) \cap V(B^*)$ and $C \subseteq A^*$. 
	And $\lvert V(A^*) \cap V(B^*)-W' \rvert \leq \lvert V(C^*) \cap V(D^*) \rvert - \lvert V(A) \cap V(B) \rvert \leq h_C$.
\end{pf}

\bigskip

For every tangle $\T$ of order $\theta$ in a graph $G$ and every $Z \subseteq V(G)$ with $\lvert Z \rvert < \theta$, we define $\T-Z$ to be the set of separations $(A,B)$ of $G-Z$ such that $(A', B') \in \T$ for some subgraphs $A',B'$ of $G$ with $V(A')=V(A) \cup Z$ and $V(B') = V(B) \cup Z$.
Note that $\T-Z$ is a tangle in $G-Z$ of order $\theta-\lvert Z \rvert$ by \cite[Theorem 6.2]{rs X}.

The following is a stronger form of \cite[Theorem 3.3]{mw} and its proof uses ideas similar to that used in \cite[Theorem 3.3]{mw}.

\begin{lemma} \label{colorful vertex spider}
	Let $G$ be a graph and $\T$ a tangle in $G$ of order $\theta$, and let $c$ be a positive integer.
	For every $i \in [c]$, let $d_i,k_i$ be positive integers, and let $\{X_{i,j} \subseteq V(G): j \in J_i\}$ be a family of subsets of $V(G)$ indexed by a set $J_i$.
	Let $d,k$ be integers such that $\theta \geq (kcd)^{d+1}+d$, $d_i \leq d$ and $k_i \leq k$ for $i \in [c]$.
	Let $J_i^* \subseteq J_i$ with $\lvert J_i^* \rvert \leq k_i$ for each $i \in [c]$, such that $\bigcup_{i=1}^c \bigcup_{j \in J_i^*} X_{i,j}$ is free with respect to $\T$ and $X_{i,j} \cap X_{i',j'} = \emptyset$ for distinct pairs $(i,j),(i',j')$ with $1 \leq i \leq i' \leq c$, $j \in J_i^*$ and $j' \in J_{i'}^*$.
	If $\lvert X_{i,j} \rvert \leq d_i$ for every $i \in [c]$ and $j \in J_i$, then either 
	\begin{enumerate}
		\item there exist $J_1',J_2',...,J_c'$ with $J_i^* \subseteq J_i' \subseteq J_i$ and $\lvert J_i' \rvert = k_i$ for each $i \in [c]$ such that $\bigcup_{i \in [c]} \bigcup_{j \in J'_i} X_{i,j}$ is free with respect to $\T$, and $X_{i,j} \cap X_{i',j'} = \emptyset$ for all distinct pairs $(i,j),(i',j')$ with $1 \leq i \leq i' \leq c$, $j \in J_i'$ and $j'\in J_{i'}'$, or
		\item there exist $Z \subseteq V(G)$ with $\lvert Z \rvert \leq (kcd)^{d+1}$ and integer $i^*\in [c]$ with $\lvert J_{i^*}^* \rvert < k_{i^*}$ such that for every $j \in J_{i^*}$, either $X_{i^*,j} \cap Z \neq \emptyset$, or $X_{i^*,j}$ is not free with respect to $\T-Z$.
	\end{enumerate}
\end{lemma}

\begin{pf}
	For every $i \in [c]$, pick $J_i'$ with $J_i^* \subseteq J'_i \subseteq J_i$ and $\lvert J'_i \rvert \leq k_i$ such that
	\begin{enumerate}
		\item[(i)] $\bigcup_{i \in [c]} \bigcup_{j \in J'_i} X_{i,j}$, denoted by $W$, is free with respect to $\T$, 
		\item[(ii)] $X_{i,j}$ and $X_{i',j'}$ are disjoint for all distinct pairs $(i,j),(i',j')$ with $1 \leq i \leq i' \leq c$, $j \in J_i'$ and $j'\in J_{i'}'$, and
		\item[(iii)] subject to (i) and (ii), the sequence $(k_1-\lvert J_1' \rvert, k_2-\lvert J'_2 \rvert, ..., k_c-\lvert J'_c \rvert)$, denoted by $s$, is lexicographically minimal.
	\end{enumerate}
	Note that such a set $W$ exists since $\bigcup_{i=1}^c \bigcup_{j \in J_i^*} X_{i,j}$ is free and $X_{i,j} \cap X_{i',j'} = \emptyset$ for all distinct pairs $(i,j),(i',j')$ with $1 \leq i \leq i' \leq c'$, $j \in J_i^*$ and $j'\in J_{i'}^*$.
	
	Assume that the first conclusion of this lemma does not hold.
	So $s$ contains a non-zero entry.
	Let $i^*$ be the smallest number $i$ such that $k_i-\lvert J'_i \rvert >0$.
	Note that $\lvert W \rvert \leq (kc-1)d$.
	Applying Lemma \ref{stronger rs spider} by taking $(h,w)=(d,(kc-1)d)$, there exists $Z \subseteq V(G)$ with $W \subseteq Z$ and $\lvert Z \rvert \leq (kcd)^{d+1}$ such that for every $(C,D) \in \T$ of order $\lvert W \rvert + h_C$ with $W \subseteq V(C)$ for some $h_C < d$, there exists $(A',B') \in \T$ with $Z \subseteq V(A' \cap B')$, $\lvert V(A' \cap B') - Z \rvert \leq h_C$ and $C \subseteq A'$.
	
	We shall prove that $Z$ and $i^*$ satisfy the second conclusion of this lemma.
	Assume that $j \in J_{i^*}$ such that $X_{i^*,j} \cap Z = \emptyset$.
	Since $W \subseteq Z$, $X_{i^*,j}$ is disjoint from $W$.
	By the maximality of $W$, $W \cup X_{i^*,j}$ is not free with respect to $\T$.
	So there exists a separation $(C,D) \in \T$ of order at most $\lvert W \rvert + \lvert X_{i^*,j} \rvert-1$ with $W \cup X_{i^*,j} \subseteq V(C)$.
	By the choice of $Z$, there exists $(A',B') \in \T$ with $Z \subseteq V(A' \cap B')$, $\lvert V(A' \cap B')-Z \rvert \leq \lvert X_{i^*,j} \rvert-1$ and $C \subseteq A'$.
	That is, $X_{i^*,j} \subseteq V(A')-Z$ and the order of $(A'-Z,B'-Z)$ is less than $\lvert X_{i^*,j} \rvert$.
	So $X_{i^*,j}$ is not free with respect to $\T-Z$.
	This proves the lemma.
\end{pf}

\bigskip

What we really need in this paper is a version of Lemma \ref{colorful vertex spider} with respect to edge-tangles.

\begin{lemma} \label{colorful edge spider}
	Let $G$ be a graph and $\E$ an edge-tangle in $G$ of order $\theta$, and let $c$ be a positive integer.
	For every $i \in [c]$, let $d_i,k_i$ be positive integers, and let $\{X_{i,j} \subseteq E(G): j \in J_i\}$ be a family of subsets of $E(G)$ indexed by a set $J_i$.
	Let $d,k$ be integers such that $\theta \geq 3(kcd)^{d+1}+3d$, $d_i \leq d$ and $k_i \leq k$ for $i \in [c]$.
	Let $J_i^* \subseteq J_i$ with $\lvert J_i^* \rvert \leq k_i$ for each $i \in [c]$, such that $\bigcup_{i=1}^c \bigcup_{j \in J_i^*} X_{i,j}$ is free with respect to $\E$ and $X_{i,j} \cap X_{i',j'} = \emptyset$ for distinct pairs $(i,j),(i',j')$ with $1 \leq i \leq i' \leq c$, $j \in J_i^*$ and $j' \in J_{i'}^*$.
	If $\lvert X_{i,j} \rvert \leq d_i$ for every $i \in [c]$ and $j \in J_i$, then either 
	\begin{enumerate}
		\item there exist $J_1',J_2',...,J_c'$ with $J_i^* \subseteq J_i' \subseteq J_i$ and $\lvert J_i' \rvert = k_i$ for each $i \in [c]$ such that $\bigcup_{i=1}^c \bigcup_{j \in J'_i} X_{i,j}$ is free with respect to $\E$, and $X_{i,j}$ and $X_{i',j'}$ are disjoint for all distinct pairs $(i,j),(i',j')$ with $1 \leq i \leq i' \leq c$, $j \in J_i'$ and $j'\in J_{i'}'$, or
		\item there exist $Z \subseteq E(G)$ with $\lvert Z \rvert \leq (kcd)^{d+1}$ and integer $i^* \in [c]$ with $\lvert J_{i^*}^* \rvert < k_{i^*}$ such that for every $j \in J_{i^*}$, either $X_{i^*,j} \cap Z \neq \emptyset$, or $X_{i^*,j}$ is not free with respect to $\E-Z$.
	\end{enumerate}
\end{lemma}

\begin{pf}
	Since $\E$ is an edge-tangle of order $\theta$ in $G$, $\Eb$ is a tangle of order at least $\lfloor \theta/3 \rfloor \geq (kcd)^{d+1}+d$ in $L(G)$ by Lemma \ref{conjugate of edge-tangle}.
	Note that for every $i \in [c]$ and $j \in J_i$, $X_{i,j}$ is a subset of $E(G)$ so it is a subset of $V(L(G))$. 
	Since $\bigcup_{i=1}^c \bigcup_{j \in J_i^*} X_{i,j}$ is free with respect to $\E$, $\bigcup_{i=1}^c\bigcup_{j \in J_i^*} X_{i,j}$ is free with respect to $\Eb$ by Lemma \ref{free wrt edge free wrt vertex}.
	So by Lemma \ref{colorful vertex spider},  either
	\begin{itemize}
		\item[(i)] there exist $J_1',J_2',...,J_c'$ with $J_i^* \subseteq J_i' \subseteq J_i$ and $\lvert J_i' \rvert = k_i$ for each $i \in [c]$ such that $\bigcup_{i=1}^c \bigcup_{j \in J'_i} X_{i,j}$ is free with respect to $\Eb$, and $X_{i,j}$ and $X_{i',j'}$ are disjoint for all distinct pairs $(i,j),(i',j')$ with $1 \leq i \leq i' \leq c$, $j \in J_i'$ and $j'\in J_{i'}'$, or
		\item[(ii)] there exist $Z \subseteq V(L(G))=E(G)$ with $\lvert Z \rvert \leq (kcd)^{d+1}$ and integer $i^* \in [c]$ with $\lvert J_{i^*}^* \rvert <k_{i^*}$ such  that for every $j \in J_{i^*}$, either $X_{i^*,j} \cap Z \neq \emptyset$, or $X_{i^*,j}$ is not free with respect to $\Eb-Z$.
	\end{itemize}
	
	If (i) holds, then $\bigcup_{i=1}^c \bigcup_{j \in J'_i} X_{i,j}$ is a set of size at most $ckd \leq \lfloor \theta/3 \rfloor$ that is free with respect to $\Eb$, so $\bigcup_{i=1}^c\bigcup_{j \in J'_i} X_{i,j}$ is free with respect to $\E$ by Lemma \ref{free wrt vertex free wrt edge}, and hence Statement 1 of this lemma holds.
	
	So we may assume that $Z$ and $i^*$ mentioned in (ii) exist.
	We shall prove that Statement 2 of this lemma holds.
	Suppose to the contrary that there exists $j \in J_{i^*}$ such that $X_{i^*,j} \cap Z = \emptyset$ and $X_{i^*,j}$ is free with respect to $\E-Z$.
	So $X_{i^*,j} \cap Z=\emptyset$ and $X_{i^*,j}$ is free with respect to $\Eb-Z$ by Lemma \ref{free wrt edge free wrt vertex}, contradicting (ii).
	This proves the lemma.
\end{pf}

\section{Excluding immersions}
\label{sec:excluding immersions structure}

Given a simple graph $H$, an $H$-{\it minor} of a graph $G$ is a map $\alpha$ with domain $V(H)$ such that
\begin{itemize}
	\item $\alpha(h)$ is a nonempty connected subgraph of $G$, for every $h \in V(H)$;
	\item if $h_1$ and $h_2$ are different vertices in $H$, then $\alpha(h_1)$ and $\alpha(h_2)$ are disjoint;
	\item if $h_1h_2$ is an edge in $H$, then there exists an edge of $G$ with one end in $\alpha(h_1)$ and one end in $\alpha(h_2)$.
\end{itemize}
We say that {\it $G$ contains an $H$-minor} if such a function $\alpha$ exists.
And for every $h \in V(H)$, $\alpha(h)$ is called a {\it branch set} of $\alpha$.

Given a simple graph $H$, an {\it $H$-thorns} of a graph $G$ is a map $\alpha$ with domain $V(H)$ such that
\begin{itemize}
	\item $\alpha(h)$ is a connected subgraph of $G$ with at least one edge, for every $h \in V(H)$;
	\item if $h_1$ and $h_2$ are different vertices in $H$, then $\alpha(h_1)$ and $\alpha(h_2)$ are edge-disjoint;
	\item if $h_1h_2$ is an edge in $H$, then $V(\alpha(h_1)) \cap V(\alpha(h_2)) \neq \emptyset$;
	\end{itemize}
We say that {\it $G$ contains an $H$-thorns} if such a function $\alpha$ exists.
And for every $h \in V(H)$, $\alpha(h)$ is called a {\it branch set} of $\alpha$.

Note that if a graph contains a vertex $v$ incident with $d$ edges, then it contains a $K_d$-thorns whose branch sets are the edges incident with $v$.
Another example of thorns is that every $r \times r$-grid contains a $K_r$-thorns by defining $\alpha(v_i)$ to be the union of the $i$-th row and the $i$-th column.

\begin{lemma} \label{edge-minor and minor}
If $H$ is a simple graph, then a graph $G$ contains an $H$-thorns if and only if $L(G)$ contains an $H$-minor.
\end{lemma}

\begin{pf}
Let $\alpha$ be an $H$-thorns in $G$.
For every $h \in V(H)$, define $\beta(h)$ to be the subgraph of $L(G)$ induced by $E(\alpha(h))$.
It is clear that $\beta$ is an $H$-minor in $L(G)$.

Let $\beta'$ be an $H$-minor in $L(G)$.
For every $h \in V(H)$, define $\alpha'(h)$ to be the connected subgraph of $G$ with $E(\alpha'(h)) = V(\beta'(h))$.
Then it is obvious that $\alpha'$ is an $H$-thorns in $G$.
\end{pf}

\bigskip

The following was proved by Robertson and Seymour \cite{rs XIII}.

\begin{lemma}[{\cite[Theorem 5.4]{rs XIII}}] \label{vertex linkage}
Let $G$ be a graph, and let $Z$ be a subset of $V(G)$ with $\lvert Z \rvert = \xi$.
Let $k \geq \lfloor \frac{3}{2} \xi \rfloor$, and let $\alpha$ be a $K_k$-minor in $G$.
If there is no separation $(A,B)$ of $G$ of order less than $\lvert Z \rvert$ such that $Z \subseteq V(A)$ and $A \cap \alpha(h)=\emptyset$ for some $h \in V(K_k)$, then for every partition $(Z_1,...,Z_n)$ of $Z$ into non-empty subsets, there are $n$ connected subgraphs $T_1, ..., T_n$ of $G$, mutually disjoint and with $V(T_i) \cap Z = Z_i$ for $1 \leq i \leq n$.
\end{lemma}

Now, we prove an edge-variant of Lemma \ref{vertex linkage}.

\begin{lemma} \label{edge linkage}
Let $G$ be a graph, and let $X$ be a subset of $E(G)$ with $\lvert X \rvert = \xi$.
Let $k \geq \lfloor \frac{3}{2} \xi \rfloor$, and let $\alpha$ be a $K_k$-thorns in $G$.
If there exist no $Y \subseteq X$ and edge-cut $[A,B]$ of $G-Y$ of order less than $\xi-\lvert Y \rvert$ such that every edge in $X-Y$ is incident with some vertex in $A$ and $A \cap V(\alpha(h))=\emptyset$ for some $h \in V(K_k)$, then for every partition $(X_1,...,X_n)$ of $X$ into non-empty subsets, there are $n$ connected subgraphs $T_1,...,T_n$ of $G$, mutually edge-disjoint and with $E(T_i) \cap X = X_i$ for $1 \leq i \leq n$.
\end{lemma}

\begin{pf}
Let $\beta$ be the $K_k$-minor in $L(G)$ corresponding to $\alpha$ mentioned in Lemma \ref{edge-minor and minor}.

\noindent{\bf Claim 1:} There does not exist a separation $(A',B')$ of $L(G)$ of order less than $\xi$ such that $X \subseteq V(A')$ and $A' \cap \beta(h) = \emptyset$ for some $h \in V(K_k)$.

\noindent{\bf Proof of Claim 1:}
Suppose to the contrary that there exists a separation $(A',B')$ of $L(G)$ of order less than $\xi$ such that $X \subseteq V(A')$ and $A' \cap \beta(h) = \emptyset$ for some $h \in V(K_k)$.
We may assume that the order of $(A',B')$ is as small as possible.
So every vertex in $V(A') \cap V(B') - X$ must have an neighbor in $V(A')-V(B')$ and a neighbor in $V(B')-V(A')$, and every vertex in $V(A') \cap V(B') \cap X$ has a neighbor in $V(B')-V(A')$.
Define $B = \{v \in V(G): \cl(v) \subseteq V(B')\}$ and $A = V(G)-B$.
Then $[A,B]$ is an edge-cut of $G$.
Let $Y$ be the subset of $X$ consisting of the edges in $X$ with every end in $B$.
Note that the order of $[A,B]$ equals $\lvert V(A') \cap V(B') \rvert - \lvert \{v \in V(A') \cap V(B'): v$ has no neighbor in $V(A')-V(B')\} \rvert = \lvert V(A') \cap V(B') \rvert - \lvert Y \rvert < \xi-\lvert Y \rvert$.
Furthermore, every edge in $X-Y$ has an end in $A$.
In addition, every vertex of $\beta(h)$ is in $V(B')-V(A')$, so every edge of $\alpha(h)$ has every end in $B$.
That is, $V(\alpha(h)) \cap A = \emptyset$, a contradiction.
$\Box$

\smallskip

Let $(X_1,X_2,...,X_n)$ be a partition of $X$ into nonempty sets.
By Lemma \ref{vertex linkage} and Claim 1, there exist mutually disjoint connected subgraphs $T_1',...,T_n'$ of $L(G)$ such that $V(T'_i) \cap X = X_i$ for every $1 \leq i \leq n$.
For every $1 \leq i \leq n$, define $T_i$ to be the connected subgraph of $G$ with $E(T_i) = V(T_i')$. 
Then $T_1,...,T_n$ are mutually edge-disjoint and $E(T_i) \cap X = X_i$.
\end{pf}

\bigskip

A tangle $\T$ in $G$ {\it controls} an $H$-minor $\alpha$ if there do not exist $(A,B) \in \T$ of order less than $\lvert V(H) \rvert$ and $h \in V(H)$ such that $V(\alpha(h)) \subseteq V(A)$.
An edge-tangle $\E$ in $G$ {\it controls} an $H$-thorns $\alpha$ if $V(\alpha(h)) \cap B \neq \emptyset$ for every $h \in V(H)$ and every $[A,B] \in \E$ of order less than $\lvert V(H) \rvert$.

The {\it degree sequence} of a graph $G$ is the non-increasing sequence of the degrees of the vertices of $G$.

\begin{lemma} \label{edge spider immersion}
Let $G$ be a graph and $H$ be a graph on $h$ vertices with degree sequence $(d_1,d_2,...,d_h)$.
Let $d=d_1$ and $t \geq 3hd$.
Let $V(H)=\{u_1,u_2,...,u_h\}$, where $\deg_H(u_i)=d_i$ for every $i \in [h]$.
Let $\E$ be an edge-tangle of order at least $2hd$ in $G$ that controls a $K_t$-thorns.
Let $\ell$ be the number of loops of $H$.
Assume that there exist pairwise disjoint subsets $X_0,X_1,X_2,...,X_h$ of $E(G)$ such that $\bigcup_{i=0}^h X_i$ is free with respect to $\E$, $X_0$ can be partitioned into $\ell$ 2-element subsets $S_1,S_2,...,S_\ell$ where the two edges in each $S_j$ share at least one common end $s_j$, and for each $i \in [h]$, $X_i$ consists of $d_i$ edges incident with a common vertex $v_i$.
If $v_1,v_2,...,v_h$ are distinct and there exists a partition of $\{S_1,S_2,...,S_\ell\}$ into sets $D_1,D_2,...,D_h$ such that for every $i \in [h]$, $\lvert D_i \rvert$ equals the number of loops incident with $u_i$ and $v_i \not \in \{s_j: S_j \in D_i\}$, then $G$ has an $H$-immersion $(\pi_V,\pi_E)$ with $\pi_V(V(H)) = \{v_1,v_2,...,v_h\}$.
\end{lemma}

\begin{pf}
Let $\alpha$ be a $K_t$-thorns in $G$ controlled by $\E$, and let $X = \bigcup_{i=0}^h X_i$.
Note that $\lvert X \rvert \leq 2\ell + hd \leq 2hd$.

\noindent{\bf Claim 1:} For every positive integer $r$ and every partition $(Z_1,Z_2,...,Z_r)$ of $X$ into non-empty subsets, there exist pairwise edge-disjoint connected subgraphs $T_1,T_2,...,T_r$ of $G$ such that $E(T_i) \cap X = Z_i$ for every $1 \leq i \leq r$.

\noindent{\bf Proof of Claim 1:}
Suppose that there exist $Y \subseteq X$ and an edge-cut $[A,B]$ of $G-Y$ of order less than $\lvert X-Y \rvert$ such that every edge in $X-Y$ is incident with some vertex in $A$ and $A \cap V(\alpha(u)) = \emptyset$ for some $u \in V(K_t)$.
We assume that $Y$ is maximal, so every edge in $X-Y$ has every end in $A$.
Since $X$ is free with respect to $\E$, $[A,B] \not \in \E-Y$.
But the order of $[A,B]$ in $G-Y$ is less than $\lvert X-Y \rvert \leq 2hd-\lvert Y \rvert$.
So $[B,A] \in \E-Y$ by (E1).
Hence $[B,A] \in \E$ is an edge-cut of $G$ of order less than $2hd \leq t$.
However, $\E$ controls $\alpha$, so $A \cap V(\alpha(u)) \neq \emptyset$, a contradiction.
Therefore, this claim follows from Lemma \ref{edge linkage}. 
$\Box$

\smallskip

Let $E(H) = \{e_1,e_2,...,e_{\lvert E(H) \rvert}\}$, where $e_j$ is not a loop for every $j \in [\lvert E(H) \rvert-\ell]$.
For every $i \in [h]$, let $Y_i$ be a subset of $X_i$ such that $\lvert Y_i \rvert$ equals the number of non-loops incident with $u_i$.
For every $i \in [h]$, define a bijection $f_i$ from $Y_i$ to the set of non-loop edges of $H$ incident with $u_i$, and define an onto function $f_i'$ from $X_i-Y_i$ to the set of loops of $H$ incident with $u_i$ such that the preimage of every loop incident with $u_i$ has size two.
Define $f_E$ to be a bijection from the set of loops of $H$ to $\{S_1,S_2,...,S_\ell\}$ such that if $e$ is a loop of $H$ incident with $u_i$ for some $i \in [h]$, then $f_E(e) \in D_i$.

For each $i \in [\lvert E(H) \rvert-\ell]$, define $Z_i$ to be the subset of $X$ consisting of the two edges in $\bigcup_{j=1}^hY_j$ mapped to $e_i$ by $f_1,f_2,...,f_h$.
For each loop $e$ of $H$, let $\{Z_{e,1},Z_{e,2}\}$ be a partition of the union of $f_E(e)$ and the preimage of $e$ by $f_1',f_2',...,f_h'$ into two sets of size two such that $\lvert Z_{e,1} \cap f_E(e) \rvert = \lvert Z_{e,2} \cap f_E(e) \rvert=1$.
So $\{Z_1,Z_2,...,Z_{\lvert E(H) \rvert-\ell}, Z_{e,1},Z_{e,2}: e$ is a loop of $H\}$ is a partition of $X$ into non-empty sets.

By Claim 1, there exist pairwise edge-disjoint connected subgraphs $T_1,T_2,...,T_{\lvert E(H) \rvert-\ell}, T_{e,1}, \allowbreak T_{e,2}$ of $G$ (for every loop $e$ of $H$) such that $E(T_i) \cap X = Z_i$ for every $1 \leq i \leq \lvert E(H) \rvert-\ell$, and $E(T_{e,j}) \cap X = Z_{e,j}$ for every loop $e$ of $H$ and $j \in [2]$.
Note that for every $i \in [\lvert E(H) \rvert-\ell]$, there exists a path in $T_i$ connecting $v_j,v_k$, where $j,k$ are the indices such that $e_i$ belongs to the image of $f_j$ and $f_k$.
For each loop $e$ of $H$, there exists a cycle contained in $T_{e,1} \cup T_{e,2}$ containing $v_i$, where $i$ is the index such that $e$ is incident with $u_i$, since $v_i \not \in \{s_j: S_j \in D_i\}$. 

Define $\pi_V: V(H) \rightarrow V(G)$ such that $\pi_V(u_i)=v_i$ for every $i \in [h]$.
Define $\pi_E$ to be a function that maps each non-loop edge $e_i$ of $H$ (for some $i \in [\lvert E(H) \rvert-\ell])$) to be a path in $T_i$ from $\pi_V(u)$ to $\pi_V(u')$, where $u,u'$ are the ends of $e_i$, and maps each loop $e$ of $H$ to a cycle in $T_{e,1} \cup T_{e,2}$ containing $\pi_V(u'')$, where $u''$ is the end of $e$.
Then $(\pi_V,\pi_E)$ is an $H$-immersion in $G$ with $\pi_V(V(H))=\{v_1,v_2,...,v_h\}$.
\end{pf}

\bigskip

A family $\D$ of edge-cuts of a graph is {\it cross-free} if $A \cap C = \emptyset$ for every pair of distinct edge-cuts $[A,B],[C,D]$ in $\D$.

\begin{lemma} \label{making cross-free}
Let $k,\theta$ be integers.
Let $G$ be a graph and $\E$ an edge-tangle in $G$ of order at least $\theta$.
If there exist $C \subseteq E(G)$ with $\lvert C \rvert \leq \theta-k$ and a subset $\D$ of $\E-C$ such that every member of $\D$ is an edge-cut of $G-C$ of order less than $k$, then there exists a cross-free family $\D^* \subseteq \E-C$ such that every member of $\D^*$ is an edge-cut of $G-C$ of order less than $k$ such that $\bigcup_{[A,B] \in \D}A = \bigcup_{[A,B] \in \D^*}A$.
\end{lemma}

\begin{pf}
Define $\D^*$ to be a subset of $\E-C$ with $\bigcup_{[A,B] \in \D}A=\bigcup_{[A,B] \in \D^*}A$ such that every member of $\D^*$ is an edge-cut of $G-C$ of order less than $k$, and subject to that, $\sum_{[A,B] \in \D^*} \lvert A \rvert$ is as small as possible.
Note that such a family $\D^*$ exists as $\D$ is a candidate.
To prove this lemma, it suffices to show that ${\mathcal D}^*$ is cross-free.

Suppose that ${\mathcal D}^*$ is not cross-free.
Then there exist $[A_1, B_1], [A_2,B_2] \in {\mathcal D}^*$ such that $A_1 \neq A_2$ and $A_1 \cap A_2 \neq \emptyset$. 
By the submodularity, $\lvert [A_1 \cap B_2, B_1 \cup A_2] \rvert + \lvert [A_1 \cup B_2, B_1 \cap A_2] \rvert \leq \lvert [A_1, B_1] \rvert + \lvert [B_2, A_2] \rvert \leq 2(k-1)$, so one of $[A_1 \cap B_2, B_1 \cup A_2]$ and $[B_1 \cap A_2, A_1 \cup B_2]$ has order at most $k-1$.
By symmetry, we may assume that $[A_1 \cap B_2, B_1 \cup A_2]$ has order at most $k-1$.
Note that the order of $\E-C$ is at least $\theta-\lvert C \rvert \geq k$.
By Lemma \ref{easy edge-tangle}, $[A_1 \cap B_2, B_1 \cup A_2] \in \E-C$, since $[A_1,B_1] \in \E-C$.
Let ${\mathcal D}' = ({\mathcal D}^* - \{[A_1,B_1]\}) \cup \{[A_1 \cap B_2, A_2 \cup B_1]\}$.
Since $A_1 \subseteq (A_1 \cap B_2) \cup A_2$, ${\mathcal D}'$ is contained in $\E-C$ and is a family of edge-cuts of $G-C$ of order at most $k-1$ such that $\bigcup_{[A,B] \in \D'}A = \bigcup_{[A,B] \in \D^*}A = \bigcup_{[A,B] \in \D}A$.
Hence, by the minimality of ${\mathcal D}^*$, $\lvert A_1 \cap B_2 \rvert \geq \lvert A_1 \rvert$.
This implies that $A_1 = A_1 \cap B_2 \subseteq B_2$, so $A_1 \cap A_2 = \emptyset$, a contradiction.
Therefore, $\D^*$ is cross-free.
\end{pf}

\bigskip

A graph is {\it exceptional} if it contains exactly one vertex of degree at least two, and this vertex is incident with a loop.

Theorem \ref{excluding immersion} is a structure theorem for excluding a fixed non-exceptional graph as an immersion in a graph with an edge-tangle controlling a big complete graph-thorns.

\begin{theorem} \label{excluding immersion}
For any positive integers $d,h$, there exist positive integers $\theta=\theta(d,h)$ and $\xi=\xi(d,h)$ such that the following holds.
If $H$ is a non-exceptional graph with degree sequence $(d_1,d_2,...,d_h)$, where $d_1=d$, and $G$ is a graph that does not contain an $H$-immersion, then for every edge-tangle $\E$ of order at least $\theta$ in $G$ controlling a $K_{3dh}$-thorns, there exist $C \subseteq E(G)$ with $\lvert C \rvert \leq \xi$, $U \subseteq V(G)$ with $\lvert U \rvert \leq h-1$ and a cross-free family ${\mathcal D} \subseteq \E-C$ such that for every vertex $v \in V(G)-U$, there exists $[A,B] \in {\mathcal D}$ of order at most $d_{\lvert U \rvert+1}-1$ with $v \in A$.
\end{theorem}

Now we sketch the proof of Theorem \ref{excluding immersion}.
We greedily pick a vertex $v$ and a set $X_v$ of sufficiently many edges incident with it such that $v$ is not picked before, $X_v$ is disjoint from all previously picked sets of edges, and the union of $X_v$ and all previously picked sets is free with respect to $\E$, until we cannot find such a vertex or a such set of edges.
We first assume that $H$ has no loops.
If we picked at least $\lvert V(H) \rvert$ vertices in the process, then we can construction an $H$-immersion by Lemma \ref{edge spider immersion}, a contradiction.
So the set $U$ of picked vertices has size at most $\lvert V(H) \rvert-1$.
If we can further repeatedly pick a vertex $v$ and a set $X_v$ of $d_{\lvert U \rvert+1}$ edges incident with $v$ such that $v$ is not picked before, $X_v$ is disjoint from all previously picked sets of edges, and the union of $X_v$ and all previously picked sets is free with respect to $\E$, until we get $\lvert V(H) \rvert$ vertices, then again we can construct an $H$-immersion, a contradiction.
So Lemma \ref{colorful edge spider} implies that one can delete a bounded number of edges such that each vertex in $V(G)-U$ is contained in the first entry of an edge-cut in the edge-tangle, and we are done.
The case that $H$ has loops is similar but takes extra work.
Lemma \ref{edge spider immersion} implies that we cannot further pick many disjoint set of two edges sharing a common end such that the union of those sets is free, so that Lemma \ref{colorful edge spider} implies that for every vertex, there exists an edge-cut of order at most one such that this vertex belongs to the first entry of the edge-cut.

\bigskip

\begin{pf1}{Theorem \ref{excluding immersion}}
For any positive integers $d,h$, define $\xi(d,h)=(h+1)((3hd^2)^{d+1} + dh)$ and $\theta(d,h) = 3(2hd^2)^{d+1}+3dh+\xi$. 

Let $d,h$ be positive integers.
Denote $\xi(d,h)$ and $\theta(d,h)$ by $\xi$ and $\theta$, respectively.
Let $H$ be a non-exceptional graph on $h$ vertices with degree sequence $(d_1,d_2,...,d_h)$ and $d_1=d$.
Since there exists no graph on one vertex with maximum degree one, this theorem holds if $d=h=1$.
Suppose that $(d,h)$ is a pair of positive integers with $d+h$ minimum such that this theorem does not hold.
That is, there exists a graph $G$ that does not contain an $H$-immersion and there exists an edge-tangle $\E$ in $G$ of order at least $\theta$ controlling a $K_{3dh}$-thorns such that there do not exist $C \subseteq E(G)$ with $\lvert C \rvert \leq \xi$, $U \subseteq V(G)$ with $\lvert U \rvert \leq h-1$ and a cross-free family ${\mathcal D} \subseteq \E-C$ such that for every vertex $v \in V(G)-U$, there exists $[A,B] \in {\mathcal D}$ of order at most $d_{\lvert U \rvert+1}-1$ with $v \in A$.

\noindent{\bf Claim 1:} $\lvert V(G) \rvert \geq h$ and $H$ does not contain an isolated vertex.

\noindent{\bf Proof of Claim 1:}
It is clear that $\lvert V(G) \rvert \geq h$, for otherwise choosing $C=\emptyset$, $U=V(G)$ and $\D=\emptyset$ leads to a contradiction.
Suppose that $H$ contains an isolated vertex $u$.
Let $H'=H-u$.
Note that the degree sequence of $H'$ is $(d_1,d_2,...,d_{h-1})$, and $H'$ is non-exceptional.
Since $\lvert V(G) \rvert \geq h$ and $G$ does not contain an $H$-immersion, $G$ does not contain an $H'$-immersion.
By the minimality of $d+h$, there exist $C \subseteq E(G)$ with $\lvert C \rvert \leq \xi(d,h-1) \leq \xi(d,h)$, $U \subseteq V(G)$ with $\lvert U \rvert \leq (h-1)-1$ and a cross-free family ${\mathcal D} \subseteq \E-C$ such that for every vertex $v \in V(G)-U$, there exists $[A,B] \in {\mathcal D}$ of order at most $d_{\lvert U \rvert+1}-1$ with $v \in A$, a contradiction.
$\Box$

\smallskip

\noindent{\bf Claim 2:} There do not exist $C \subseteq E(G)$ with $\lvert C \rvert \leq \xi$ and $U \subseteq V(G)$ with $\lvert U \rvert \leq h-1$ such that for every $v \in V(G)-U$, there exists $[A_v,B_v] \in \E-C$ of order at most $d_{\lvert U \rvert+1}-1$ such that $v \in A_v$.

\noindent{\bf Proof of Claim 2:} Suppose to the contrary that there exist $C \subseteq E(G)$ with $\lvert C \rvert \leq \xi$ and $U \subseteq V(G)$ with $\lvert U \rvert \leq h-1$ such that for every $v \in V(G)-U$, there exists $[A_v,B_v] \in \E-C$ of order at most $d_{\lvert U \rvert+1}-1$ such that $v \in A_v$.
That is, there exists a family $\D' \subseteq \E-C$ of edge-cuts of $G-C$ of order at most $d_{\lvert U \rvert+1}-1$ such that for every $v \in V(G)-U$, there exists $[A,B] \in \D'$ such that $v \in A$.
In particular, $V(G)-U \subseteq \bigcup_{[A,B] \in \D'}A$.
By Lemma \ref{making cross-free}, there exists a cross-free family $\D \subseteq \E-C$ of edge-cuts of $G-C$ of order at most $d_{\lvert U \rvert+1}-1$ such that $\bigcup_{[A,B] \in \D}A=\bigcup_{[A,B] \in \D'}A \supseteq V(G)-U$.
Hence for every $v \in V(G)-U$, there exists $[A,B] \in \D$ with $v \in A$, a contradiction.
$\Box$

\smallskip

For each $i \in [d]$, define $U_i$ to be a subset of $V(G)$ and define $\Se_i^*$ to be a collection of subsets of $E(G)$ such that $U_i$ and $\Se_i^*$ satisfy the following properties.
	\begin{itemize}
		\item[(i)] For every $S \in \Se_i^*$, $S$ consists of $d-i+1$ edges of $G$ with a common end $v_S \not \in \bigcup_{j=1}^{i-1} U_j$, and $S$ is disjoint from $S'$ for every $S' \in \bigcup_{j=1}^{i-1}\Se_j^*$.
		\item[(ii)] For every pair of distinct sets $S,S' \in \Se_i^*$, we have $S \cap S' = \emptyset $ and $v_S \neq v_{S'}$.
		\item[(iii)] $\bigcup_{j=1}^i\bigcup_{S \in \Se_j^*} S$ is free with respect to $\E$.
		\item[(iv)] Subject to (i)-(iii), $\Se_i^*$ is maximal.
		\item[(v)] $U_i = \{v_S: S \in \Se_i^*\}$.
	\end{itemize}

If there exists $k \in [d]$ such that $\lvert \bigcup_{i=1}^k U_i \rvert < \lvert \{u \in V(H): \deg_H(u) \geq d-k+1\} \rvert$, then define $r$ to the minimum such $k$; if there does not exist such $k$, then $\lvert \bigcup_{i=1}^d U_i \rvert \geq \lvert \{u \in V(H): \deg_H(u) \geq 1\} \rvert = h$ since $H$ has no isolated vertex, and we define $r=d$.

If $\lvert \bigcup_{i=1}^r U_i \rvert \leq h$, then define $U^*=\bigcup_{i=1}^r U_i$; otherwise, $r=d$ and we define $U^*$ to be a set with $\bigcup_{i=1}^{j-1}U_i \subseteq U^* \subseteq \bigcup_{i=1}^jU_i$ with $\lvert U^* \rvert=h$, where $j$ is the minimum such that $\sum_{i=1}^j \lvert U_i \rvert \geq h$.
Define $\Se^* = \{S \in \bigcup_{i=1}^r\Se_i^*: v_S \in U^*\}$. 
Note that $\lvert \Se^* \rvert = \lvert U^* \rvert \leq h$, and members of $\Se^*$ are pairwise disjoint and $\bigcup_{S \in \Se^*}S$ is free with respect to $\E$.
For every $v \in V(G)$, define $\Se_v$ to be the collection of the sets of $d_{\lvert U^* \rvert+1}$ edges of $G-\bigcup_{S \in \Se^*}S$ incident with $v$, where we define $d_i=0$ if $i>h$.
Note that $\Se_v = \emptyset$ if $v$ is incident with less than $d_{\lvert U^* \rvert+1}$ edges.

\noindent{\bf Claim 3:} There exist distinct $v_1,v_2,...,v_{h-\lvert U^* \rvert} \in V(G)-U^*$ and pairwise disjoint sets $X_1,X_2,...,X_{h-\lvert U^* \rvert}$ such that $\bigcup_{S \in \Se^*}S \cup \bigcup_{i=1}^{h-\lvert U^* \rvert}X_i$ is free with respect to $\E$, and for each $i \in [h-\lvert U^* \rvert]$, $X_i \in \Se_{v_i}$.

\noindent{\bf Proof of Claim 3:}
There is nothing to prove if $\lvert U^* \rvert \geq h$.
So we may assume that $\lvert U^* \rvert<h$.

Let $\Se_1=\Se^*$ and let $\Se_2=\bigcup_{v \in V(G)-U^*}\Se_v$.
For $i \in [2]$, let $J_i$ be a set such that we can write $\Se_i = \{X_{i,j}: j \in J_i\}$.
Let $J_1^*=J_1$ and $J_2^*=\emptyset$.
So $\bigcup_{i=1}^2\bigcup_{j \in J_i^*}X_{i,j} = \bigcup_{X \in \Se^*}X$ is free with respect to $\E$, and $X_{i,j} \cap X_{i',j'} =\emptyset$ for every distinct pairs $(i,j),(i',j')$ with $1 \leq i \leq i' \leq 2$, $j \in J_i^*$ and $j' \in J_{i'}^*$.
Let $k_1=\lvert \Se_1 \rvert$ and $k_2=(d-1)(h-1)+1$.
Let $k=dh$.
Note that $k_1 = \lvert U^* \rvert \leq h-1$.
So $k \geq \max\{k_1,k_2\}$ and $\theta \geq 3(2kd)^{d+1}+3d$.
Since every member of $\Se_1$ has size at most $d$ and every member of $\Se_2$ has size $d_{\lvert U^* \rvert+1} \leq d$, by Lemma \ref{colorful edge spider}, either 
	\begin{itemize}
		\item[(i')] there exist $J_1',J_2'$ with $J_i^* \subseteq J_i' \subseteq J_i$ and $\lvert J_i' \rvert=k_i$ for $i \in [2]$ such that $\bigcup_{i=1}^2\bigcup_{j \in J_i'}X_{i,j}$ is free with respect to $\E$, and $X_{i,j} \cap X_{i',j'} = \emptyset$ for every distinct pairs $(i,j),(i',j')$ with $1 \leq i \leq i' \leq 2$, $j \in J_i'$ and $j' \in J_{i'}'$, or
		\item[(ii')] there exist $Z \subseteq E(G)$ with $\lvert Z \rvert \leq (2dk)^{d+1}$ and integer $i^* \in [2]$ with $\lvert J_{i^*}^* \rvert < k_{i^*}$ such that for every $j \in J_{i^*}$, either $X_{i^*,j} \cap Z \neq \emptyset$, or $X_{i^*,j}$ is not free with respect to $\E-Z$.
	\end{itemize}

Suppose that (ii') holds.
Define $C = Z \cup \bigcup_{X \in \Se^*}X$.
Then $\lvert C \rvert \leq \lvert Z \rvert + d\lvert \Se^* \rvert \leq (2dk)^{d+1}+dh \leq \xi$.
Since $\lvert J_1^* \rvert = \lvert J_1 \rvert=k_1$, $i^*=2$.
So every member $X$ of $\Se_1 \cup \Se_2$ disjoint from $C$ belongs to $\Se_2=\Se_{i^*}$ and hence is not free with respect to $\E-Z$ and hence is not free with respect to $\E-C$ by Lemma \ref{subset of free set is free}.
By Claim 2, there exists $v \in V(G)-U^*$ such that there does not exist $[A_v,B_v] \in \E-C$ of order at most $d_{\lvert U^* \rvert+1}-1$ such that $v \in A_v$, for otherwise choosing $C=C$ and $U=U^*$ contradicts Claim 2.
In particular, $v \in V(G)-U^*$ is incident with at least $d_{\lvert U^* \rvert+1}$ edges in $G-C$.
Hence there exists $X \in \Se_v \subseteq \Se_2$ such that every edge in $X$ is incident with $v$ and $X \cap C = \emptyset$.
So $X$ is not free with respect to $\E-C$.
Hence there exist $Y \subseteq X$ and an edge-cut $[A,B] \in \E-(C \cup Y)$ of $G-(C \cup Y)$ of order less than $\lvert X-Y \rvert$ such that every edge in $X-Y$ has every end in $A$.
Since every edge in $X-Y$ is incident with $v$, we have that $v \in A$ and $[A,B] \in \E-C$ is an edge-cut of $G-C$ of order less than $\lvert X \rvert=d_{\lvert U^* \rvert+1}$, a contradiction.

So (i') holds.
Note that $J_1=J_1'=J_1^*$.
If there exist $h-\lvert U^* \rvert$ distinct vertices $v_1,v_2,...,v_{h-\lvert U^* \rvert} \in V(G)-U^*$ and $j_1,j_2,...,j_{h-\lvert U^* \rvert} \in J_2'$ such that $X_{2,j_i} \in \Se_{v_i}$ for each $i \in [h-\lvert U^* \rvert]$, then the claim holds. 

So we may assume that there exist at most $h-\lvert U^* \rvert-1 \leq h-1$ vertices $v_1,v_2,...,v_{h-\lvert U^* \rvert-1}$ in $V(G)-U^*$ such that $\{X_{2,j}: j \in J_2'\} \subseteq \bigcup_{i=1}^{h-\lvert U^* \rvert-1} \Se_{v_i}$.
Since $\lvert J_2' \rvert = (d-1)(h-1)+1$, there exists $j^* \in [h-\lvert U^* \rvert-1]$ such that $\Se_{v_{j^*}}$ contains at least $d$ members of $\{X_{2,j}: j \in J_2'\}$.
Let $W$ be a subset of $\bigcup_{X \in \Se_{v_{j^*}} \cap \{X_{2,j}: j \in J_2'\}}X$ of size $d$.
Note that such a set $W$ exists since members of $\{X_{2,j}: j \in J_2'\}$ are pairwise disjoint and non-empty.
Since $\bigcup_{i=1}^2\bigcup_{j \in J_i'}X_{i,j}$ is free with respect to $\E$, $W \cup (\bigcup_{S \in \Se_1^*}S)$ is free with respect to $\E$.
But $W$ is disjoint from $\bigcup_{S \in \Se_1^*}S$ and consists of $d$ edges incident with $v_{j^*} \not \in U_1$, contradicting the maximality of $\Se_1^*$.
This proves the claim.
$\Box$

\smallskip

\noindent{\bf Claim 4:} $H$ contains a loop, $h \geq 2$ and $d_2 \geq 2$.

\noindent{\bf Proof of Claim 4:}
Let $Y_1,Y_2,...,Y_{\lvert \Se^* \rvert}$ be the members of $\Se^*$ such that $\lvert Y_j \rvert \geq \lvert Y_k \rvert$ for every $1 \leq j \leq k \leq \lvert \Se^* \rvert$.
Let $Y_{\lvert \Se^* \rvert+i}=X_i$ for every $i \in [h-\lvert U^* \rvert]$, where $X_i$ is defined in the statement of Claim 3.

We first show that $\lvert Y_j \rvert \geq d_j$ for every $j \in [h]$.
Suppose to the contrary that there exists $j \in [h]$ with $\lvert Y_j \rvert <d_j$.
It is clear that $\lvert Y_j \rvert \geq d_{\lvert U^* \rvert+1} \geq d_j$ when $\lvert U^* \rvert+1 \leq j \leq h$.
So there exists $i_j \in [r]$ such that $Y_j \in \Se_{i_j}^*$.
Since $Y_j \in \Se_{i_j}^*$, $d-i_j+1=\lvert Y_j \rvert<d_j$, so $d-d_j+1<i_j$.
Since $\lvert Y_j \rvert < d_j$, $\lvert \bigcup_{i=1}^{d-d_j+1} U_i \rvert = \lvert \bigcup_{i=1}^{d-d_j+1} \Se_i^* \rvert \leq j-1 < \lvert \{u \in V(H): \deg_H(u) \geq d-(d-d_j+1)+1\} \rvert$, so $d-d_j+1 \geq r \geq i_j$ by the definition of $r$, a contradiction.

So $\lvert Y_j \rvert \geq d_j$ for every $j \in [h]$.
Hence for every $j \in [h]$, there exists $Y_j' \subseteq Y_j$ with $\lvert Y_j' \rvert = d_j$.
Since $Y_1,Y_2,...,Y_h$ are pairwise disjoint and $\bigcup_{j=1}^h Y_j$ is free with respect to $\E$ by Claim 3, we know $Y_1',Y_2',...,Y_h'$ are pairwise disjoint and $\bigcup_{j=1}^h Y_j'$ is free with respect to $\E$.
If $H$ does not contain a loop, then $G$ contains an $H$-immersion by Lemma \ref{edge spider immersion}, a contradiction.

So $H$ contains a loop.
Since $H$ is not exceptional, $h \geq 2$ and $d_2 \geq 2$.
$\Box$

\smallskip

\noindent{\bf Claim 5:} For every $v \in U^* \cup \{v_i: 1 \leq i \leq h-\lvert U^* \rvert\}$, there exist $Z_v \subseteq E(G)$ with $\lvert Z_v \rvert \leq (3hd^2)^{d+1} + dh$ and $[A_v,B_v] \in \E-Z_v$ of order at most one such that $v \in A_v$.

\noindent{\bf Proof of Claim 5:}
Let $v$ be a vertex in $U^* \cup \{v_i: 1 \leq i \leq h-\lvert U^* \rvert\}$.
Let $\Se_1 = \Se^* \cup \{X_i: 1 \leq i \leq h-\lvert U^* \rvert\}$.
Let $\Se_2 = \{W \subseteq E(G): W$ consists of two edges of $G$ incident with $v\}$.
Let $\Se_3 = \{W \subseteq E(G): W$ consists of two edges of $G$ sharing at least one common end $u \in V(G)-\{v\}\}$.

For $i \in [3]$, let $J_i$ be a set such that $\Se_i$ can be written as $\{Y_{i,j}: j \in J_i\}$.
Let $J_1^*=J_1$, $J_2^*=\emptyset$, and $J_3^*=\emptyset$.
So $\bigcup_{i=1}^3 \bigcup_{j \in J_i^*} Y_{i,j}$ is free with respect to $\E$.
Let $k_1 = \lvert \Se_1 \rvert$.
Let $k_2 = dh$ and let $k_3=dh$.
So $\max\{k_1,k_2,k_3\} \leq hd$.
Note that for every $S \in \bigcup_{i=1}^3\Se_i$, $\lvert S \rvert \leq \max\{d,2\} \leq d$ since $d_2 \geq 2$.
Since $\theta \geq 3(hd \cdot 3 \cdot d)^{d+1}+3d$, by Lemma \ref{colorful edge spider}, either 
	\begin{itemize}
		\item[(i')] there exist $J_1',J_2',J_3'$ with $J_i^* \subseteq J_i' \subseteq J_i$ and $\lvert J_i' \rvert = k_i$ for each $i \in [3]$ such that $\bigcup_{i=1}^3\bigcup_{j \in J_i'}Y_{i,j}$ is free with respect to $\E$, and $X_{i,j} \cap X_{i',j'}=\emptyset$ for all distinct pairs $(i,j),(i',j')$ with $1 \leq i \leq i' \leq 3$, $j \in J_i'$ and $j' \in J_{i'}'$, or
		\item[(ii')] there exists $Z_v' \subseteq E(G)$ with $\lvert Z_v' \rvert \leq (3hd^2)^{d+1}$ and $i^* \in [3]$ with $\lvert J_{i^*}^* \rvert <k_{i^*}$ such that for every $j \in J_{i^*}$, either $Y_{i^*,j} \cap Z_v' \neq \emptyset$ or $Y_{i^*,j}$ is not free with respect to $\E-Z_v'$.
	\end{itemize}

Suppose that (i') holds.
Let $\Se_2'=\{Y_{2,j}: j \in J_2'\}$ and $\Se_3' = \{Y_{3,j}: j \in J_3'\}$.
Let $X_0=\bigcup_{Y \in \Se_2' \cup \Se_3'}Y$.
Let $D_v = \Se_3'$.
Since $\lvert J_2' \rvert = hd$, there exists a partition $(D_u: u \in (U^* \cup \{v_i: 1 \leq i \leq h-\lvert U^* \rvert\})-\{v\})$ of $\Se_2'$ into subsets of size at least $d$.
Note that each member of $D_v$ consists of two edges incident with a vertex in $V(G)-\{v\}$, and for every $u \in (U^* \cup \{v_i: 1 \leq i \leq h-\lvert U^* \rvert\})-\{v\}$, each member of $D_u$ consists of two edges incident with $v \in V(G)-\{u\}$.
Hence $G$ contains an $H$-immersion by Lemma \ref{edge spider immersion}, a contradiction.

Therefore, (ii') holds.
Since $k_1=\lvert J_1^* \rvert$, $i^* \in \{2,3\}$.
Let $Z_v = Z_v ' \cup \bigcup_{S \in \Se_1}S$.
So $\lvert Z_v \rvert \leq (3hd^2)^{d+1} + dh \leq \xi$.

If $i^*=2$, then let $u=v$; otherwise, let $u$ be a vertex in $V(G)-\{v\}$ such that there exists $X \in \Se_{i^*}$ such that $u$ is a common end of all edges in $X$.
If there exists at most one edge of $G-Z_v$ incident with $u$, then there exists $[A_u,B_u] \in \E-Z_v$ of order at most one such that $u \in A_u$.
If there exist at least two edges of $G-Z_v$ incident with $u$, then let $W$ be a set of two edges of $G-Z_v$ incident with $u$, so $W \in \Se_{i^*}$.
Since $W$ is disjoint from $Z_v'$, $W$ is not free with respect to $\E-Z_v$, so there exists $W' \subseteq W$ and $[A_W,B_W] \in \E-(Z_v \cup W')$ of order less than $\lvert W -W' \rvert =2-\lvert W' \rvert$ such that every edge in $W-W'$ has every end in $A_W$.
Since $u$ is an end of any edge of $W$, $u \in A_W$.
Note that $[A_W,B_W] \in \E-Z_v$ has order at most one.

So the claim follows if $i^*=2$.
Hence we may assume that $i^*=3$.
Define $U=\{v\}$ and $C=Z_v$.
Then for every $w \in V(G)-U$, either $w$ is incident with at most one edge in $G-Z_v$, or $w \in V(G)-\{v\}$ is a common end of all edges in $X$ for some $X \in \Se_3$.
But in either case, there exists $[A_w,B_w] \in \E-Z_v$ of order at most $1 \leq d_2-1 = d_{\lvert U \rvert+1}-1$ such that $w \in A_w$, contradicting Claim 2.
$\Box$

\noindent{\bf Claim 6:} There exists $Z_0 \subseteq E(G)$ with $\lvert Z_0 \rvert \leq (2hd^2)^{d+1}+hd$ such that for every $v \in V(G)-(U^* \cup \{v_i: 1 \leq i \leq h-\lvert U^* \rvert\})$, there exists $[A_v,B_v] \in \E-Z_0$ of order at most one such that $v \in A_v$.

\noindent{\bf Proof of Claim 6:}
Let $\Se_1=\Se^* \cup \{X_i: 1 \leq i \leq h-\lvert U^* \rvert\}$.
Let $\Se_2 = \{W \subseteq E(G): W$ consists of two edges of $G$ sharing at least one common end $u \in V(G)-(U^* \cup \{v_i: 1 \leq i \leq h-\lvert U^* \rvert\})\}$.
For $i \in [2]$, let $J_i$ be a set such that $\Se_i$ can be written as $\{Y_{i,j}: j \in J_i\}$.
Let $J_1^*=J_1$ and $J_2^*=\emptyset$.
Let $k_1=\lvert J_1^* \rvert$ and $k_2 = hd$.
Note that $\max\{k_1,k_2\} \leq hd$ and $\lvert S \rvert \leq \max\{d,2\}=d$ for every $S \in \Se_1 \cup \Se_2$, since $d_2 \geq 2$.
Since $\theta \geq 3(hd \cdot 2 \cdot d)^{d+1}+3d$ and $\lvert J_1^* \rvert=k_1$, by Lemma \ref{colorful edge spider}, either
	\begin{itemize}
		\item[(i')] there exist $J_1',J_2'$ with $J_i^* \subseteq J_i' \subseteq J_i$ and with $\lvert J_i' \rvert = k_i$ for each $i \in [2]$ such that the members of $\{Y_{1,j}: j \in J_1'\} \cup \{Y_{2,j}: j \in J_2'\}$ are pairwise disjoint, and $\bigcup_{i=1}^2\bigcup_{j \in J_i'}Y_{i,j}$ is free with respect to $\E$, or
		\item[(ii')] there exists $Z_0' \subseteq E(G)$ with $\lvert Z_0' \rvert \leq (2hd^2)^{d+1}$ such that for every $Y_{2,j} \in \Se_2$, either $Y_{2,j} \cap Z_0' \neq \emptyset$, or $Y_{2,j}$ is not free with respect to $\E-Z_0'$.
	\end{itemize}
If (i') holds, then $G$ contains an $H$-immersion by Lemma \ref{edge spider immersion}.
So (ii') holds.
Let $Z_0 = Z_0' \cup \bigcup_{Y \in \Se_1}Y$.
Then $\lvert Z_0 \rvert \leq (2hd^2)^{d+1}+hd$.
For every $v \in V(G)-(U^* \cup \{v_i: 1 \leq i \leq h-\lvert U^* \rvert\})$, either $v$ is incident with at most one edge of $G-Z_0$, or there exists $W \in \Se_2$ consisting of two edges of $G-Z_0$ incident with $v$ such that $W$ is not free with respect to $\E-Z_0$.
For the former, there exists $[A_v,B_v] \in \E-Z_0$ of order at most one such that $v \in A_v$ and we are done.
So we may assume that there exists $W \in \Se_2$ consisting of two edges of $G-Z_0$ incident with $v$ such that $W$ is not free with respect to $\E-Z_0$.
Since $W$ is disjoint from $Z_0'$, $W$ is not free with respect to $\E-Z_0$, so there exists $W' \subseteq W$ and $[A_W,B_W] \in \E-(Z_0 \cup W')$ of order less than $\lvert W -W' \rvert =2-\lvert W' \rvert$ such that every edge in $W-W'$ has every end in $A_W$.
Since $v$ is an end of any edge of $W$, $v \in A_W$.
Note that $[A_W,B_W] \in \E-Z_0$ has order at most one.
This proves the claim.
$\Box$

Define $C= Z_0 \cup \bigcup_{v \in U^* \cup \{v_i: 1 \leq i \leq h-\lvert U^* \rvert\}}Z_v$.
So $\lvert C \rvert \leq (2hd^2)^{d+1}+hd + h \cdot ((3hd^2)^{d+1} + dh) \leq \xi$.
By Claims 5 and 6, for every $v \in V(G)$, there exists $[A_v,B_v] \in \E-C$ of order at most $1 \leq d_1-1$ such that $v \in A_v$.
It is a contradiction to Claim 2 by choosing $U=\emptyset$.
This proves the theorem.
\end{pf1}

\section{Isolating an immersion}
\label{sec:isolating an immersion}

The main result of this section is Lemma \ref{isolating immersion} which states that if a graph that does not contain many edge-disjoint $H$-immersions has an edge-tangle $\E$ of large order that controls a large complete graph-thorns, then one can delete a bounded number of edges to push all $H$-immersions in the remaining graph into the first entry of an edge-cut belonging to $\E$.
The proof of Lemma \ref{isolating immersion} follows from an induction on the number of components of $H$.
The main difficulty lies at the base case, namely the case that $H$ is connected.
This base case will be proved in Lemma \ref{isolating connected immersion}.

Now we sketch the proof of Lemma \ref{isolating connected immersion}.
Assume that $G$ does not contain $k$ edge-disjoint $H$-immersions.
Then $G$ does not contain a large graph $H'$ with $\lvert V(H') \rvert =\lvert V(H) \rvert$ as an immersion.
So Theorem \ref{excluding immersion} implies that one can delete a bounded number of edges from $G$ such that for every vertex $v$ in the remaining graph not contained in a set $U$ with $\lvert U \rvert \leq \lvert V(H) \rvert-1$, there exists an edge-cut $[A_v,B_v]$ in $\E$ of small order such that $v \in A_v$.
Hence for every $H$-immersion $\Pi$ in the remaining graph, there exists an edge-cut $[A_\Pi,B_\Pi]$ in $\E$ of bounded order such that all its branch vertices are contained in $U \cup A_\Pi$, and $A_\Pi$ contains at least one branch vertex.

Assume that we can further delete a set of edges of bounded size from $G$ to either decrease $\lvert U \rvert$ or decrease the order of $[A_\Pi,B_\Pi]$ for all $H$-immersions $\Pi$.
As $\lvert U \rvert$ is bounded, by repeatedly deleting small sets of edges a bounded number of times, at some point we will keep decreasing the order of $[A_\Pi,B_\Pi]$ for all $H$-immersions $\Pi$, so that eventually $[A_\Pi,B_\Pi]$ will have order 0 for all $H$-immersions $\Pi$.
It implies that $\Pi(H)$ is contained in $G[A_\Pi]$ as $H$ is connected and $A_\Pi$ contains at least one branch vertex.
So Lemma \ref{isolating connected immersion} is proved.

So it suffices to show that we can further delete a set of edges of bounded size to either decrease $\lvert U \rvert$ or decrease the order of $[A_\Pi,B_\Pi]$ for all $H$-immersions $\Pi$.

First we show that if we cannot delete a small number of edges to reduce $\lvert U \rvert$, then one can find a set $S_u$ of many edges incident with $u$ for every vertex $u \in U$, such that $S_{u'}$ and $S_{u''}$ are pairwise disjoint for distinct $u',u'' \in U$, and $\bigcup_{u \in U}S_u$ is free with respect to $\E$.
This is the purpose of Lemma \ref{making large degree free} and can be proved by an application of Menger's theorem.

We may assume that we cannot delete a small set of edges to decrease $\lvert U \rvert$, for otherwise we are done.
So we may assume that those sets $S_u$'s exist.
Then we shall show that we can delete a bounded number of edges to reduce the order of $[A_\Pi,B_\Pi]$ for all $H$-immersions $\Pi$.
This is the purpose of Lemma \ref{isolating connected immersion one step}.
Then Lemma \ref{isolating connected immersion} follows from repeatedly applying Lemma \ref{isolating connected immersion one step}.

Now we sketch the proof of Lemma \ref{isolating connected immersion one step}.
For each $H$-immersion $\Pi=(\pi_V,\pi_E)$, since all branch vertices of $\Pi$ are contained in $A_\Pi \cup U$, for every edge $e \in E(H)$ with $V(\pi_E(e)) \cap U = \emptyset$, $V(\pi_E(e)) \cap A_\Pi \neq \emptyset$.
We say that $G[A_\Pi]$ ``fully realizes'' an edge $e$ of $H$ if $\pi_E(e) \subseteq G[A_\Pi]$; $G[A_\Pi]$ ``partially realizes'' an edge $e$ of $H$ if $\pi_E(e) \not \subseteq G[A_\Pi]$ but $V(\pi_E(e)) \cap A_\Pi \neq \emptyset$.
Since $A_\Pi$ contains at least one branch vertex, at least one edge of $H$ is fully realized or partially realized by $G[A_\Pi]$.
This leads to the notion of ``shell'' defined right above the statement of Lemma \ref{isolating connected immersion one step}, which is a collection of subgraphs of $H$ indicating what the vertices of $H$ whose corresponding branch vertices are contained in $U$ are, what the edges fully realized by $G[A_\Pi]$ are, and what the edges partially realized by $G[A_\Pi]$ are.
Note that for each partially realized edge $e$ of $H$, one can find an edge between $A_\Pi$ and $B_\Pi$ contained in $\pi_E(e)$.

Hence, if we can find many $H$-immersions $\Pi_1,\Pi_2,...$, where $G[A_{\Pi_i}] \cap \Pi_i(H)$ are pairwise edge-disjoint, such that the union of $\bigcup_{u \in U}S_u$ and the set of edges between $A_{\Pi_i}$ and $B_{\Pi_i}$ over all $i$ is free, then we can ``link'' those edges to create $k$ edge-disjoint $H$-immersions by using Lemma \ref{edge linkage} to obtain a contradiction.
So it implies that we cannot find such $H$-immersions.
Then Lemma \ref{colorful edge spider} implies one can delete a bounded number of edges to reduce the order of $[A_\Pi,B_\Pi]$ for all $H$-immersions $\Pi$, so that Lemma \ref{isolating connected immersion one step} is proved.

Now we formally prove all results in this section.

\begin{lemma} \label{making large degree free}
For any positive integers $h,w$, there exists a nonnegative integer $\xi^*=\xi^*(h,w)$ such that the following holds.
Let $\theta,\xi,p$ be positive integers with $\theta \geq \xi+\xi^*+1$.
Let $G$ be a graph and $\E$ an edge-tangle in $G$ of order at least $\theta$.
Assume that there exist $Y_0 \subseteq E(G)$ with $\lvert Y_0 \rvert \leq \xi$, $U_0 \subseteq V(G)$ with $\lvert U_0 \rvert \leq h-1$ and a family $\F_0 \subseteq \E-Y_0$ of edge-cuts of $G-Y_0$ of order less than $p$. 
Then there exist $U \subseteq U_0$, a set $Z$ with $Y_0 \subseteq Z \subseteq E(G)$ with $\lvert Z \rvert \leq \lvert Y_0 \rvert+\xi^*$, a family $\F \subseteq \E-Z$ of edge-cuts of $G-Z$ of order less than $p$ and a collection $\{S_u: u \in U\}$ such that the following hold.
	\begin{enumerate}
		\item $U_0-U \subseteq \bigcup_{[A,B] \in \F}A$ and $\F_0 \subseteq \F$. 
		\item For every $u \in U$, $S_u$ consists of $w$ edges of $G-Z$ incident with $u$.
		\item $S_u \cap S_{u'} = \emptyset$ for distinct $u,u' \in U$.
		\item $\bigcup_{u \in U}S_u$ is free with respect to $\E-Z$.
	\end{enumerate}
\end{lemma}

\begin{pf}
Let $a_0=0$, and for every positive integer $i$, let $a_i=a_{i-1}+2h^2w$.
Define $\xi^*=a_{h-1}$.

Let $r$ be an integer with $0 \leq r \leq h-1$ such that there exist $Y_r \subseteq E(G)$ with $Y_0 \subseteq Y_r$ and $\lvert Y_r \rvert \leq \lvert Y_0 \rvert + a_r$, $U_r \subseteq U_0$ with $\lvert U_r \rvert \leq \lvert U_0 \rvert-r$ and a family $\F_r \subseteq \E-Y_r$ of edge-cuts of $G-Y_r$ of order less than $p$ with $U_0-U_r \subseteq \bigcup_{[A,B] \in \F_r}A$ and $\F_0 \subseteq \F_r$.
Note that such a number $r$ exists as we can choose $r=0$.
We assume that $r$ is as large as possible.

We shall prove that this lemma is true if we take $U=U_r$, $Z=Y_r$ and $\F=\F_r$.
It suffices to prove the existence of a collection $\{S_u: u \in U_r\}$ satisfying Statements 2-4.

If $r \geq \lvert U_0 \rvert$, then $U_r = \emptyset$, so Statements 2-4 hold.
So we may assume that $r \leq \lvert U_0 \rvert-1 \leq h-2$.

\noindent{\bf Claim 1:} There exists a collection $\{S_u: u \in U_r\}$ of pairwise disjoint sets such that for every $u \in U_r$, $S_u$ consists of $w$ edges of $G-Y_r$ incident with $u$.

\noindent{\bf Proof of Claim 1:}
Let $U_r'$ be a minimal subset of $U_r$ such that there exists a collection $\{S_u: u \in U_r-U_r'\}$ of pairwise disjoint sets such that for every $u \in U_r-U_r'$, $S_u$ consists of $w$ edges incident with $u$ of $G-Y_r$ whose every end is in $U_r-U_r'$.
So for every two distinct $u,u' \in U_r'$, there exist at most $2w-1$ edges of $G-Y_r$ between $u,u'$ by the minimality of $U_r'$.
Let $Y$ be the set consists of the non-loop edges whose both ends are in $U_r'$.
Note that $\lvert Y \rvert \leq {\lvert U_r' \rvert \choose 2}(2w-1) \leq (h-1)^2(2w-1)$.
Let $G'=G-(Y_r \cup Y)$.

To prove this claim, it suffices to show that there exists a collection $\{S_u:u \in U_r'\}$ of pairwise disjoint sets such that for every $u \in U_r'$, $S_u$ consists of $w$ edges of $G'$ incident with $u$.

Define $H'$ to be the directed graph such that the following hold.
	\begin{itemize}
		\item $V(H')$ is the disjoint union of a set $Q$ and a set $R$, where $Q$ is a copy of $U_r'$ and $R$ is a copy of $V(G)$.
		For each $u \in U_r'$, we denote the copy of $u$ in $Q$ by $u'$; for each $v \in V(G)$, we denote the copy of $v$ in $R$ by $v'$.
		\item Every edge of $H'$ is from $Q$ to $R$.
		\item For every $u' \in Q$ and $v' \in R$ with $u \neq v$, the number of edges of $H'$ from $u'$ to $v'$ equals the number of edges of $G'$ with ends $u,v$.
		\item For every $u' \in Q$ and $v' \in R$ with $u=v$, the number of edges of $H'$ from $u'$ to $v'$ equals the number of loops of $G'$ incident with $u$.
	\end{itemize}
Note that no two distinct vertices in $U_r'$ are adjacent in $G'$.
So there exists a bijection $g$ between $E(H')$ and the set of edges of $G'$ incident with $U_r'$ such that for every edge $e$ of $E(H')$, the ends of $g(e)$ are exactly the originals of ends of $e$.
Define $H$ to be the directed graph obtained from $H'$ by adding two new vertices $s,t$ and adding $w$ edges from $s$ to $u'$ and $w\lvert U_r' \rvert$ edges from $v'$ to $t$ for each $u' \in Q$ and $v' \in R$.

Assume that there exist $w\lvert U_r' \rvert$ edge-disjoint directed paths $P_1,P_2,...,P_{w\lvert U_r' \rvert}$ in $H$ from $s$ to $t$.
So every edge of $H$ incident with $s$ belongs to $\bigcup_{i=1}^{w\lvert U_r' \rvert}P_i$.
Hence for every $u' \in Q$, there exist $w$ edges in $\bigcup_{i=1}^{w\lvert U_r' \rvert}P_i$ from $u'$ to $R$.
For each $u \in U_r'$, define $S_u = \{g(e): e \in E(\bigcup_{i=1}^{w\lvert U_r' \rvert}P_i), e$ is from $u$ to $R\}$.
Then the collection $\{S_u: u \in U_r'\}$ consists of pairwise disjoint sets, and $S_u$ consists of $w$ edges of $G'$ incident with $u$ for each $u \in U_r'$.
So the claim holds.

Hence we may assume that there do not exist $w\lvert U_r' \rvert$ edge-disjoint directed paths in $H$ from $s$ to $t$.
By Menger's Theorem, there exists $X \subseteq E(H)$ with $\lvert X \rvert < w\lvert U_r' \rvert$ such that there exists no directed path in $H-X$ from $s$ to $t$.
We assume that $\lvert X \rvert$ is minimum.
Since for each $v' \in R$, there exist $w\lvert U_r' \rvert>\lvert X \rvert$ edges from $v'$ to $t$, we know there exists an edge in $H-X$ from $v'$ to $t$.
So $X$ does not contain any edge incident with $t$, for otherwise removing any edge incident with $t$ from $X$ does not create a directed path from $s$ to $t$, contradicting the minimality of $X$.
Let $T$ be the subset of $Q$ consisting of the vertices that can be reached from $s$ by a directed path in $H-X$.
So there exists no directed path in $H-X$ from $T$ to $t$.
Note that $T \neq \emptyset$ since there are more than $\lvert X \rvert$ edges incident with $s$.
Since $X$ does not contain any edge incident with $t$, there exist no directed path in $H-X$ from $T$ to $R$.
Let $X' = \{g(e): e \in E(H-s) \cap X\}$.
Then $X'$ contains all the edges of $G'$ incident with $\{u \in U_r': u' \in T\}$.

Let $[A,B] = [\{u \in U_r': u'\in T\}, V(G)-\{u \in U_r': u' \in T\}]$.
Define $Y_{r+1}=Y_r \cup Y \cup X'$, $U_{r+1}=U_r - \{u \in U_r': u'\in T\}$ and $\F_{r+1} = \F_r \cup \{[A,B]\}$.
Then $\lvert Y_{r+1} \rvert \leq \lvert Y_r \rvert + \lvert Y \rvert + \lvert X' \rvert \leq (\lvert Y_0 \rvert + a_r) + (h-1)^2(2w-1) + w\lvert U_r' \rvert -1 \leq (\lvert Y_0 \rvert + a_r) + (h-1)^2(2w-1) + (w(h-1)-1) \leq \lvert Y_0 \rvert + a_{r+1}$.
And $\lvert U_{r+1} \rvert = \lvert U_r \rvert - \lvert T \rvert \leq \lvert U_0 \rvert - (r+1)$.
Note that $[A,B]$ is an edge-cut of $G-Y_{r+1}$ of order zero.
So $\F_{r+1}$ is a family of edge-cuts of $G-Y_{r+1}$ of order less than $p$ such that $U_0-U_{r+1} = (U_0-U_r) \cup \{u \in U_r': u'\in T\} \subseteq \bigcup_{[A',B'] \in \F_{r+1}}A'$.
Since there are at most $\lvert Y \cup X' \rvert < (h-1)^2(2w-1) + (w(h-1)-1) \leq \theta-\lvert Y_r \rvert$ edges of $G-Y_r$ incident with $A$, $[A,B] \in \E-Y_r$ by (E1) and (E3).
So $[A,B] \in \E-Y_{r+1}$.
This contradicts the maximality of $r$.
$\Box$

\smallskip

Let $\{S_u: u \in U_r\}$ be a collection mentioned in Claim 1.
To prove the claim, it suffices to prove that $\bigcup_{u \in U_r}S_u$ is free with respect to $\E-Y_r$.

Suppose to the contrary that $\bigcup_{u \in U_r}S_u$ is not free with respect to $\E-Y_r$.
Then there exist $X \subseteq \bigcup_{u \in U_r}S_u$ and $[A,B] \in \E-(Y_r \cup X)$ of order less than $\lvert (\bigcup_{u \in U_r}S_u) - X \rvert$ such that every edge in $(\bigcup_{u \in U_r}S_u)-X$ has every end in $A$.
Note that $(\bigcup_{u \in U_r}S_u)-X \neq \emptyset$ since there exists no edge-cut of order less than 0.
Let $X'$ be the union of $X$ and the set of edges of $G-Y_r$ with one end in $A$ and one end in $B$.
So $\lvert X' \rvert \leq \lvert X \rvert + \lvert \bigcup_{u \in U_r}S_u \rvert \leq 2\lvert \bigcup_{u \in U_r}S_u \rvert \leq 2(h-1)w$.

Define $Y_{r+1}=Y_r \cup X'$, $U_{r+1}=U_r-A$, and $\F_{r+1} = \F_r \cup \{[A,B]\}$.
So $\lvert Y_{r+1} \rvert \leq \lvert Y_r \rvert + \lvert X' \rvert \leq \lvert Y_0 \rvert + a_r + 2(h-1)w \leq \lvert Y_0 \rvert + a_{r+1}$.
Note that every edge in $(\bigcup_{u \in U_r}S_u)-X \neq \emptyset$ has every end in $A$, so $U_r \cap A \neq \emptyset$.
Hence $\lvert U_{r+1} \rvert \leq \lvert U_r \rvert-1 \leq \lvert U_0 \rvert - (r+1)$.
Since $[A,B] \in \E-Y_r$, $[A,B] \in \E-Y_{r+1}$ is an edge-cut of $G-Y_{r+1}$ of order 0.
So $\F_r \subseteq \F_{r+1}$ and $U_0-U_{r+1} \subseteq \bigcup_{[A',B'] \in \F_{r+1}}A'$.
This contradicts the maximality of $r$ and proves the lemma.
\end{pf}

\bigskip

Let $G$ be a graph and $S$ a subgraph of $G$.
We define $S_G^+$ to be the graph obtained from $S$ by attaching $\deg_G(v)-\deg_S(v)$ leaves to $v$, for each $v \in V(S)$.
So every vertex in $V(S_G^+)-V(S)$ corresponds to an edge in $E(G)-E(S)$.
Note that if $e$ is an edge in $E(G)-E(S)$ with both ends $u,v$ in $V(S)$, then $e$ contributes two leaves to $S_G^+$, where one is adjacent to $u$ and one is adjacent to $v$.
In particular, if $e \in E(G)-E(S)$ is a loop incident with a vertex $v$ in $S$, then $e$ contributes two leaves adjacent to $v$ in $S_G^+$.

Let $G$ and $H$ be graphs, and let $S,R$ be subgraphs of $G,H$, respectively.
We say that $S_G^+$ {\it realizes} $R_H^+$ if $S_G^+$ contains a $R_H^+$-immersion $(\pi_V,\pi_E)$ such that $\pi_V(V(R_H^+)-V(R)) \subseteq V(S_G^+)-V(S)$ and $\pi_V(V(R)) \subseteq V(S)$.

Let $H$ be a graph.
A {\it shell} of $H$ is a collection of disjoint connected subgraphs of $H$ such that every vertex of $H$ is contained in a member of the collection.
For any $H$-immersion $\Pi=(\pi_V,\pi_E)$ in a graph $G$, we denote the subgraph $\bigcup_{e \in E(H)}\pi_E(e) \cup \bigcup_{v \in V(H)}\pi_V(v)$ of $G$ by $\Pi(H)$.

\begin{lemma} \label{isolating connected immersion one step}
For every connected graph $H$ and for positive integers $k,p,\xi_0'$, there exist integers $\theta^*=\theta^*(H,k,p,\xi_0'),w^*=w^*(H,k,p,\xi_0'),\xi^*=\xi^*(H,k,p,\xi_0')$ such that the following holds.
Assume that $G$ is a graph that does not contain $k$ edge-disjoint $H$-immersions and $\E$ is an edge-tangle in $G$ of order at least $\theta^*$ controlling a $K_w$-thorns for some $w \geq w^*$.
If there exist $U' \subseteq V(G)$ with $\lvert U' \rvert \leq \lvert V(H) \rvert-1$, $Z_0' \subseteq E(G)$ with $\lvert Z_0' \rvert \leq \xi_0'$ and a family $\F' \subseteq \E-Z_0'$ of edge-cuts of $G-Z_0'$ of order less than $p$ such that for every $H$-immersion $L=(\pi_V,\pi_E)$ in $G-Z_0'$, there exists $[A_L',B_L'] \in \F'$ such that $\pi_V(V(H)) \subseteq U' \cup A_L'$, then there exist $U \subseteq U'$, $Z^* \subseteq E(G)$ with $\lvert Z^* \rvert \leq \xi^*$ and a family $\F^* \subseteq \E-Z^*$ of edge-cuts of $G-Z^*$ such that either
	\begin{enumerate}
		\item $U \subset U'$, every member of $\F^*$ has order less than $\lvert V(H) \rvert p$, and for every $H$-immersion $\Pi=(\pi_V,\pi_E)$ in $G-Z^*$, there exists $[A^*_\Pi,B^*_\Pi] \in \F^*$ such that $\pi_V(V(H)) \subseteq U \cup A^*_\Pi$, or
		\item every member of $\F^*$ has order less than $p-1$, and for every $H$-immersion $\Pi=(\pi_V,\pi_E)$ in $G-Z^*$, either there exists $[A,B] \in \E-Z^*$ of order zero such that $\Pi(H) \subseteq G[A]$, or there exists $[A_\Pi^*,B_\Pi^*] \in \F^*$ such that $\pi_V(V(H)) \subseteq U \cup A_\Pi^*$.
	\end{enumerate}
\end{lemma}

\begin{pf}
Let $H$ be a connected graph with degree sequence $(d_1,d_2,...,d_h)$, where $h=\lvert V(H) \rvert$.
Let $k,p,\xi_0'$ be positive integers.
We define the following.
	\begin{itemize}
		\item Let $\xi_0 = \xi_0' + \xi_{\ref{making large degree free}}$, where $\xi_{\ref{making large degree free}}$ is the number $\xi^*(h,kd_1)$ mentioned in Lemma \ref{making large degree free}.
		\item Let $\xi_0''=(khd_1(kh+2)(kd_1+p))^{kd_1+p+1}$.
		\item Define $\xi^*=\xi_0+(2h^2d_1)^h\xi_0''$, $w^* = 6k^2h^2d_1p+\xi^*$ and $\theta^* = \xi^*+w^*+hp$.
	\end{itemize}
Let $G$ be a graph that does not contain $k$ edge-disjoint $H$-immersions, and let $\E$ be an edge-tangle of order at least $\theta^*$ in $G$ controlling a $K_w$-thorns $\alpha$ for some $w \geq w^*$.
Assume there exist $U' \subseteq V(G)$ with $\lvert U' \rvert \leq \lvert V(H) \rvert-1$, $Z_0' \subseteq E(G)$ with $\lvert Z_0' \rvert \leq \xi_0'$ and a family $\F' \subseteq \E-Z_0'$ of edge-cuts of $G-Z_0'$ of order less than $p$ such that for every $H$-immersion $L=(\pi_V,\pi_E)$ in $G-Z_0'$, there exists $[A_L',B_L'] \in \F'$ such that $\pi_V(V(H)) \subseteq U' \cup A_L'$.

By Lemma \ref{making large degree free}, there exist $U \subseteq U'$, a set $Z_0$ with $Z_0' \subseteq Z_0 \subseteq E(G)$ with $\lvert Z_0 \rvert \leq \xi_0$, a family $\F \subseteq \E-Z_0$ of edge-cuts of $G-Z_0$ of order less than $p$ and a collection $\{S_u: u \in U\}$ such that the following hold.
	\begin{itemize}
		\item $U'-U \subseteq \bigcup_{[A,B] \in \F}A$ and $\F' \subseteq \F$.
		\item For every $u \in U$, $S_u$ consists of $kd_1$ edges of $G-Z_0$ incident with $u$.
		\item $S_u \cap S_{u'} = \emptyset$ for distinct $u,u' \in U$.
		\item $\bigcup_{u \in U}S_u$ is free with respect to $\E-Z_0$.
	\end{itemize}
We may assume that $U$ is inclusion-wise minimal subject to the conditions above.
So for every $u \in U$, there exists no $[A,B] \in \E-Z_0$ of order less than $p$ such that $u \in A$, for otherwise, we may add $[A,B] \in \F$ and remove $u$ from $U$.

For any subset $Z$ of $E(G)$, we say an $H$-immersion $\Pi$ in $G-Z$ is {\it active} (with respect to $Z$) if there does not exist $[A,B] \in \E-Z$ of order zero such that $\Pi(H) \subseteq G[A]$.

Suppose that this lemma does not hold.

\noindent{\bf Claim 1:} $G-Z_0$ contains an active $H$-immersion with respect to $Z_0$, and $U=U'$.

\noindent{\bf Proof of Claim 1:}
If $G-Z_0$ contains no active $H$-immersion with respect to $Z_0$, then Statement 2 of this lemma holds by taking $Z^*=Z_0$ and $\F^*=\emptyset$, a contradiction.
So $G-Z_0$ contains an active $H$-immersion with respect to $Z_0$.

Now we suppose that $U \subset U'$.
Since $U'-U \subseteq \bigcup_{[A,B] \in \F}A$, for each $u \in U'-U$, there exists $[A_u,B_u] \in \F$ with $u \in A_u$.
Since the order of $\E-Z_0$ is at least $hp$, we know for every $[A,B] \in \F$, $[A \cup \bigcup_{u \in U'-U}A_u, B \cap \bigcap_{u \in U'-U}B_u]$ is an edge-cut of order less than $\lvert [A,B] \rvert + (h-1)p<hp$ and hence belongs to $\E-Z_0$ by Lemma \ref{easy edge-tangle}.
Since for every $H$-immersion $L=(\pi_V,\pi_E)$ in $G-Z_0 \subseteq G-Z_0'$, there exists $[A_L,B_L] \in \F' \subseteq \F$ such that $\pi_V(V(H)) \subseteq U' \cup A_L$, we know $\pi_V(V(H)) \subseteq U \cup A_L \cup \bigcup_{u \in U'-U}A_u$.
So Statement 1 of this lemma follows if we take $Z^*=Z_0$ and $\F^*=\{[A \cup \bigcup_{u \in U'-U}A_u, B \cap \bigcap_{u \in U'-U}B_u]: [A,B] \in \F\}$, a contradiction.
$\Box$

\smallskip

For any member $S$ of some shell of $H$ and any active $H$-immersion $L=(\pi_V,\pi_E)$ in $G-Z_0$ with respect to $Z_0$, we say that an edge-cut $[A,B]$ of $G-Z_0$ is {\it useful for $L,S$} if the following hold.
	\begin{itemize}
		\item $[A,B] \in \E-Z_0$ and the order of $[A,B]$ is less than $p$.
		\item $\pi_V(V(H)) \subseteq A \cup U$.
		\item $G[A]^+_G$ realizes $S^+_H$.
		\item For every vertex in $A$, there exists a path in $G[A]-Z_0$ from this vertex to an end of an edge between $A$ and $B$. 
	\end{itemize}

\noindent{\bf Claim 2:} For every active $H$-immersion $L=(\pi_V,\pi_E)$ in $G-Z_0$ with respect to $Z_0$, there exist a shell $\P_L$ of $H$ and $[A_L,B_L] \in \E-Z_0$ of order less than $p$ such that $\{v\} \in \P_L$ for every $v \in V(H)$ with $\pi_V(v) \in U$, and $[A_L,B_L]$ is useful for $L,S$ for every member $S$ of $\P_L-\{\{v\}: \pi_V(v) \in U\}$.

\noindent{\bf Proof of Claim 2:}
Let $L=(\pi_V,\pi_E)$ be an active $H$-immersion in $G-Z_0$ with respect to $Z_0$.
Note that $L$ is an $H$-immersion in $G-Z_0'$, so there exists $[A_L,B_L] \in \F' \subseteq \F$ such that $\pi_V(V(H)) \subseteq U' \cup A_L=U \cup A_L$.
Since $\F' \subseteq \E-Z_0$ and every member of $\F'$ is an edge-cut $G-Z_0'$ of order less than $p$, we know $[A_L,B_L] \in \E-Z_0$ is an edge-cut of $G-Z_0$ of order less than $p$.

Let $S'$ be the subgraph of $H$ such that $V(S')=\{v \in V(H): \pi_V(v) \in A_L\}$ and $E(S')=\{e \in E(H): \pi_E(e) \subseteq G[A_L]\}$.
Let $\P_L$ be the shell of $H$ that is the union of the set $\{\{v\}: v \in V(H)-V(S')\}$ and the collection consisting of the components of $S'$.
Since for every $u \in U$, there exists no $[A,B] \in \E-Z_0$ of order less than $p$ such that $u \in A$, we know that $\{\{v\}: v \in V(H)-V(S')\} = \{\{v\}: v \in V(H), \pi_V(v) \in U\}$.
Since $\lvert U \rvert \leq h-1$, $\pi_V(v) \not \in U$ for some $v \in V(H)$, so $S'$ contains at least one vertex.
Then $G[A_L]^+_G$ realizes $S^+_H$ for every member $S$ of $\P_L-\{\{v\}: v \in V(H), \pi_V(v) \in U\}$.

So there exists $[A_L,B_L] \in \E-Z_0$ satisfying the first three conditions of being useful for $L,S$, for every member $S$ of $\P_L-\{\{v\}: v \in V(H), \pi_V(v) \in U\}$.
We further choose such $[A_L,B_L]$ such that the order of $[A_L,B_L]$ is as small as possible, and subject to that, $A_L$ is minimal.
To show that $[A_L,B_L]$ is useful for $L,S$ for every member $S$ of $\P_L-\{\{v\}: v \in V(H), \pi_V(v) \in U\}$, it suffices to show that for every vertex in $A_L$, there exists a path in $G[A_L]-Z_0$ from this vertex to an end of an edge between $A_L$ and $B_L$.

Since $L$ is active, the order of $[A_L,B_L]$ is greater than zero.
Suppose that there exists a vertex in $A_L$ such that there exists no path in $G[A_L]-Z_0$ from this vertex to an edge between $A_L$ and $B_L$.
Then there exists a component $C$ of $G[A_L]-Z_0$ such that there exists no path in $G-Z_0$ from $V(C)$ to any edge between $A_L$ and $B_L$.
We define $[A_L',B_L']=[A_L-V(C),B_L \cup V(C)]$.
Since the order of $[A_L',B_L']$ is the same as the order of $[A_L,B_L]$, we know $[A_L',B_L'] \in \E-Z_0$ by Lemma \ref{easy edge-tangle}.
By the minimality of $A_L$, $V(C)$ contains $\pi_V(v)$ for some $v \in V(H)$.
Since $H$ is connected and $C$ is a component of $G[A_L]-Z_0$, $C$ contains $\pi_E(E(H))$.
So $[V(C), V(G)-V(C)]$ is an edge-cut of $G-Z_0$ of order zero such that $C$ contains $\pi_E(E(H))$ and hence $G[V(C)]_G^+$ realizes $H_H^+$.
Note that it implies that $\pi_V(V(H)) \cap U=\emptyset$, so $\{H\}$ is a shell $\P'$ of $H$ with $\{v\} \in \P'$ for each $v \in V(H)$ with $\pi_V(v) \in U$.
In addition, $[V(C),V(G)-V(C)] \in \E-Z_0$ by Lemma \ref{easy edge-tangle}.
Since the order of $[A_L,B_L]$ is greater than 0, it contradicts the minimality of the order of $[A_L,B_L]$.
$\Box$

\smallskip

For every shell $\P$ of $H$ and every subset $D$ of $\{v \in V(H): \{v\} \in \P\}$ of size at most $\lvert U \rvert$, we define the following.
	\begin{itemize}
		\item Define $H_{\P,D}$ to be the graph obtained from the disjoint union of $k$ copies of $H$ by for each $v \in D$, identifying the $k$ copies of $v$ into a vertex.
			Note that $\lvert V(H_{\P,D}) \rvert = k(\lvert V(H) \rvert-\lvert D \rvert)+\lvert D \rvert$, and for any two (not necessarily distinct) vertices in $D$, if there are $\ell$ edges of $H$ between then, then there are $k\ell$ edges of $H_{\P,D}$ between them.
			Note that $G$ contains no $H_{\P,D}$-immersion, for otherwise $G$ contains $k$ edge-disjoint $H$-immersions.
		\item Define $\Q_{\P,D}$ to be the shell of $H_{\P,D}$ consisting of $\{v\}$, for each $v \in D$, and the members of $\P-\{\{v\}: v \in D\}$ in each copy of $H$.
		\item For each $S \in \Q_{\P,D}$, define $\X_{S}=\{X: X$ is the set of edges between $A_L$ and $B_L$, for some active $H$-immersion $L$ in $G-Z_0$ and some edge-cut $[A_L,B_L]$ of $G-Z_0$ that is useful for $L,S\}$. 
			Note that each member of $\X_S$ has size at most $p$.
	\end{itemize}
Define $\X_0 = \{S_u: u \in U\}$.
Recall that $\bigcup_{X \in \X_0} X$ is free with respect to $\E-Z_0$.
Define $\X_E$ to be the collection of the 2-element subsets of $E(G-Z_0)$ each consisting of two edges having at least one common end.
Define $\X_0^*=\X_0$, $k_0=\lvert U \rvert$, $\X_E^*=\emptyset$, $k_E = khd_1$, $\X_{S}^*=\emptyset$ and $k_S=kh$ for each $S \in \Q_{\P,D}-\{\{v\}: v \in D\}$.
Note that $\lvert \Q_{\P,D} \rvert \leq kh$.

\noindent{\bf Claim 3:} For every shell $\P$ of $H$ and every subset $D$ of $\{v \in V(H): \{v\} \in \P\}$ of size at most $\lvert U \rvert$, there exist $Z_{\P,D} \subseteq E(G)-Z_0$ with $\lvert Z_{\P,D} \rvert \leq \xi_0''$ and $S \in \Q_{\P,D}-\{\{v\}: v\in D\}$ such that for every $X \in \X_{S}$, either $X \cap Z_{\P,D} \neq \emptyset$, or $X$ is not free with respect to $\E-(Z_0 \cup Z_{\P,D})$.

\noindent{\bf Proof of Claim 3:}
Let $\P$ be a shell of $H$ and $D$ a subset of $\{v \in V(H): \{v\} \in \P\}$ of size at most $\lvert U \rvert$.
Notice that $\lvert \X_0^* \rvert = k_0$.
By Lemma \ref{colorful edge spider}, one of the following holds. 
	\begin{itemize}
		\item[(i)] There exist a collection $\X_0'$ of size $k_0$ with $\X_0^* \subseteq \X_0' \subseteq \X_0$, a collection $\X_E'$ of size $k_E=khd_1$ with $\X_E^* \subseteq \X_E' \subseteq \X_E$ and collections $\X_S'$ of size $k_S=kh$ with $\X_S^* \subseteq \X_S' \subseteq \X_S$ for each $S \in \Q_{\P,D}-\{\{v\}: v \in D\}$ such that $\X_0' \cup \X'_E \cup \bigcup_{S \in \Q_{\P,D}-\{\{v\}: v \in D\}} \X_S'$ consists of pairwise disjoint members, and the union of its members is free with respect to $\E-Z_0$.
		\item[(ii)] There exist $Z_{\P,D} \subseteq E(G)-Z_0$ with $\lvert Z_{\P,D} \rvert \leq (khd_1(\lvert \Q_{\P,D} \rvert-\lvert D \rvert+2)(kd_1+p))^{kd_1+p+1} \leq \xi_0''$ and $S \in \Q_{\P,D}-\{\{v\}: v \in D\}$ such that for every $X \in \X_S$, either $X \cap Z_{\P,D} \neq \emptyset$, or $X$ is not free with respect to $\E-(Z_0 \cup Z_{\P,D})$.
		\item[(iii)] There exists $Z_{\P,D}\subseteq E(G)-Z_0$ with $\lvert Z_{\P,D} \rvert \leq (khd_1(\lvert \Q_{\P,D} \rvert-\lvert D \rvert+2)(kd_1+p))^{kd_1+p+1} \leq \xi_0''$ such that every set of two edges of $G-(Z_0 \cup Z_{\P,D})$ sharing at least one common end is not free with respect to $\E-(Z_0 \cup Z_{\P,D})$.
	\end{itemize}

Note that Statement (iii) cannot hold by Lemma \ref{no 2 edges} since $\theta \geq \xi_0+\xi_0''+2$.
To prove this claim, it suffices to show that Statement (i) does not hold.

Suppose to the contrary that Statement (i) holds.
We shall derive a contradiction by showing that $G$ contains $k$ edge-disjoint $H$-immersions.

Let $X$ be the union of the members of $\X_0' \cup \X_E' \cup \bigcup_{S \in \Q_{\P,D}-\{\{v\}: v \in D_\P\}} \X_S'$.
So $X$ is free with respect to $\E-Z_0$, and $\lvert X \rvert \leq hkd_1+2hkd_1+kh \cdot kh \cdot p \leq 4h^2k^2d_1p$.
Since $\alpha$ is a $K_w$-thorns controlled by $\E$, there exists a $K_{w-\xi_0}$-thorns $\alpha'$ in $G-Z_0$ controlled by $\E-Z_0$.
Note that $w-\xi_0 \geq w^*-\xi_0 \geq \frac{3}{2} \lvert X \rvert$.

Suppose that there exist $Y \subseteq X$ and an edge-cut $[A,B]$ of $G-(Z_0 \cup Y)$ of order less than $\lvert X \rvert-\lvert Y \rvert$ such that every edge in $X-Y$ is incident with vertices in $A$ and $A \cap V(\alpha'(t)) = \emptyset$ for some $t \in V(K_{w-\xi_0})$.
By (E1), $[A,B]$ or $[B,A]$ is in $\E-(Z_0 \cup Y)$, and hence $[A,B]$ or $[B,A]$ is in $\E-Z_0$.
Since $X$ is free with respect to $\E-Z_0$, $[A,B] \not \in \E-Z_0$.
Since $\E-Z_0$ controls $\alpha'$, $[B,A] \not \in \E-Z_0$, a contradiction.

Therefore, there do not exist $Y \subseteq X$ and an edge-cut $[A,B]$ of $G-(Z_0 \cup Y)$ of order less than $\lvert X \rvert-\lvert Y \rvert \leq \frac{2}{3}(w^*-\xi_0)-\lvert Y \rvert$ such that every edge in $X-Y$ is incident with vertices in $A$ and $A \cap V(\alpha'(t)) = \emptyset$ for some $t \in V(K_{w-\xi_0})$.

For each $S \in \Q_{\P,D}-\{\{v\}: v \in D\}$ and $X_S \in \X_S'$, 
	\begin{itemize}
		\item define $L_S$ to be an $H$-immersion $(\pi^{(S)}_V,\pi^{(S)}_E)$ in $G-Z_0$ such that $X_S$ is the set of edges between $A_{L_S}$ and $B_{L_S}$, where $[A_{L_S},B_{L_S}]$ is a useful edge-cut of $G-Z_0$ for $L_S,S$, 
		\item let $(\pi^{(S,1)}_V,\pi^{(S,1)}_E)$ be an $S_H^+$-immersion in $G[A_{L_S}]_G^+$ such that $\pi_V(V(S_H^+)-V(S)) \subseteq V(G[A_{L_S}]_G^+)-A_{L_S}$ and $\pi_V(V(S)) \subseteq A_{L_S}$, and 
		\item let $f_{X_S}$ be the injection from $V(S^+)-V(S)$ to $X_S$ such that for every $x \in V(S^+)-V(S)$, $f_{X_S}(x)$ is the edge in $X_S$ contained in $\pi^{(S,1)}_E(e)$, where $e$ is the edge in $S_H^+$ incident with $x$.
	\end{itemize}
Define $\iota_D$ to be an injection from $D$ to $U$, and for every $v \in D$, define $f'_{X_{\{v\}}}$ to be an injection from the set of edges of $H_{\P,D}$ incident with $v$ to the set $S_{\iota_D(v)}$, and define $f_{X_{\{v\}}}$ to be the injection from $V(H_{\P,D}[\{v\}]^+_{H_{\P,D}})-\{v\}$ to $S_{\iota_D(v)}$ such that $f_{X_{\{v\}}}(e)=f'_{X_{\{v\}}}(e')$ for every $e \in V(H_{\P,D}[\{v\}])^+_{H_{\P,D}}-\{v\}$ and edge $e'$ of $H_{\P,D}[\{v\}]^+_{H_{\P,D}}$ (so $e'$ is an edge of $H_{\P,D}$) incident with $e$.
Note that for every $v \in D$, $\{v\}$ is a member of $\Q_{\P,D}$.
So for every $S \in \Q_{\P,D}$, $f_{X_S}$ is defined.
In addition, $k_E \geq \lvert E(H_{\P,D}) \rvert$, so there exists an injection $\iota$ from $E(H_{\P,D})$ to $\X_E'$.

For each edge $e$ of $H_{\P,D}$ not contained in any member of $\Q_{\P,D}$, we define the following.
	\begin{itemize}
		\item Say $e$ has one end in $S_1 \in \Q_{\P,D}$ and one end in $S_2 \in \Q_{\P,D}$.
			Note that $S_1$ and $S_2$ are not necessarily distinct, and $e$ corresponds to a leaf $e_1$ in $S_1^+$ and a leaf $e_2$ in $S_2^+$, where $e_1 \neq e_2$ even if $S_1=S_2$ or $e$ is a loop.
			We define $W_e = \{f_{X_{S_1}}(e_1), f_{X_{S_2}}(e_2)\}$.
		\item Define $W'_e=\iota(e)$.
			Note that $W_e'$ is a member of $\X_E'$. 
		\item Define $\{W_{e,1},W_{e,2}\}$ to be a partition of $W_e \cup W'_e$ into two sets of size two each containing exactly one element in $W_e$.
	\end{itemize}
Let $W$ be the union of $W_{e,1} \cup W_{e,2}$ over all edges $e$ of $H_{\P,D}$ not contained in any member of $\Q_{\P,D}$.
Let $\W=\{W_{e,1},W_{e,2}: e \in E(H_{\P,D})$ not contained in any member of $\Q_{\P,D}\}$.
Note that $W$ is a subset of $X$ and $\W$ is a partition of $W$.
Let $\R=\{\{x\}: x \in X-W\}$, and let $\R^*=\R \cup \W$.
Note that $\R$ is a partition of $X-W$ and $\R^*$ is a partition of $X$.

Recall that there do not exist $Y \subseteq X$ and an edge-cut $[A,B]$ of $G-(Z_0 \cup Y)$ of order less than $\lvert X \rvert-\lvert Y \rvert \leq \frac{2}{3}(w-\xi_0)-\lvert Y \rvert$ such that every edge in $X-Y$ is incident with vertices in $A$ and $A \cap V(\alpha'(t)) = \emptyset$ for some $t \in V(K_{w-\xi_0})$.
So by Lemma \ref{edge linkage} (where the partition mentioned in Lemma \ref{edge linkage} is taken to be $\R^*$), there exists a collection $\{T_x: x \in X-W\} \cup \{T_{e,i}: i \in [2], e \in E(H_{\P,D})$ not contained in any member of $\Q_{\P,D}\}$ of pairwise edge-disjoint connected subgraphs in $G-Z_0$ such that $x \in E(T_x)$ for every $x \in X-W$, and $E(T_{e,i}) \cap X = E(T_{e,i}) \cap W = W_{e,i}$ for each $e$ and $i$.

Let $\Q_{\P,D}' = \Q_{\P,D}-\{\{v\}: v \in D\}$.
Note that $\bigcup_{S \in \Q_{\P,D}', X_S \in \X'_S}X_S$ contains all edges between $\bigcup_{S \in \Q_{\P,D}',X_S \in \X'_S}A_{L_S}$ and $\bigcap_{S \in \Q_{\P,D}',X_S \in \X'_S}B_{L_S}$.
Suppose there exists an edge $x \in X$ whose every end is in $\bigcup_{S \in \Q_{\P,D}', X_S \in \X'_S}A_{L_S}$.
Let $Y=X-\{x\}$.
Hence, $[\bigcup_{S \in \Q_{\P,D}',X_S \in \X'_S}A_{L_S}, \allowbreak \bigcap_{S \in \Q_{\P,D}',X_S \in \X'_S}B_{L_S}]$ is an edge-cut of $G-(Z_0 \cup Y)$ of order $0<\lvert X \rvert - \lvert Y \rvert$ such that every edge in $X-Y$ has every end in $\bigcup_{S \in \Q_{\P,D}',X_S \in \X'_S}A_{L_S}$.
Note that $[\bigcup_{S \in \Q_{\P,D}',X_S \in \X'_S}A_{L_S}, \allowbreak \bigcap_{S \in \Q_{\P,D}',X_S \in \X'_S}B_{L_S}] \in \E-(Z_0 \cup Y)$ by Lemma \ref{easy edge-tangle}.
So $X$ is not free with respect to $\E-Z_0$, a contradiction.

Hence every edge in $X$ has at most one end in $\bigcup_{S \in \Q_{\P,D}',X_S \in \X'_S}A_{L_S}$.
In particular, since $\bigcup_{S \in \Q_{\P,D}', X_S \in \X'_S}X_S$ contains all edges between $\bigcup_{S \in \Q_{\P,D}',X_S \in \X'_S}A_{L_S}$ and $\bigcap_{S \in \Q_{\P,D}',X_S \in \X'_S}B_{L_S}$, every edge in a member of $\X_E'$ has every end in $\bigcap_{S \in \Q_{\P,D}',X_S \in \X'_S}B_{L_S}$.
Therefore, $X$ consists of all edges between $\bigcup_{S \in \Q_{\P,D}',X_S \in \X'_S}A_{L_S}$ and $\bigcap_{S \in \Q_{\P,D}',X_S \in \X'_S}B_{L_S}$, and some edge whose every end is in $\bigcap_{S \in \Q_{\P,D}',X_S \in \X'_S}B_{L_S}$.
It follows that each subgraph $T_{e,i}$ does not contain an edge whose every end is in $\bigcup_{S \in \Q_{\P,D}',X_S \in \X'_S}A_{L_S}$.

Since $X$ is free with respect to $\E-Z_0$ and $X_{S'} \cap X_{S''}=\emptyset$ for distinct $S',S'' \in \Q_{\P,D}'$, we have $A_{L_{S'}} \cap A_{L_{S''}}=\emptyset$ for distinct $S',S'' \in \Q_{\P,D}'$ by Lemma \ref{A disjoint}.
So the subgraphs $T_{e,i}$ together with the intersection of the image of $\pi^{(S)}_E$ and $G[A_{L_S}]$, for each $S \in \Q_{\P,D}'$, define a subgraph of $G-Z_0$ containing an $H_{\P,D}$-immersion $(\pi_V,\pi_E)$ in $G-Z_0$ with $\pi_V(v)=\iota_D(v)$ for every $v \in D$, a contradiction.
This proves the claim.
$\Box$

\smallskip

Let $Z^*$ be the union of $Z_0$ and the sets $Z_{\P,D}$ over all shells $\P$ of $H$ and subsets $D$ of $\{v \in V(H): \{v\} \in \P\}$ of size at most $\lvert U \rvert$ mentioned in Claim 3.
Note that there are at most $h^h(hd_1)^h$ different shells of $H$, and for each shell $\P$ of $H$, there are at most $2^h$ different subsets of $\{v \in V(H): \{v\} \in \P\}$.
So $\lvert Z^* \rvert \leq \xi_0+(2h^2d_1)^h\xi_0'' \leq \xi^*$. 

Note that we may assume that there exists an active $H$-immersion in $G-Z^*$ with respect to $Z^*$, for otherwise Statement 2 of this lemma holds by taking $\F^*=\emptyset$.
Note that every $H$-immersion in $G-Z^*$ is an immersion in $G-Z_0'$.

\noindent{\bf Claim 4:} For every active $H$-immersion $L=(\pi_V,\pi_E)$ in $G-Z^*$ with respect to $Z^*$, there exists $[A^*_L,B^*_L] \in \E-Z^*$ of order less than $p-1$ such that $\pi_V(V(H)) \subseteq U \cup A^*_L$. 

\noindent{\bf Proof of Claim 4:}
Let $L=(\pi_V,\pi_E)$ be an active $H$-immersion in $G-Z^*$ with respect to $Z^*$.
So $L$ is an active $H$-immersion in $G-Z_0$ with respect to $Z_0$.
By Claim 2, there exist a shell $\P_L$ of $H$ and $[A_L,B_L] \in \E-Z_0$ of order less than $p$ such that $\{v\} \in \P_L$ for every $v \in V(H)$ with $\pi_V(v) \in U$, and $[A_L,B_L]$ is useful for $L,S$ for every member $S$ of $\P_L-\{\{v\}: \pi_V(v) \in U\}$. 
Let $X$ be the set of edges of $G-Z_0$ between $A_L$ and $B_L$.
So $\lvert X \rvert \leq p-1$.
Let $D=\{v \in V(H): \pi_V(v) \in U\}$.
So $D$ is a subset of $\{v \in V(H): \{v\} \in \P_\L\}$ of size at most $\lvert U \rvert$.
Since $\lvert U \rvert \leq \lvert V(H) \rvert-1$, $X \in \X_S$ for every member $S$ of $\Q_{\P_L,D}-\{\{v\}: v\in D\}$.

By Claim 3, either $X \cap Z^* \neq \emptyset$ or $X$ is not free with respect to $\E-Z^*$.
If $X \cap Z^* \neq \emptyset$, then $[A_L,B_L]$ is an edge-cut of $G-Z^*$ of order less than its order in $G-Z_0$, so $[A_L,B_L] \in \E-Z^*$ is an edge-cut of $G-Z^*$ of order less than $p-1$ such that $\pi_V(V(H)) \subseteq U \cup A_L$ and we are done.
So we may assume that $X \cap Z^*=\emptyset$.

Hence $X$ is not free with respect to $\E-Z^*$.
So there exist $Y \subseteq X$ and $[A,B] \in \E-(Z^* \cup Y)$ of order less than $\lvert X \rvert - \lvert Y \rvert$ such that every edge in $X-Y$ has every end in $A$.
We assume that the order of $[A,B]$ is as small as possible, and subject to that, $A$ is maximal.

Since $[A,B]$ is an edge-cut of $G-Z^*$ of order less than $\lvert X \rvert \leq p-1$, $[A,B] \in \E-Z^*$.
So we are done if $A_L \subseteq A$.

So we may assume that $A_L \not \subseteq A$.
Let $A^*=A \cup A_L$ and $B^*=B \cap B_L$.
Since every edge of $G-(Z^* \cup Y)$ between $A_L,B_L$ is an edge in $X-Y$, it is not incident with $B$.
So $[A^*,B^*]$ is an edge-cut of $G-(Z^* \cup Y)$ of order at most the order of $[A,B]$ with $A^* \supset A$.
Hence $[A^*,B^*] \in \E-(Z^* \cup Y)$ by Lemma \ref{easy edge-tangle}.
But this contradicts the choice of $[A,B]$.
This proves the claim.
$\Box$

Define $\F^*=\{[A^*_L,B^*_L]: L$ is an active $H$-immersion in $G-Z^*$ with respect to $Z^*\}$, where $[A^*_L,B^*_L]$ is the edge-cut mentioned in Claim 4.
Then Statement 2 of this lemma follows.
\end{pf}

\begin{lemma} \label{isolating connected immersion big step}
For every connected graph $H$ and for positive integers $k,p,\xi$, there exist integers $\theta^*=\theta^*(H,k,p,\xi),w^*=w^*(H,k,p,\xi),\xi^*=\xi^*(H,k,p,\xi),p^*=p^*(H,k,p,\xi)$ such that the following holds.
Assume that $G$ is a graph that does not contain $k$ edge-disjoint $H$-immersions and $\E$ is an edge-tangle of $G$ of order at least $\theta^*$ controlling a $K_w$-thorns for some $w \geq w^*$.
If there exist $U' \subseteq V(G)$ with $\lvert U' \rvert \leq \lvert V(H) \rvert-1$, $Z \subseteq E(G)$ with $\lvert Z \rvert \leq \xi$ and a family $\F' \subseteq \E-Z$ of edge-cuts of $G-Z$ of order less than $p$ such that for every $H$-immersion $L=(\pi_V,\pi_E)$ in $G-Z$, there exists $[A_L',B_L'] \in \F'$ such that $\pi_V(V(H)) \subseteq U' \cup A_L'$, then there exist $U \subseteq U'$, $Z^* \subseteq E(G)$ with $\lvert Z^* \rvert \leq \xi^*$ and a family $\F^* \subseteq \E-Z^*$ of edge-cuts of $G-Z^*$ such that either
	\begin{enumerate}
		\item $U \subset U'$, every member of $\F^*$ has order less than $p^*$, and for every $H$-immersion $\Pi=(\pi_V,\pi_E)$ in $G-Z^*$, there exists $[A^*_\Pi,B^*_\Pi] \in \F^*$ such that $\pi_V(V(H)) \subseteq U \cup A^*_\Pi$, or
		\item for every $H$-immersion $\Pi$ in $G-Z^*$, there exists $[A,B] \in \E-Z^*$ of order zero such that $\Pi(H) \subseteq G[A]$.
	\end{enumerate}
\end{lemma}

\begin{pf}
Let $H$ be a connected graph and $k,p,\xi$ be positive integers.
Let $h=\lvert V(H) \rvert$.
We define the following.
	\begin{itemize}
		\item Let $\xi_0=\xi$, $\theta_0=0$ and $w_0=0$.
		\item For every positive integer $i$ with $1 \leq i \leq p$, define $\theta_i = \theta_{i-1} + \theta_{\ref{isolating connected immersion one step}}(H,k,p-i+1,\xi_{i-1})$, $w_i=w_{i-1}+w_{\ref{isolating connected immersion one step}}(H,k,p-i+1,\xi_{i-1})$, and $\xi_i = \xi_{\ref{isolating connected immersion one step}}(H,k,p-i+1,\xi_{i-1})$, where $\theta_{\ref{isolating connected immersion one step}}, w_{\ref{isolating connected immersion one step}},\xi_{\ref{isolating connected immersion one step}}$ are the numbers $\theta^*,w^*,\xi^*$ mentioned in Lemma \ref{isolating connected immersion one step}, respectively.
		\item Define $\theta^* = \sum_{i=1}^p\theta_i$, $w^*=\sum_{i=1}^pw_i$, $\xi^*=\sum_{i=1}^p\xi_i$ and $p^*=hp$.
	\end{itemize}

Let $G$ be a graph that does not contain $k$ edge-disjoint $H$-immersions, and let $\E$ be an edge-tangle of $G$ of order at least $\theta^*$ controlling a $K_w$-thorns for some $w \geq w^*$.
Assume that there exist $U' \subseteq V(G)$ with $\lvert U' \rvert \leq \lvert V(H) \rvert-1$, $Z \subseteq E(G)$ with $\lvert Z \rvert \leq \xi$ and a family $\F' \subseteq \E-Z$ of edge-cuts of $G-Z$ of order less than $p$ such that for every $H$-immersion $L=(\pi_V,\pi_E)$ in $G-Z$, there exists $[A_L',B_L'] \in \F'$ such that $\pi_V(V(H)) \subseteq U' \cup A_L'$.

Let $Z_0=Z$ and $\F_0=\F'$.
Let $r$ be an integer with $0 \leq r \leq p-1$ such that there exist a set $Z_r \subseteq E(G)$ with $\lvert Z_r \rvert \leq \xi_r$ and a family $\F_r \subseteq \E-Z_r$ of edge-cuts of $G-Z_r$ of order less than $p-r$ such that for every $H$-immersion $\Pi=(\pi_V,\pi_E)$ in $G-Z_r$, there exists $[A,B] \in \F_r$ such that $\pi_V(V(H)) \subseteq U' \cup A$.
Note that such an integer $r$ exists as $r=0$ is a candidate.
We assume that $r$ is as large as possible.

Applying Lemma \ref{isolating connected immersion one step} by taking $(H,k,p,\xi_0',U',Z_0',\F')=(H,k,p-r,\xi_r,U',Z_r,\F_r)$, there exist $U \subseteq U'$, $Z^* \subseteq E(G)$ with $\lvert Z^* \rvert \leq \xi_{r+1}$ and a family $\F^* \subseteq \E-Z^*$ of edge-cuts of $G-Z^*$ such that either 
	\begin{itemize}
		\item[(i)] $U \subset U'$, every member of $\F^*$ has order less than $h \cdot (p-r)$, and for every $H$-immersion $\Pi=(\pi_V,\pi_E)$ in $G-Z^*$, there exists $[A,B] \in \F^*$ such that $\pi_V(V(H)) \subseteq U \cup A$, or
		\item[(ii)] every member of $\F^*$ has order less than $p-r-1$ and for every $H$-immersion $\Pi=(\pi_V,\pi_E)$ in $G-Z^*$, either there exists $[A,B] \in \E-Z^*$ of order zero such that $\Pi(H) \subseteq G[A]$, or there exists $[A^*_\Pi,B^*_\Pi] \in \F^*$ such that $\pi_V(V(H)) \subseteq U \cup A^*_\Pi$.
	\end{itemize}

If (i) holds, then since $\xi^* \geq \xi_{r+1}$ and $p^* \geq h(p-r)$, Statement 1 of this lemma holds.
So we may assume that (ii) holds.

Assume that $r=p-1$.
Since there exists no edge-cut of order less than zero, $\F^*=\emptyset$.
So for every $H$-immersion $\Pi$ in $G-Z^*$, there exists $[A,B] \in \E-Z^*$ of order zero such that $\Pi(H) \subseteq G[A]$.
Hence Statement 2 of this lemma holds.

So we may assume that $r \leq p-2$.
Define $Z_{r+1}=Z^*$, and define $\F_{r+1}=\F^* \cup \{[A,B] \in \E-Z^*: [A,B]$ is an edge-cut of $G-Z^*$ of order zero$\}$.
Since $r \leq p-2$, every member of $\F_{r+1}$ has order less than $p-r-1$.
This contradicts the maximality of $r$ and completes the proof.
\end{pf}

\begin{lemma} \label{isolating connected immersion two vertices}
For every connected graph $H$ on at least two vertices and for every positive integer $k$, there exist integers $\theta=\theta(H,k),w=w(H,k),\xi=\xi(H,k)$ with $\theta > w+\xi$ such that the following holds.
If $G$ is a graph that does not contain $k$ edge-disjoint $H$-immersions and $\E$ is an edge-tangle in $G$ of order at least $\theta$ controlling a $K_{w'}$-thorns for some $w' \geq w$, then there exist $Z \subseteq E(G)$ with $\lvert Z \rvert \leq \xi$ and $[A,B] \in \E-Z$ of order zero such that $G[A]$ contains all $H$-immersions in $G-Z$.
\end{lemma}

\begin{pf}
Let $H$ be a connected graph with degree sequence $(d_1,d_2,...,d_h)$, where $h=\lvert V(H) \rvert \geq 2$.
Let $k$ be a positive integer.
We define the following.
	\begin{itemize}
		\item Let $\theta_0=\theta_{\ref{excluding immersion}}(kd_1,h)$ and let $\xi_0=\xi_{\ref{excluding immersion}}(kd_1,h)$, where $\theta_{\ref{excluding immersion}}$ and $\xi_{\ref{excluding immersion}}$ are the numbers $\theta$ and $\xi$ mentioned in Theorem \ref{excluding immersion}, respectively.
		\item Let $p_0=kd_1h$ and $w_0=3kd_1h$.
		\item For every positive integer $i$, let $\theta_i = \theta_{i-1} + \theta_{\ref{isolating connected immersion big step}}(H,k,p_{i-1},\xi_{i-1})$, $w_i = w_{\ref{isolating connected immersion big step}}(H,k,p_{i-1},\xi_{i-1})$, $\xi_i = \xi_{\ref{isolating connected immersion big step}}(H,k,p_{i-1},\xi_{i-1})$ and $p_i = p_{\ref{isolating connected immersion big step}}(H,k,p_{i-1},\xi_{i-1})$, where $\theta_{\ref{isolating connected immersion big step}}, w_{\ref{isolating connected immersion big step}}, \xi_{\ref{isolating connected immersion big step}}$ and $p_{\ref{isolating connected immersion big step}}$ are the numbers $\theta^*,w^*,\xi^*$ and $p^*$ mentioned in Lemma \ref{isolating connected immersion big step}, respectively.
		\item Define $\xi = \sum_{i=0}^{h}\xi_i$, $w=w_{h}+\xi$ and $\theta=\theta_{h} +w + \xi + \sum_{i=0}^{h}p_i$.
	\end{itemize}
Let $G$ be a graph that does not contain $k$ edge-disjoint $H$-immersions, and let $\E$ be an edge-tangle of order at least $\theta$ in $G$ controlling a $K_{w'}$-thorns $\alpha$ for some $w' \geq w$.
We may assume that $k \geq 2$, for otherwise the lemma holds by choosing $Z=\emptyset$.

Define $H_k$ to be the graph obtained from $H$ by duplicating each edge $k$ times.
Note that $H_k$ is a graph on $h$ vertices with maximum degree $kd_1$.
Since $h \geq 2$ and $H$ is connected, $d_2 \geq 1$.
So $H_k$ contains at least two vertices of degree at least $k \geq 2$ and hence is not an exceptional graph.
Since $G$ does not contain $k$ edge-disjoint $H$-immersions, $G$ does not contain an $H_k$-immersion.
By Theorem \ref{excluding immersion}, there exist $Z_0 \subseteq E(G)$ with $\lvert Z_0 \rvert \leq \xi_0$, $U_0 \subseteq V(G)$ with $\lvert U_0 \rvert \leq h-1$ and a family $\F'_0 \subseteq \E-Z_0$ of edge-cuts such that for each $v \in V(G)-U_0$, there exists $[A_v,B_v] \in \F'_0$ of order less than $kd_1$ with $v \in A_v$.

For every $H$-immersion $L=(\pi_V,\pi_E)$ in $G-Z_0$, define $[A_L^*,B_L^*] = [\bigcup_{v \in V(H), \pi_V(v) \not \in U_0}A_{\pi_V(v)}, \allowbreak \bigcap_{v \in V(H), \pi_V(v) \not \in U_0}B_{\pi_V(v)}]$.
Note that each $[A_L^*,B_L^*]$ has order less than $kd_1h$ and hence belongs to $\E-Z_0$ by Lemma \ref{easy edge-tangle}.
Define $\F_0=\{[A^*_L,B^*_L]: L$ is an $H$-immersion in $G-Z_0\}$ to be a collection of edge-cuts of $G-Z_0$.
Note that for every $H$-immersion $L=(\pi_V,\pi_E)$ in $G-Z_0$, $\pi_V(V(H)) \subseteq U_0 \cup A_L^*$.

Let $r$ be an integer with $0 \leq r \leq h-1$ such that there exist $U_r \subseteq V(G)$ with $\lvert U_r \rvert \leq h-1-r$, $Z_r \subseteq E(G)$ with $\lvert Z_r \rvert \leq \xi_r$ and a family $\F_r \subseteq \E-Z_r$ of edge-cuts of $G-Z_r$ of order less than $p_r$ such that for every $H$-immersion $L=(\pi_V,\pi_E)$ in $G-Z_r$, there exists $[A_L',B_L'] \in \F_r$ such that $\pi_V(V(H)) \subseteq U_r \cup A_L'$.
Note that such an integer $r$ exists since $r=0$ is a candidate.
We assume that $r$ is as large as possible.

Apply Lemma \ref{isolating connected immersion big step} by taking $(H,k,p,\xi,U',Z,\F')=(H,k,p_r,\xi_r,U_r,Z_r,\F_r)$, there exist $U \subseteq U_r$, $Z^* \subseteq E(G)$ with $\lvert Z^* \rvert \leq \xi_{r+1}$ and a family $\F^* \subseteq \E-Z^*$ of edge-cuts of $G-Z^*$ such that either
	\begin{itemize}
		\item[(i)] $U \subset U_r$, every member of $\F^*$ has order less than $p_{r+1}$, and for every $H$-immersion $\Pi=(\pi_V,\pi_E)$ in $G-Z^*$, there exists $[A^*_\Pi,B^*_\Pi] \in \F^*$ such that $\pi_V(V(H)) \subseteq U \cup A^*_\Pi$, or
		\item[(ii)] for every $H$-immersion $\Pi$ in $G-Z^*$, there exists $[A_\Pi^*,B_\Pi^*] \in \E-Z^*$ of order zero such that $\Pi(H) \subseteq G[A_\Pi^*]$.
	\end{itemize}

We first suppose that (i) holds.
If $r=h-1$, then $U_r=\emptyset$, so (i) does not hold, a contradiction.
So $r \leq h-2$.
But it is a contradiction to the maximality of $r$ by defining $U_{r+1}=U$, $Z_{r+1}=Z^*$ and $\F_{r+1}=\F^*$.

Hence (i) does not hold.
So (ii) holds.
Let $[A,B]=[\bigcup_L A_L^*, \bigcap_L B_L^*]$, where the union and intersection are over all $H$-immersions $L$ in $G-Z^*$ and each $[A_L^*,B_L^*]$ is the member of $\E-Z^*$ mentioned in (ii).
Note that $[A,B]$ has order zero.
So $[A,B] \in \E-Z^*$ by Lemma \ref{easy edge-tangle}.
Note that for every $H$-immersion $\Pi$ in $G-Z^*$, $\Pi(H) \subseteq G[A]$.
Then this lemma follows by taking $Z=Z^*$.
\end{pf}

\bigskip

Now we drop the requirement of the number of vertices of $H$ from Lemma \ref{isolating connected immersion two vertices}.

\begin{lemma} \label{isolating connected immersion}
	For every connected graph $H$ and for every positive integer $k$, there exist integers $\theta=\theta(H,k),w=w(H,k),\xi=\xi(H,k)$ with $\theta > w+\xi$ such that the following holds.
	If $G$ is a graph that does not contain $k$ edge-disjoint $H$-immersions and $\E$ is an edge-tangle in $G$ of order at least $\theta$ controlling a $K_{w'}$-thorns for some $w' \geq w$, then there exist $Z \subseteq E(G)$ with $\lvert Z \rvert \leq \xi$ and $[A,B] \in \E-Z$ of order zero such that $G[A]$ contains all $H$-immersions in $G-Z$.
\end{lemma}

\begin{pf}
Let $H$ be a connected graph and let $k$ be a positive integer.
By Lemma \ref{isolating connected immersion two vertices}, we may assume $\lvert V(H) \rvert=1$.
Note that we may assume $\lvert E(H) \rvert \geq 1$, for otherwise every graph on at least one vertex contains arbitrarily many edge-disjoint $H$-immersions.
Hence $H$ is a one-vertex graph with at least one loop.
Let $H'$ be the graph obtained from $H$ by subdividing one edge of $H$ once.

Define $\xi=\xi_{\ref{isolating connected immersion two vertices}}(H',k)+(k-1)(k\lvert E(H) \rvert-1)$, $w=w_{\ref{isolating connected immersion two vertices}}(H',k)$ and $\theta=\theta_{\ref{isolating connected immersion two vertices}}(H',k)+w+\xi$, where $\xi_{\ref{isolating connected immersion two vertices}},w_{\ref{isolating connected immersion two vertices}},\theta_{\ref{isolating connected immersion two vertices}}$ are the numbers $\xi,w,\theta$ mentioned in Lemma \ref{isolating connected immersion two vertices}, respectively.

Let $G$ be a graph that does not contain $k$ edge-disjoint $H$-immersions and $\E$ an edge-tangle in $G$ of order at least $\theta$ controlling a $K_{w'}$-thorns for some $w' \geq w$.
Since $G$ does not contain $k$ edge-disjoint $H$-immersions, no vertex of $G$ is incident with at least $k\lvert E(H) \rvert$ loops, and there are at most $k-1$ vertices of $G$ incident with at least $\lvert E(H) \rvert$ loops.
Hence there exists $Z_0 \subseteq E(G)$ with $\lvert Z_0 \rvert \leq (k-1)(k\lvert E(H) \rvert-1)$ such that no vertex in $G-Z_0$ is incident with at least $\lvert E(H) \rvert$ loops in $G-Z_0$.

Since $G$ does not contain $k$ edge-disjoint $H$-immersions, $G$ does not contain $k$ edge-disjoint $H'$-immersions.
By Lemma \ref{isolating connected immersion two vertices}, there exist $Z' \subseteq E(G)$ with $\lvert Z' \rvert \leq \xi_{\ref{isolating connected immersion two vertices}}(H',k)$ and $[A,B] \in \E-Z'$ of order zero such that $G[A]$ contains all $H'$-immersions in $G-Z'$.
Let $Z=Z_0 \cup Z'$.
So $\lvert Z \rvert \leq \xi$ and $[A,B] \in \E-Z$ is an edge-cut of $G-Z$ of order zero.

Suppose that there exists an $H$-immersion $\Pi$ in $G-Z$ such that $\Pi(H) \not \subseteq G[A]$.
Since $H$ is connected, $\Pi(H) \subseteq G[B]$.
So $\Pi(H)$ does not contain an $H'$-immersion.
Hence $\Pi(H)$ consists of one vertex and $\lvert E(H) \rvert$ loops.
But no vertex in $G-Z$ is incident with at least $\lvert E(H) \rvert$ loops in $G-Z$, a contradiction.
This proves the lemma.
\end{pf}

\bigskip

The following is the main result of this section, which says that the connectivity of the graph $H$ in Lemma \ref{isolating connected immersion} can be dropped.

\begin{lemma} \label{isolating immersion}
For every graph $H$ and for every positive integer $k$, there exist integers $\theta=\theta(H,k),w=w(H,k),\xi=\xi(H,k)$ with $\theta > w + \xi$ such that the following holds.
If $G$ is a graph that does not contain $k$ edge-disjoint $H$-immersions and $\E$ is an edge-tangle in $G$ of order at least $\theta$ controlling a $K_{w'}$-thorns for some $w' \geq w$, then there exist $Z \subseteq E(G)$ with $\lvert Z \rvert \leq \xi$ and $[A,B] \in \E-Z$ of order zero in $G-Z$ such that $G[B]-Z$ contains no $H$-immersion.
\end{lemma}

\begin{pf}
Let $H$ be a graph and let $k$ be a positive integer.
Let $p$ be the number of components of $H$.
We shall prove this lemma by induction on $p$.
When $p=1$, this lemma holds by taking $\theta,w,\xi$ to be the numbers $\theta_{\ref{isolating connected immersion}}(H,k),w_{\ref{isolating connected immersion}}(H,k),\xi_{\ref{isolating connected immersion}}(H,k)$, where $\theta_{\ref{isolating connected immersion}},w_{\ref{isolating connected immersion}},\xi_{\ref{isolating connected immersion}}$ are the integers $\theta,w,\xi$ mentioned in Lemma \ref{isolating connected immersion}.
So we may assume that $p \geq 2$ and the lemma holds for every graph with less than $p$ components.

Define $\F_1$ to be the set of graphs that can be obtained from $H$ by adding an edge between different components.
Define $\F_2$ to be the set of graphs that can be obtained from $H$ by subdividing an edge and adding an edge between this new vertex and another component of $H$.
Define $\F_3$ to be the set of graphs that can be obtained from $H$ by subdividing two edges in different components and either adding an edge between those two new vertices or identifying the two new vertices.
Define $\F_4$ to be the set of graphs that can be obtained from $H$ by subdividing an edge and identify this new vertex with a vertex in another component.
Let $\F = \F_1 \cup \F_2 \cup \F_3 \cup \F_4$.
Note that $\lvert \F \rvert \leq \lvert V(H) \rvert^2 + \lvert E(H) \rvert \lvert V(H) \rvert + 2\lvert E(H) \rvert^2 + \lvert E(H) \rvert \lvert V(H) \rvert \leq 5(\lvert V(H) \rvert^2+\lvert E(H) \rvert^2)$.
Since every graph in $\F$ contains less than $p$ components, by the induction hypothesis, for every graph $F \in \F$, there exist integers $\theta(F,k), w(F,k),\xi(F,k)$ such that the lemma holds.

Define $w(H,k)=\sum_{F \in \F}w(F,k)$, $\xi(H,k) = \sum_{F \in \F}\xi(F,k)$, and $\theta(H,k) = w(H,k) + \xi(H,k) + \sum_{F \in \F}\theta(F,k)$.
We shall prove that the numbers $\theta(H,k),w(H,k)$ and $\xi(H,k)$ satisfy the lemma.
Let $\theta=\theta(H,k)$, $w=w(H,k)$ and $\xi=\xi(H,k)$.

Let $G$ be a graph that does not contain $k$ edge-disjoint $H$-immersions and $\E$ an edge-tangle in $G$ of order at least $\theta$ controlling a $K_{w'}$-thorns for some $w' \geq w$.
Note that for every $F \in \F$ and every $Z \subseteq E(G)$, any subgraph of $G-Z$ containing an $F$-immersion contains an $H$-immersion.
So for every $F \in \F$, $G$ does not contain $k$ edge-disjoint $F$-immersions.
By the induction hypothesis, for every $F \in \F$, there exist $Z_F \in E(G)$ with $\lvert Z_F \rvert \leq \xi(F,k)$ and $[A_F,B_F] \in \E-Z_F$ of order zero in $G-Z_F$ such that $G[B_F]-Z_F$ contains no $F$-immersion.
Define $Z=\bigcup_{F \in \F} Z_F$ and $[C,D]=[\bigcup_{F \in \F} A_F, \bigcap_{F \in \F} B_F]$.
So $[C,D]$ has order zero in $G-Z$ and $G[D]-Z$ contains no $F$-immersion for each $F \in \F$.
Since $[C,D]$ has order zero, $[C,D] \in \E-Z$ by Lemma \ref{easy edge-tangle}.
Define $[A,B]$ to be the edge-cut of $G-Z$ of order zero such that $[A,B] \in \E-Z$ and $C \subseteq A$, and subject to those, $A$ is maximal.
Then the maximality of $A$ implies that $G[B]-Z$ is connected by Lemma \ref{easy edge-tangle}.

Suppose that $G[B]-Z$ contains an $H$-immersion.
Then for each $i \in [p]$, $G[B]-Z$ contains an $H_i$-immersion $\Pi_i=(\pi_V^{(i)},\pi_E^{(i)})$, where $H_i$ is the $i$-th component of $H$, such that the images of $\pi_V^{(1)},...,\pi_V^{(p)}$ are pairwise disjoint and the images of $\pi_E^{(1)},...,\pi_E^{(p)}$ are pairwise edge-disjoint.
If there exist distinct $i,j \in [p]$ such that $V(\Pi_i(H_i)) \cap V(\Pi_j(H_j)) \neq \emptyset$, then $G[B]-Z$ contains an $F'$-immersion for some $F' \in \F_3 \cup \F_4 \subseteq \F$, a contradiction.
Since $G[B]-Z$ is connected, there exist distinct $i,j \in [p]$ and a path $P$ in $G[B]-Z$ of length at least one from $V(\Pi_i(H_i))$ to $V(\Pi_j(H_j))$ internally disjoint from $\bigcup_{\ell=1}^pV(\Pi_\ell(H_\ell))$.
But it implies that $G[B]-Z$ contains an $F'$-immersion for some $F' \in \F_1 \cup \F_2 \cup \F_3 \subseteq \F$, a contradiction.
Therefore, $G[B]-Z$ contains no $H$-immersion.
\end{pf}

\section{Edge-tangles in 4-edge-connected graphs}
\label{sec:4-edge-connected edge-tangle}
A $m \times n$ {\it grid} is the graph with vertex-set $\{1,2,...,n\} \times \{1,2,...,m\}$ and two vertices $(x,y),(x',y')$ are adjacent if and only if $\lvert x-x' \rvert + \lvert y-y' \rvert=1$.
For every $i \in [m]$, the {\it $i$-th row} of a $m \times n$ grid is the subgraph induced by $\{(x,i): x \in [n]\}$.
For every $j \in [n]$, the {\it $j$-th column} of a $m \times n$ grid is the subgraph induced by $\{(j,y): y \in [m]\}$.

For every positive integer $r$, the {\it diagonal vertices} of the $r \times 2r$ wall are the vertices $\{(2i-1,i): 1 \leq i \leq r\}$.

\begin{lemma}[{\cite[Theorem (1.5)]{cdks}}] \label{wall-subdiv grid-immersion}
For every $g >1$, there exists $b \geq 0$ such that the following holds.
Let $(\pi_V,\pi_E)$ be a wall-subdivision in a graph $G$, and let $S$ be a subset of the image of $\pi_V$ of the diagonal vertices of the wall such that for every pair of distinct vertices $x,y$ in $S$, $G$ contains four edge-disjoint paths from $x$ to $y$.
If $\lvert S \rvert \geq b$, then there exists a $g \times g$ grid-immersion $(\pi_V',\pi_E')$ in $G$ such that the image of $\pi_V'$ is contained in $S$.
\end{lemma}

In fact, in \cite{cdks}, Chudnovsky et al.\ proved that the grid-immersion $(\pi_V',\pi_E')$ mentioned in Lemma \ref{wall-subdiv grid-immersion} is a ``strong immersion.''
We omit the definition of strong immersions as we do not need this notion in the rest of the paper.
But we remark that every $H$-strong immersion is an $H$-immersion.
On the other hand, the following lemma shows that if we do not require $(\pi_V',\pi_E')$ to be a strong immersion, we can strengthen Lemma \ref{wall-subdiv grid-immersion} by showing that the mentioned wall-subdivision $(\pi_V,\pi_E)$ can be replaced by a wall-immersion.

\begin{lemma} \label{wall-immersion grid-immersion}
For every $g >1$, there exists $b \geq 0$ such that the following holds.
Let $(\pi_V,\pi_E)$ be a wall-immersion in a graph $G$, and let $S$ be a subset of the image of $\pi_V$ of the diagonal vertices of the wall such that for every pair of distinct vertices $x,y$ in $S$, $G$ contains four edge-disjoint paths from $x$ to $y$.
If $\lvert S \rvert \geq b$, then there exists a $g \times g$ grid-immersion $(\pi_V',\pi_E')$ in $G$ such that the image of $\pi_V'$ is contained in $S$.
\end{lemma}

\begin{pf}
Let $g$ be an integer with $g>1$.
Define $b$ to be the number $b$ mentioned in Lemma \ref{wall-subdiv grid-immersion}.

Let $G$ be a graph and let $W$ be a wall such that $(\pi_V,\pi_E)$ is a $W$-immersion in $G$.
Let $S$ be a subset of the image of $\pi_V$ of the diagonal vertices of $W$.
Assume that $\lvert S \rvert \geq b$, and for every pair of distinct vertices $x,y$ in $S$, $G$ contains four edge-disjoint paths from $x$ to $y$.

Since $W$ is simple, $\pi_E(e)$ does not contain any loop of $G$ for every $e \in E(W)$.
In addition, for every pair of distinct vertices $x,y$ in $S$, any path from $x$ to $y$ does not contain any loop.
So we may assume that $G$ is loopless by deleting all loops of $G$.

Let $G'$ be the graph obtained from $G$ by subdividing every edge once.
Let $H$ be the graph obtained from $L(G')$ by for each $v \in V(G) \subseteq V(G')$, adding a vertex $u_v$ adjacent to every vertex in $\cl(v)$ in $G'$.
Then it is clear that $H$ admits a $W$-subdivision $(\pi_V'',\pi_E'')$ such that $\pi_V''(x) = u_{\pi_V(x)}$ for every $x \in V(W)$. 
In particular, $\pi_V''(\pi_V^{-1}(S))=\{u_s: s \in S\}$.
Note that since $G$ is loopless, for every edge $e$ of $G$ with ends $x,y$, there exists an edge in $H$ with one end in $V(\cl(x))$ and one end in $V(\cl(y))$, and we also denote this edge in $H$ as $e$.
Since every wall does not contain a loop, the image of $\pi_E''$ of each edge of $W$ is path in $H$.

For every pair of distinct vertices $x,y$ of $S$, there exist four edge-disjoint paths $P_1,P_2,P_3,P_4$ in $G$ from $x$ to $y$, so there exist four paths $Q_1,Q_2,Q_3,Q_4$ in $H$ from $u_x$ to $u_y$ such that $E(Q_i)$ contains $E(P_i)$ for $1 \leq i \leq 4$.
If we choose those paths $Q_1,...,Q_4$ such that the sum of their length is minimum, then $Q_1,...,Q_4$ are pairwise edge-disjoint.
Therefore, by Lemma \ref{wall-subdiv grid-immersion}, there exists a $g \times g$ grid-immersion $(\pi_V''',\pi_E''')$ in $H$ such that the image of $\pi_V'''$ is contained in $\{u_s: s \in S\}$.

Note that if for every $v \in V(G)$, we identify the vertices in $\{u_v\} \cup \cl(v)$ into a vertex and delete all resulting loops, then we obtain $G$.
By the same procedure, we know there exists a $g \times g$ wall-immersion $(\pi_V^*,\pi_E^*)$ in $G$ such that the image of $\pi_V^*$ is $\{s \in S: u_s$ is in the image of $\pi_V'''\}$.
This proves the lemma.
\end{pf}

\bigskip

Recall that every large wall has a natural edge-tangle in it by Lemma \ref{wall has edge-tangle}, and every immersion induces an edge-tangle by Lemma \ref{immersion induces edge-tangle}.
The following lemma shows that every graph with no edge-cut of order three but with an edge-tangle induced by an immersion of a large wall has an edge-tangle controlling a large complete graph-thorns.

\begin{lemma} \label{4-edge-connected wall controlling clique}
For any positive integers $\theta$ and $t$, there exists a positive integer $w=w(\theta,t)$ with $w>\theta$ such that the following holds.
If $G$ is a graph with no edge-cut of order three, and $\E$ is an edge-tangle in $G$ of order $w$ induced by a $2w \times 4w$ wall-immersion and the natural tangle of order $w$ in the $2w \times 4w$ wall, then there exists an edge-tangle $\E' \subseteq \E$ of order at least $\theta$ in $G$ controlling a $K_t$-thorns.
\end{lemma}

\begin{pf}
Let $\theta$ and $t$ be positive integers.
Let $\theta'=\theta+t$.
Let $b$ be the number mentioned in Lemma \ref{wall-immersion grid-immersion} by taking $g=4\theta'$.
Define $w=b+2$.
Note that $w \geq (4\theta')^2+2>\theta'+2>\theta+2$.

Denote the $2w \times 4w$ wall by $W$ and denote the $4\theta' \times 4\theta'$ grid by $R$.
Let $S$ be the set of diagonal vertices of $W$ not contained in the first and the last column of $W$.
So $\lvert S \rvert \geq b$.
Let $G$ be a graph with no edge-cut of order three, and let $\E$ be an edge-tangle in $G$ of order $w$ induced by a $W$-immersion $(\pi_V,\pi_E)$ in $G$ and the natural edge-tangle in $W$ of order $w$.
By Lemma \ref{summary edge-tanlge induced by wall immersion}, for every edge-cut $[A,B]$ of $G$ of order less than $w$, $[A,B] \in \E$ if and only if $B$ contains the image of $\pi_V$ of all vertices of a column of $W$.

\noindent{\bf Claim 1:} For any two vertices $x,y$ in $\pi_V(S)$, there exist four edge-disjoint paths in $G$ between $x$ and $y$.

\noindent{\bf Proof of Claim 1:}
Let $x,y \in \pi_V(S)$.
So there exist $x',y' \in S$ such that $x=\pi_V(x')$ and $y=\pi_V(y')$.
Since $x'$ and $y'$ are diagonal vertices in $W$ not belong to the first and last column of $W$, there exist three edge-disjoint paths in $W$ from $x'$ to $y'$.
Hence there exist three edge-disjoint paths in $G$ from $x$ to $y$.
Suppose that there do not exist four edge-disjoint paths in $G$ from $x$ to $y$.
Then there exists an edge-cut $[A,B]$ of $G$ of order at most three such that $x \in A$ and $y \in B$.
Since there are three edge-disjoint paths in $G$ from $x$ to $y$, the order of $[A,B]$ is exactly three, contradiction that $G$ has no edge-cut of order three.
$\Box$

\smallskip

By Lemma \ref{wall-immersion grid-immersion} and Claim 1, $G$ admits an $R$-immersion $(\pi_V',\pi_E')$ such that $\pi_V'(V(R))$ is contained in $\pi_V(S)$.
Define $\E'$ to be the collection of all edge-cuts $[A,B]$ of $G$ of order less than $\theta'$ such that $B$ contains the image of $\pi_V'$ of all vertices of a row of $R$.

A wall that is a subgraph of $R$ is {\it canonical} if every its row is a subgraph of a row of $R$ and every its column is a subgraph of the union of two consecutive columns of $R$.
Note that for every canonical $2\theta' \times 4\theta'$ wall $W'$ and for every $[A,B] \in \E'$, $B$ contains the image of $\pi_V'$ of all vertices of a row of $W'$, so $B$ intersects the image of $\pi_V'$ of vertices in at least $\theta'$ columns of $W'$.

Since $R$ contains a canonical $2\theta' \times 4\theta'$ wall $W^*$ as a subgraph, for every $[A,B] \in \E'$, $B$ intersects the image of $\pi_V'$ of vertices in at least $\theta'$ columns of $W^*$.
By Lemma \ref{summary edge-tanlge induced by wall immersion}, $\E'$ is the edge-tangle in $G$ of order $\theta' \geq \theta$ induced by an $W^*$-immersion and the nature edge-tangle in $W^*$ of order $\theta'$. 

For every $i$ with $1 \leq i \leq t \leq \theta'$, define $\alpha(v_i)$ to be the union of the image of $\pi_E'$ of the edges in the $i$-th column and the edges in the $i$-th row of $R$, where we write $V(K_t)=\{v_j: 1 \leq j \leq t\}$.
So $\alpha$ is a $K_t$-thorns.

We claim that $\E'$ controls $\alpha$.
Suppose to the contrary that there exist $[A,B] \in \E'$ with order less than $t$ and $v \in V(K_t)$ such that $V(\alpha(v)) \cap B = \emptyset$.
Since $[A,B] \in \E'$, $B$ contains the image of $\pi_V'$ of all vertices of a row of $R$.
Since $\alpha(v)$ intersects the image of $\pi_V'$ of each row, $B \cap V(\alpha(v)) \neq \emptyset$, a contradiction.
Hence $\E'$ controls a $K_t$-thorns $\alpha$. 

It suffices to prove that $\E' \subseteq \E$ to complete the proof.
Let $[A,B] \in \E'$.
So the order of $[A,B]$ is less than $\theta'$.
Since $\pi_V'(V(R)) \subseteq \pi_V(S)$ and $B$ contains the image of $\pi_V'$ of all vertices of a row of $R$, we know $B$ contains at least $\theta'$ vertices in $\pi_V(S)$.
Since different vertices in $S$ belong to different columns of $W$, $B$ intersects the image of $\pi_V$ of vertices in at least $\theta'$ columns of $W$.
Since $\theta'<w$, $[A,B] \in \E$ by Lemma \ref{summary edge-tanlge induced by wall immersion}.
This proves that $\E' \subseteq \E$.
\end{pf}

\bigskip

The following theorem is the main result of this section.

\begin{theorem} \label{edge-tangle in 4-edge-connected control edge-minor}
For any positive integers $k$ and $\theta$ with $\theta > k$, there exists a positive integer $w=w(k,\theta)$ such that if $G$ is a graph with no edge-cut of order three, and $\E$ is an edge-tangle in $G$ of order at least $w$, then $\E_\theta$ controls a $K_k$-thorns, where $\E_\theta$ is the edge-tangle in $G$ of order $\theta$ such that $\E_\theta \subseteq \E$.
\end{theorem}

\begin{pf}
Let $k$ and $\theta$ be positive integers with $\theta>k$.
Note that $\theta \geq 2$.
Let $w_1=w_{\ref{4-edge-connected wall controlling clique}}(\theta,k)$, where $w_{\ref{4-edge-connected wall controlling clique}}$ is the integer $w$ mentioned in Lemma \ref{4-edge-connected wall controlling clique}.
Note that $w_1>\theta \geq 2$ by Lemma \ref{4-edge-connected wall controlling clique}.
Define $w=w_{\ref{truncation large deg or wall}}(w_1,k)$, where $w_{\ref{truncation large deg or wall}}$ is the integer $w$ mentioned in Lemma \ref{truncation large deg or wall}.

For every integer $t$ and for every edge-tangle $\E$ in a graph of order at least $t$, let $\E_t$ be the edge-tangle in the same graph of order $t$ such that $\E_t \subseteq \E$.

Let $G$ be a graph with no edge-cut of order three, and let $\E$ be an edge-tangle in $G$ of order at least $w$.
By Lemma \ref{truncation large deg or wall}, either 
	\begin{itemize}
		\item[(i)] there exists $v \in V(G)$ incident with at least $k$ edges in $G$ such that $v \in B$ for every $[A,B] \in \E_{w_1}$, or
		\item[(ii)] $\E_{w_1}$ is induced by a $2w_1 \times 4w_1$ wall-immersion and the natural edge-tangle of order $w_1$ in the $2w_1 \times 4w_1$ wall.
	\end{itemize}

We first assume that (i) holds.
Define $\alpha$ to be a $K_k$-thorns such that $\alpha(h)$ is an edge of $G$ incident with $v$ for each $h \in V(K_k)$.
We shall prove that $\E_{\theta}$ controls $\alpha$.
Let $[A,B] \in \E_{\theta}$ with order less than $k$.
Since $\theta<w_1$, $[A,B] \in \E_{w_1}$.
Hence, $v \in B \cap V(\alpha(h))$ for every $h \in V(K_k)$.
Therefore, $\E_{\theta}$ controls a $K_k$-thorns.

So we may assume that (ii) holds.
That is, $\E_{w_1}$ is induced by a $2w_1 \times 4w_1$ wall-immersion and the natural edge-tangle of order $w_1$ in the $2w_1 \times 4w_1$ wall.
By Lemma \ref{4-edge-connected wall controlling clique}, there exists an edge-tangle $\E' \subseteq \E_{w_1}$ of order at least $\theta$ in $G$ controlling a $K_k$-thorns.
Therefore, $\E_\theta = \E'_\theta$ controls a $K_k$-thorns.
\end{pf}

\section{Erd\H{o}s-P\'{o}sa property}
\label{sec:EP main}

We say that a graph $G$ is {\it nearly 3-cut free} if $\lvert V(G) \rvert \geq 2$, $G$ is connected and for every edge-cut of $G$ of order three, the edges between $A$ and $B$ are parallel with the same ends.

\begin{lemma} \label{decomp into 4-edge-conn}
If $G$ is either a nearly 3-cut free graph or a graph with $\lvert V(G) \rvert=1$, then there exist a tree $T$ and a partition $\{X_t: t \in V(T)\}$ of $V(G)$ such that the following hold.
	\begin{enumerate}
		\item For every $t \in V(T)$, either $\lvert X_t \rvert=1$ or $G[X_t]$ does not have an edge-cut of order three. 
		\item If there is an edge of $G$ with one end in $X_{t_1}$ and one end in $X_{t_2}$ for some distinct $t_1,t_2 \in V(T)$, then $t_1$ is adjacent to $t_2$ in $T$.
		\item For every edge $t_1t_2$ of $T$, there are exactly three edges with one end in $X_{t_1}$ and one end in $X_{t_2}$, and those edges are parallel with the same ends.
	\end{enumerate}
\end{lemma}

\begin{pf}
We prove this lemma by induction on $\lvert V(G) \rvert$.
If either $G$ does not have an edge-cut of order three or $G$ has only one vertex, then we are done by taking the tree on one vertex and the partition of $V(G)$ with one part.
This proves the base case and we may assume that $\lvert V(G) \rvert \geq 2$ and the lemma holds for every nearly 3-cut free graph on less than $\lvert V(G) \rvert$ vertices.
And we may assume that there exists an edge-cut $[A,B]$ of $G$ of order three.
Since $G$ is nearly 3-cut free with $\lvert V(G) \rvert \geq 2$, the edges between $A$ and $B$ are three parallel edges with the same ends $u,v$, say $u \in A$ and $v \in B$.

Suppose that $\lvert A \rvert \geq 2$ and $G[A]$ is not nearly 3-cut free.
Then there exists an edge-cut $[A',B']$ of $G[A]$ of order zero or three such that either there is no edge between $A'$ and $B'$, or the edges between $A',B'$ are not parallel with the same ends.
By symmetry, we may assume that $u \in B'$.
So $[A',B' \cup B]$ is an edge-cut of $G$ such that the edges between $A',B' \cup B$ are the edges between $A',B'$.
Since $G$ is nearly 3-cut free, there is at least one edge between $A',B' \cup B$, and the edges between $A', B' \cup B$ are parallel with the same ends.
Therefore, there is at least one edge between $A',B'$, and the edges between $A',B'$ are parallel with the same ends, a contradiction.

Hence either $\lvert A \rvert=1$ or $G[A]$ is nearly 3-cut free.
Similarly, either $\lvert B \rvert=1$ or $G[B]$ is nearly 3-cut free.
By the induction hypothesis, there exist trees $T_A,T_B$, a partition $\{Y_t: t \in V(T_A)\}$ of $A$ and a partition $\{Z_t: t \in V(T_B)\}$ of $B$ satisfying the three properties mentioned in the lemma.
Let $t_u \in V(T_A)$ and $t_v \in V(T_B)$ be the vertices such that $u \in X_{t_u}$ and $v \in X_{t_v}$.
Define $T$ to be the tree obtained from the union of $T_A$ and $T_B$ by adding the edge $t_ut_v$.
For every $t \in V(T)$, define $X_t=Y_t$ if $t \in V(T_A)$, and $X_t=Z_t$ if $t \in V(T_B)$.
Then $T$ and the partition $\{X_t: t \in V(T)\}$ of $V(G)$ satisfy the three properties mentioned in the lemma.
\end{pf}

\bigskip

Now we are ready to address the Erd\H{o}s-P\'{o}sa property.
The purpose of Lemma \ref{EP not isolated} is to deal with the main difficulty of the proof of Theorem \ref{main intro}.
Lemma \ref{EP not isolated} implies Theorem \ref{main intro} for the case when $H$ has no isolated vertices and $G$ is nearly 3-cut free.

We give the intuition of the statement of Lemma \ref{EP not isolated} and sketch its proof.
We shall prove that given a nearly 3-cut free graph $G$, if $G$ does not contain $k$ edge-disjoint $H$-immersions, then we can hit all $H$-immersions in $G$ by a set of edges with bounded size.
We assume that $H$ is connected in the proof sketch, as the case that $H$ is disconnected follows from a relatively easier argument by (more or less) induction on the number of components of $H$.
We shall prove it by induction on $k$, and assume that $G$ does not contain $k$ edge-disjoint $H$-immersions.
Note that as long as there exist an edge-cut $[A,B]$ of $G$, a hitting set of $H$-immersions of $G[A]$ and a hitting set of $H$-immersions of $G[B]$, we can obtain a hitting set of $H$-immersions of $G$ by taking the union of those two hitting sets together with all edges between $A,B$, since $H$ is connected.
So there is a win if there exists an edge-cut $[A,B]$ of $G$ of small order such that each $G[A]$ and $G[B]$ has a hitting set of small size.
If both $G[A]$ and $G[B]$ contain $H$-immersions, then each of $G[A]$ and $G[B]$ does not contain $k-1$ edge-disjoint $H$-immersions, so we expect to obtain hitting sets of $H$-immersions in $G[A]$ and $G[B]$ by induction on $k$.
However, the induction does not apply, as $G[A]$ and $G[B]$ might not be nearly 3-cut free.
So instead of considering $G[A]=G-B$ and $G[B]=G-A$, we consider the graph $G_A$ obtained from $G$ by contracting $B$ and the graph $G_B$ obtained from $G$ by contracting $A$.
Note that $G_A$ and $G_B$ are nearly 3-cut free.
But contracting a subset of $V(G)$ can create more $H$-immersions.
So we should treat those new vertices obtained by contractions as special vertices.
This is the purpose of the set $S$ and function $\gamma$ stated in Lemma \ref{EP not isolated}.
It can be helpful (though not completely true) to think that each vertex $v$ in $S$ corresponds to a subset of vertices that induces a subgraph that contains $\gamma(v)$ edge-disjoint $H$-immersions.
This setting allows us to apply induction on $k-\sum_{v \in S}\gamma(v)$ for $G_A$ and $G_B$.
In addition, each of those special vertices is obtained by contracting one side of an edge-cut, so its degree equals the order of the edge-cut.
As we will only contract the sides of edge-cuts of bounded order, the degree of those special vertices in $S$ is bounded.
Since $H$ has no isolated vertices, each $H$-immersion in $G$ intersecting $S$ must intersect an edge incident with a vertex in $S$.
Hence if we can hit all $H$-immersions in $G-S$ by a set of edges of bounded size, then we can hit all $H$-immersions in $G$ by a set of edges of bounded size by further including all edges incident with $S$, as long as $\lvert S \rvert$ is bounded.
Indeed, $\lvert S \rvert$ is bounded by $\sum_{v \in S}\gamma(v)$.

Further intuition and proof sketch of Lemma \ref{EP not isolated} will be stated after we prove the following easy lemma which is the base case of Lemma \ref{EP not isolated}.

\begin{lemma} \label{one edge eEP}
For every connected graph $H$ that has exactly one edge and every function $g: {\mathbb N} \times ({\mathbb N} \cup \{0\}) \rightarrow {\mathbb N} \cup \{0\}$, there exists a function $f: {\mathbb N} \times ({\mathbb N} \cup \{0\}) \rightarrow {\mathbb N} \cup \{0\}$ such that for every graph $G$, every positive integer $k$, every $S \subseteq V(G)$ and every function $\gamma: S \rightarrow {\mathbb N}$, if $S$ does not contain any vertex of degree at least $g(k,\sum_{v \in S}\gamma(v))$, then either $G-S$ contains $k-\sum_{v \in S}\gamma(v)$ edge-disjoint $H$-immersions, or there exists $Z \subseteq E(G)$ with $\lvert Z \rvert \leq f(k, \sum_{v \in S}\gamma(v))$ such that $G-Z$ does not contain an $H$-immersion.
\end{lemma}

\begin{pf}
Let $H$ be a connected graph that has exactly one edge, and let $g: {\mathbb N} \times ({\mathbb N} \cup \{0\}) \rightarrow {\mathbb N} \cup \{0\}$ be a function.
Note that $H$ is either $K_2$ or the one-vertex graph with one loop.
It was shown in \cite[Chapter 9, Exercise 6]{d} that there exists a function $h: {\mathbb N} \rightarrow {\mathbb N}$ such that for every simple graph $G$, either $G$ contains $k$ edge-disjoint cycles, or there exists $Z \subseteq E(G)$ with $\lvert Z \rvert \leq h(k)$ such that $G-Z$ has no cycle.
Define $f$ to be the function such that $f(x,y)=h(x+yg(x,y))+2(x+yg(x,y))^2$ for any $x \in {\mathbb N}$ and $y \in {\mathbb N} \cup \{0\}$.

Let $G$ be a graph, $k$ a positive integer, $S$ a subset of $V(G)$ and $\gamma: S \rightarrow {\mathbb N}$ a function such that $S$ does not contain any vertex of degree at least $g(k,\sum_{v \in S}\gamma(v))$.
Let $d=\sum_{v \in S}\gamma(v)$.

Since every vertex of $S$ has degree less than $g(k,d)$, if $G$ contains $k+\lvert S \rvert g(k,d)$ edge-disjoint $H$-immersions, then $G-S$ contains $k \geq k-d$ edge-disjoint $H$-immersions.
Note that $\lvert S \rvert \leq d$ since $\gamma(v) \geq 1$ for every $v \in S$.
So to prove this lemma, it suffices to prove that either $G$ contains $k+dg(k,d)$ edge-disjoint $H$-immersions, or there exists $Z \subseteq E(G)$ with $\lvert Z \rvert \leq f(k,d)$ such that $G-Z$ does not contain an $H$-immersion.

We first assume that $H=K_2$.
If $G$ contains at least $k+dg(k,d)$ non-loop edges, then $G$ contains $k+dg(k,d)$ edge-disjoint $H$-immersions; if $G$ has less than $k + dg(k,d)$ non-loop edges, there exists $Z \subseteq E(G)$ with $\lvert Z \rvert \leq k+dg(k,d) \leq f(k,d)$ such that $G-Z$ has no non-loop edge and has no $H$-immersion.
So this lemma holds if $H=K_2$.

Now we assume that $H$ is the one-vertex graph with one loop.
We assume that $G$ does not contain $k+dg(k,d)$ edge-disjoint $H$-immersions and show that there exists a set $Z \subseteq E(G)$ with $\lvert Z \rvert \leq f(k,d)$ such that $G-Z$ has no $H$-immersion.
So $G$ does not contain $k+dg(k,d)$ loops, and there do not exist two distinct vertices such that there are $2(k+dg(k,d))$ parallel edges between them.
In addition, $G$ does not contain $k+dg(k,d)$ distinct pairs of distinct vertices of $G$ such that there are at least two edges between each pair.
Hence $G$ has at most $k+dg(k,d)-1$ loops, no pair of distinct vertices of $G$ has at least $2(k+dg(k,d))$ edges between them, and there are at most $k+dg(k,d)-1$ pairs of distinct vertices having parallel edges between them.
So there exists $Z_1 \subseteq E(G)$ with $\lvert Z_1 \rvert \leq k+dg(k,d)-1 + (k+dg(k,d)-1)(2(k+dg(k,d))-1) \leq 2(k+dg(k,d))^2$ such that $G-Z_1$ is simple.
Since $G$ does not contain $k+dg(k,d)$ edge-disjoint $H$-immersions, $G-Z_1$ is a simple graph that does not contain $k+dg(k,d)$ edge-disjoint cycles.
By \cite[Chapter 9, Exercise 6]{d}, there exists $Z_2 \subseteq E(G-Z_1)$ with $\lvert Z_2 \rvert \leq h(k+dg(k,d))$ such that $G-(Z_1 \cup Z_2)$ is a simple graph with no cycle and hence has no $H$-immersion.
This proves the lemma since $\lvert Z_1 \cup Z_2 \rvert \leq h(k+dg(k,d))+2(k+dg(k,d))^2 \leq f(k,d)$.
\end{pf}

\bigskip

Now we continue the intuition and proof sketch of Lemma \ref{EP not isolated}.
Recall that we aim to prove that given a nearly 3-cut free graph $G$ and a connected graph $H$, if $G-S$ does not contain $k-\sum_{v \in S}\gamma(v)$ edge-disjoint $H$-immersions, then we can hit all $H$-immersions in $G$ by a set of edges with bounded size, where $S$ is a special set of vertices whose size is bounded by $\sum_{v \in S}\gamma(v)$.
Also recall that our setting for the set $S$ of special vertices allows us to apply induction on $G_A$ and $G_B$ whenever we have an edge-cut $[A,B]$ of small order such that both $G[A]$ and $G[B]$ contain $H$-immersions, and the degree of the vertices in $S$ can be bounded if we only work on edge-cuts of bounded order.
So now we may assume that there exists no edge-cut $[A,B]$ of $G$ of small order such that each $G[A]$ and $G[B]$ contains an $H$-immersion.
This will allow us to define an edge-tangle in $G$ of large order (see Claims 3-5), by simply seeing which side of each edge-cut contains an $H$-immersion.
Note that the order of the edge-tangle is related to the degree of the vertices in $S$ and the order of the edge-cuts that we can work with.
For a technical reason, we need this number to be depend on $\lvert S \rvert$ (or more precisely, $\sum_{v \in S}\gamma(v)$).
And that is the reason why we consider the function $g$ in Lemma \ref{EP not isolated} to indicate the degree condition of the vertices of $S$.
For another technical reason, we want this function $g$ growing sufficiently quickly, and that is the motivation of the notion of ``$H$-legal'' functions defined below.
If $G$ has no edge-cut of order three, then we know this edge-tangle controls a $K_w$-thorns for some large $w$ by Theorem \ref{edge-tangle in 4-edge-connected control edge-minor}, and hence we can obtain a hitting set by Lemma \ref{isolating immersion}.
So we may assume that $G$ is nearly 3-cut free but has an edge-cut of order three.
Hence we can decompose $G$ into pieces with no edge-cut of order three in a tree-like fashion by Lemma \ref{decomp into 4-edge-conn}.
Claims 6-8 tell us that we can use the tree to reduce the problem to a piece of $G$ with no edge-cut of order three and hence complete the proof.

Recall that an isolated vertex in a graph is a vertex of degree zero.
For a graph $H$ with no isolated vertices, we say that a function $g$ is {\it $H$-legal} if $g$ is a function from ${\mathbb N} \times ({\mathbb N} \cup \{0\})$ to ${\mathbb N} \cup \{0\}$ satisfying that 
	\begin{itemize}
		\item $g(x,y) \geq g(x,y')+2y$ and $g(x,y) \geq g(x',y)$ for every $x,x' \in {\mathbb N}$, $y,y' \in {\mathbb N} \cup \{0\}$ with $x \geq x'$ and $y>y'$, and 
		\item for any positive integers $m$ and $n$, 
			\begin{align*}
				& g(m,n) \\
				\geq & \max_{H'}\big\{w_{\ref{edge-tangle in 4-edge-connected control edge-minor}} \big( w_{\ref{isolating immersion}}(H',m+(n-1) \cdot g(m,n-1)), \theta_{\ref{isolating immersion}}(H',m+(n-1) \cdot g(m,n-1)) \big) \\
				 & +3m\lvert V(H') \rvert d_{H'}\big\},
			\end{align*}
			where the maximum is over all graphs $H'$ with no isolated vertices and with $\lvert E(H') \rvert \leq \lvert E(H) \rvert$, and $w_{\ref{edge-tangle in 4-edge-connected control edge-minor}}$ is the integer $w$ mentioned in Theorem \ref{edge-tangle in 4-edge-connected control edge-minor}, and $\theta_{\ref{isolating immersion}}, w_{\ref{isolating immersion}}$ are the integers $\theta,w$ mentioned in Lemma \ref{isolating immersion}, respectively, and $d_{H'}$ is the maximum degree of $H'$. 
	\end{itemize}
Note that if $g$ is $H$-legal for some graph $H$ with no isolated vertices, then $g$ is $H''$-legal for any graph $H''$ with no isolated vertices with $\lvert E(H'') \rvert \leq \lvert E(H) \rvert$.
And it is easy to see that $H$-legal functions exist for any graph $H$ with no isolated vertices.

\begin{lemma} \label{EP not isolated}
For every graph $H$ with no isolated vertices, there exists an $H$-legal function $g^*: {\mathbb N} \times ({\mathbb N} \cup \{0\}) \rightarrow {\mathbb N} \cup \{0\}$ such that for every $H$-legal function $g: {\mathbb N} \times ({\mathbb N} \cup \{0\}) \rightarrow {\mathbb N} \cup \{0\}$ with $g \geq g^*$, there exists a function $f: {\mathbb N} \times ({\mathbb N} \cup \{0\}) \rightarrow {\mathbb N} \cup \{0\}$ such that for every nearly 3-cut free graph $G$, every positive integer $k$, every $S \subseteq V(G)$ and every function $\gamma: S \rightarrow {\mathbb N}$, if $S$ does not contain any vertex of degree at least $g(k,\sum_{v \in S}\gamma(v))$, then either $G-S$ contains $k-\sum_{v \in S}\gamma(v)$ edge-disjoint $H$-immersions, or there exists $Z \subseteq E(G)$ with $\lvert Z \rvert \leq f(k, \sum_{v \in S}\gamma(v))$ such that $G-Z$ does not contain an $H$-immersion.
\end{lemma}

\begin{pf}
Let $H$ be a graph with no isolated vertices.
Denote $\lvert V(H) \rvert$ by $h$ and the maximum degree of $H$ by $d$.
Since $H$ has no isolated vertices, $d \geq 1$.

We shall prove this lemma by induction on $\lvert E(H) \rvert$.
If $H$ contains only one edge, then $H$ is connected since $H$ has no isolated vertices, so the lemma holds by Lemma \ref{one edge eEP} by choosing $g^*$ to be any $H$-legal function.
This proves the base case of the induction.
We assume that this lemma is true for every graph $H'$ without isolated vertices with $\lvert E(H') \rvert < \lvert E(H) \rvert$ and denote the corresponding function $g^*$ and the corresponding function $f$ (when some $H'$-legal function $g$ with $g \geq g^*$ is given) by $g^*_{H'}$ and $f_{H',g}$, respectively.

We define the following.
	\begin{itemize}
		\item Define $g^*: {\mathbb N} \times ({\mathbb N} \cup \{0\}) \rightarrow {\mathbb N} \cup \{0\}$ such that the following hold.
			\begin{itemize}
				\item For every positive integer $m$, define $g^*(m,0)=w_{\ref{edge-tangle in 4-edge-connected control edge-minor}}(w_{\ref{isolating immersion}}(H,m),\theta_{\ref{isolating immersion}}(H,m))+3mhd + \sum_{H'}g^*_{H'}(m,0)$, where $w_{\ref{edge-tangle in 4-edge-connected control edge-minor}}$ is the integer $w$ mentioned in Theorem \ref{edge-tangle in 4-edge-connected control edge-minor}, and $\theta_{\ref{isolating immersion}}, w_{\ref{isolating immersion}}$ are the integers $\theta,w$ mentioned in Lemma \ref{isolating immersion}, respectively, and the last sum is over all graphs $H'$ with no isolated vertices and with less edges than $H$.
				\item For every positive integers $m,n$, define $g^*(m,n)= w_{\ref{edge-tangle in 4-edge-connected control edge-minor}}(w_{\ref{isolating immersion}}(H,m+(n-1) \cdot g^*(m,n-1)),\theta_{\ref{isolating immersion}}(H,m+(n-1) \cdot g^*(m,n-1)))+3mhd+g^*(m,n-1)+2n+\theta_{\ref{isolating immersion}}(H,m+(n-1) \cdot g^*(m,n-1)) + \sum_{H'}g^*_{H'}(m,n)$, where the last sum is over all graphs $H'$ with no isolated vertices and with less edges than $H$.
			\end{itemize} 
		Note that $\theta_{\ref{isolating immersion}}(H,t) > w_{\ref{isolating immersion}}(H,t)$ for any positive integer $t$ by Lemma \ref{isolating immersion}.
			Clearly, $g^*$ is $H$-legal.
	\item For every $H$-legal function $g$ with $g \geq g^*$, we define the following.
		\begin{itemize}
			\item Let $f'_{H,g}: {\mathbb N} \rightarrow {\mathbb N} \cup \{0\}$ be the function such that $f'_{H,g}(k) = \sum_{H'} \sum_{i=0}^{k} f_{H',g}(k,i)$ for every positive integer $k$, where the first sum is taken over all graphs $H'$ with no isolated vertices having less edges than $H$.
		Note that there are only finitely many such graphs $H'$, as $H'$ has no isolated vertices.
			\item Define $f: {\mathbb N} \times ({\mathbb N} \cup \{0\}) \rightarrow {\mathbb N} \cup \{0\}$ to be the function satisfying the following.
				\begin{itemize}
					\item $f(m,n)=0$ for every integers $m,n$ with $0<m \leq n$.
					\item $f(m,n)=2f(m,n+1)+(m+n+2)g(m,n+1)+(2^{h+1}-4)f'_{H,g}(mg(m,n+2g(m,n+1))) + mhd$ for every integers $m,n$ with $m >n \geq 0$.
				\end{itemize}
				Note that $f$ depends on $g$, but we do not add subscript $g$ to describe $f$ for simplicity of notations.
		\end{itemize}
	\end{itemize}

We shall prove that the functions $g^*$ and $f$ defined above satisfy the conclusion of this lemma for the graph $H$.
That is, we shall prove that for every $H$-legal function $g$ with $g \geq g^*$, the function $f$ satisfies the property that for every nearly 3-cut free graph $G$, every positive integer $k$, every set $S \subseteq V(G)$ and every function $\gamma: S \rightarrow {\mathbb N}$ such that $S$ does not contain any vertex of degree at least $g(k, \sum_{v \in S}\gamma(v))$, either $G-S$ contains $k-\sum_{v \in S}\gamma(v)$ edge-disjoint $H$-immersions, or there exists $Z \subseteq E(G)$ with $\lvert Z \rvert \leq f(k,\sum_{v \in S}\gamma(v))$ such that $G-Z$ does not contain an $H$-immersion.

We do induction on $k-\sum_{v \in S}\gamma(v)$.
Suppose to the contrary that there exists a tuple $(g,G,k,S,\gamma)$ such that the following hold.
	\begin{itemize}
		\item[(i)] $g$ is an $H$-legal function with $g \geq g^*$, $G$ is a nearly 3-cut free graph, $k$ is a positive integer, $S$ is a subset of $V(G)$, and $\gamma: S \rightarrow {\mathbb N}$ is a function such that $S$ does not contain any vertex of degree at least $g(k,\sum_{v \in S}\gamma(v))$.
		\item[(ii)] $G-S$ does not contain $k-\sum_{v \in S}\gamma(v)$ edge-disjoint $H$-immersions, but there does not exist $Z \subseteq E(G)$ with $\lvert Z \rvert \leq f(k,\sum_{v \in S}\gamma(v))$ such that $G-Z$ has no $H$-immersion.
		\item[(iii)] Subject to (i) and (ii), $k-\sum_{v \in S}\gamma(v)$ is minimum.
	\end{itemize}
Note that (ii) implies that $k-\sum_{v \in S}\gamma(v) \geq 1$, so the minimum mentioned in (iii) exists.

In the rest of the proof, we denote $\sum_{v \in S}\gamma(v)$ by $\bar{r}$.

For every edge-cut $[A,B]$ of $G$ with $A \neq \emptyset \neq B$, define $G_A$ (and $G_B$, respectively) to be the graphs obtained from $G$ by identifying $B$ (and $A$, respectively) into one new vertex $v_B$ (and $v_A$, respectively), and deleting all resulting loops.
Define $S_A = (S \cap A) \cup \{v_B\}$ and $S_B = (S \cap B) \cup \{v_A\}$.
Note that the degree of $v_B$ in $G_A$ and the degree of $v_A$ in $G_B$ are the order of $[A,B]$.

\noindent{\bf Claim 1:} For every edge-cut $[A,B]$ of $G$ with $A \neq \emptyset \neq B$, $G_A$ and $G_B$ are nearly 3-cut free.

\noindent{\bf Proof of Claim 1:}
Since $A \neq \emptyset \neq B$, $G_A$ and $G_B$ contain at least two vertices.
Since $G$ is nearly 3-cut free, $G$ is connected, so $G_A$ and $G_B$ are connected.
If $G_A$ is not nearly 3-cut free, then there exists an edge-cut $[X,Y]$ of $G_A$ with $v_B \in Y$ of order three such that the edges between $X$ and $Y$ are not parallel edges with the same ends.
But then $[X,(Y-\{v_B\}) \cup B]$ is an edge-cut of $G$ of order three such that the edges in between are not parallel edges with the same ends, a contradiction.
So $G_A$ is nearly 3-cut free.
Similarly, $G_B$ is nearly 3-cut free.
$\Box$

\smallskip

\noindent{\bf Claim 2:} Let $\theta$ be a positive integer.
If there exist $W \subseteq V(G)$, an edge-cut $[A,B]$ of $G$ of order less than $\theta$ and a set $Z_0 \subseteq E(G)$ containing all edges between $A$ and $B$ such that $G[A]-(W \cup Z_0)$ and $G[B]-(W \cup Z_0)$ do not contain $H$-immersions, then there exists $Z$ with $Z_0 \subseteq Z \subseteq E(G)$ and $\lvert Z \rvert \leq \lvert Z_0 \rvert + (2^{h}-2)f'_{H,g}(kg(k,\bar{r}+\theta))$ such that $G-(W \cup Z)$ has no $H$-immersion.

\noindent{\bf Proof of Claim 2:}
If $A=\emptyset$, then $G=G[B]$, so we are done by taking $Z=Z_0$.
Similarly, we are done if $B=\emptyset$.
So we may assume that $A \neq \emptyset \neq B$ and hence $G_A$ and $G_B$ contain at least two vertices and are nearly 3-cut free by Claim 1.

If $H$ is connected, then every $H$-immersion in $G-Z_0$ must be in $G[A]$ or $G[B]$, as $Z_0$ contains all edges between $A$ and $B$.
So we are done by taking $Z=Z_0$.

Now we assume that $H$ is not connected.
Let $H_1,H_2,...,H_p$ be the components of $H$, where $p \geq 2$.
For every set $I$ with $\emptyset \subset I \subset [p]$, define $Q_I$ to be the disjoint union of $H_i$ over all $i \in I$.

Since $H$ has no isolated vertices, every $H$-immersion in $G_A$ (or $G_B$, respectively) intersecting $S_A$ (or $S_B$, respectively) must intersect an edge incident with a vertex in $S_A$ (or $S_B$, respectively).
Since every vertex in $S_A-\{v_B\}$ has degree in $G_A$ at most $g(k,\bar{r})-1$ in $G_A$ and $v_B$ has degree in $G_A$ at most $\theta-1$, for every $I$ with $\emptyset \subset I \subset [p]$ and any nonnegative integer $k'$, if $G_A$ contains at least $k'+(\lvert S_A \rvert-1) (g(k,\bar{r})-1) + \theta-1$ edge-disjoint $Q_I$-immersions, then $G[A]-S_A$ contains $k'$ edge-disjoint $Q_I$-immersions.
Similarly, for every $I$ with $\emptyset \subset I \subset [p]$ and any nonnegative integer $k'$, if $G_B$ contains at least $k'+(\lvert S_B \rvert-1) (g(k,\bar{r})-1) + \theta-1$ edge-disjoint $Q_I$-immersions, then $G[B]-S_B$ contains $k'$ edge-disjoint $Q_I$-immersions.

Note that for every nonnegative integer $k'$ and set $I$ with $\emptyset \subset I \subset [p]$, if $G[A]-S_A$ contains $k'$ edge-disjoint $Q_I$-immersions and $G[B]-S_B$ contains $k'$ edge-disjoint $Q_{[p]-I}$-immersions, then $G-S$ contains $k'$ edge-disjoint $H$-immersions.
Hence, since $G-S$ has no $k-\bar{r}$ edge-disjoint $H$-immersions, for every $I$ with $\emptyset \subset I \subset [p]$, either $G_A$ does not contain $(k-\bar{r})+(\lvert S_A \rvert-1) (g(k,\bar{r})-1) + \theta-1$ edge-disjoint $Q_I$-immersions, or $G_B$ does not contain $(k-\bar{r})+(\lvert S_B \rvert-1) (g(k,\bar{r})-1) + \theta-1$ edge-disjoint $Q_{[p]-I}$-immersions.

As $k-\bar{r}  \geq 1$, $\max\{\lvert S_A \rvert, \lvert S_B \rvert\} \leq \lvert S \rvert +1 \leq \bar{r} +1 \leq k$.
Hence, for every $I$ with $\emptyset \subset I \subset [p]$, either $G_A$ does not contain $(k-1)g(k, \bar{r})-\bar{r}  + \theta$ edge-disjoint $Q_I$-immersions, or $G_B$ does not contain $(k-1)g(k,\bar{r})-\bar{r} +\theta$ edge-disjoint $Q_{[p]-I}$-immersions.

Define $\gamma_A: S_A \rightarrow {\mathbb N}$ to be the function such that $\gamma_A(v_B)=\theta+\sum_{v \in S \cap B}\gamma(v)$ and $\gamma_A(x)=\gamma(x)$ for every $x \in S \cap A$.
Define $\gamma_B:S_B \rightarrow {\mathbb N}$ to be the function such that $\gamma_B(v_A)=\theta+\sum_{v\in S \cap A}\gamma(v)$ and $\gamma_B(x)=\gamma(x)$ for every $x \in S \cap B$.

Note that $\sum_{v \in S_A}\gamma_A(v) = \theta+\bar{r}$ and $\sum_{v \in S_B}\gamma_B(v) = \theta+\bar{r}$.
So $g(k,\sum_{v \in S_A}\gamma_A(v)) = g(k,\theta+\bar{r}) \geq g(k,\bar{r}) + 2\theta$ since $g$ is $H$-legal.
Similarly, $g(k,\sum_{v \in S_B}\gamma_B(v)) \geq g(k,\bar{r}) + 2\theta$.

Let $k_A=kg(k,\sum_{v \in S_A}\gamma_A(v))$ and let $k_B=kg(k,\sum_{v\in S_B}\gamma_B(v))$.
Hence, if $I$ is a set with $\emptyset \subset I \subset [p]$ such that $G_A-S_A$ does not contain $(k-1)g(k, \bar{r})-\bar{r}  + \theta$ edge-disjoint $Q_I$-immersions, then $G_A-S_A$ does not contain $((k-1)g(k,\bar{r}) - \bar{r} + \theta + \sum_{v \in S_A}\gamma_A(v))-\sum_{v \in S_A}\gamma_A(v) = ((k-1)g(k,\bar{r})+2\theta)-\sum_{v \in S_A}\gamma_A(v) \leq kg(k,\sum_{v \in S_A}\gamma_A(v))-\sum_{v \in S_A}\gamma_A(v) = k_A-\sum_{v \in S_A}\gamma_A(v)$ edge-disjoint $Q_I$-immersions.
Similarly, if $I$ is a set with $\emptyset \subset I \subset [p]$ such that $G_B-S_B$ does not contain $(k-1)g(k,\bar{r})-\bar{r} +\theta$ edge-disjoint $Q_{[p]-I}$-immersions, then $G_B-S_B$ does not contain $k_B-\sum_{v \in S_B}\gamma_B(v)$ edge-disjoint $Q_{[p]-I}$-immersions.

Therefore, for every $I$ with $\emptyset \subset I \subset [p]$, either $G_A-S_A$ does not contain $k_A-\sum_{v \in S_A}\gamma_A(v)$ edge-disjoint $Q_I$-immersions, or $G_B$ does not contain $k_B-\sum_{v \in S_B}\gamma_B(v)$ edge-disjoint $Q_{[p]-I}$-immersions.

Note that every vertex in $S_A$ has degree in $G_A$ less than $g(k,\bar{r})+\theta \leq g(k,\sum_{v \in S_A}\gamma_A(v)) \leq g(k_A,\sum_{v \in S_A}\gamma_A(v))$, since $g$ is $H$-legal.
Similarly, every vertex in $S_B$ has degree in $G_B$ less than $g(k,\bar{r})+\theta \leq g(k_B,\sum_{v \in S_B}\gamma_B(v))$.

Recall that $G_A$ and $G_B$ are nearly 3-cut free graphs.
For every $I$ with $\emptyset \subset I \subset [p]$, $Q_I$ and $Q_{[p]-I}$ are graphs with no isolated vertices and with less edges than $H$, $g$ is $Q_I$-legal and $Q_{[p]-I}$-legal, and $g \geq g^* \geq g^*_{Q_I}+g^*_{Q_{[p]-I}}$, so by the induction hypothesis, either there exists $Z_{A,I} \subseteq E(G_A)$ with $\lvert Z_{A,I} \rvert \leq f_{Q_I,g}(k_A,\sum_{v \in S_A}\gamma_A(v))$ such that $G_A-Z_{A,I}$ does not contain an $Q_I$-immersion, or there exists $Z_{B,I} \subseteq E(G_B)$ with $\lvert Z_{B,I} \rvert \leq f_{Q_{[p]-I},g}(k_B,\sum_{v \in S_B}\gamma_B(v))$ such that $G_B-Z_{B,I}$ does not contain an $Q_{[p]-I}$-immersion.

Note that $k_A = kg(k,\sum_{v \in S_A}\gamma_A(v)) \geq \sum_{v \in S_A}\gamma_A(v)$ and $k_B=kg(k,\sum_{v\in S_B}\gamma_B(v)) \geq \sum_{v \in S_B}\gamma_B(v)$ since $g$ is $H$-legal.
Therefore, for every $I$ with $\emptyset \subset I \subset [p]$, there exists $Z_I \subseteq E(G)$ with $\lvert Z_I \rvert \leq f_{Q_I,g}(k_A, \allowbreak \sum_{v \in S_A}\gamma_A(v)) + f_{Q_{[p]-I},g}(k_B,\sum_{v \in S_B}\gamma_B(v)) \leq \sum_{i=0}^{k_A}f_{Q_I,g}(k_A,i) + \sum_{i=0}^{k_B}f_{Q_{[p]-I},g}(k_B,i) \allowbreak \leq f'_{H,g}(kg(k, \allowbreak \bar{r} +\theta))$ such that either $G_A-Z_I$ has no $Q_I$-immersion or $G_B-Z_I$ has no $Q_{[p]-I}$-immersion.

Define $Z = Z_0 \cup \bigcup_{\emptyset \subset I \subset [p]} Z_I$.
Note that $\lvert Z \rvert \leq \lvert Z_0 \rvert+(2^{p}-2) \cdot f'_{H,g}(kg(k,\bar{r} +\theta))\leq \lvert Z_0 \rvert+(2^{h}-2)f'_{H,g}(kg(k,\bar{r} +\theta))$, since $p \leq h$.

Suppose that $G-(W \cup Z)$ contains an $H$-immersion.
Since $Z$ contains all edges between $A$ and $B$, and $G[A]-(W \cup Z_0)$ and $G[B]-(W \cup Z_0)$ do not contain $H$-immersions, there exists $I$ with $\emptyset \subset I \subset [p]$ such that $G[A]-Z$ contains a $Q_I$-immersion and $G[B]-Z$ contains a $Q_{[p]-I}$-immersion, contradicting the existence of $Z_I$.
This proves the claim.
$\Box$

\smallskip

\noindent{\bf Claim 3:} There exists no edge-cut $[A,B]$ of $G$ of order less than $g(k,1+\bar{r})$ such that $G[A]-S$ contains an $H$-immersion and $G[B]-S$ contains an $H$-immersion.

\noindent{\bf Proof of Claim 3:}
Suppose to the contrary that there exists an edge-cut $[A,B]$ of $G$ of order less than $g(k,1+\bar{r})$ such that $G[A]-S$ contains an $H$-immersion and $G[B]-S$ contains an $H$-immersion.
Note that $\deg_{G_A}(v_B) < g(k,1+\bar{r})$ and $\deg_{G_B}(v_A)<g(k,1+\bar{r})$.
Since both $G[A]$ and $G[B]$ contain $H$-immersions, $A \neq \emptyset \neq B$, so $G_A$ and $G_B$ are nearly 3-cut free by Claim 1.

Since $G[B]-S$ contains an $H$-immersion, $G_A-S_A = G[A]-S$ does not contain $k-\bar{r} -1$ edge-disjoint $H$-immersions, for otherwise $G-S$ contains $k-\bar{r}$ edge-disjoint $H$-immersions, contradicting (ii).
Define $\gamma_A: S_A \rightarrow {\mathbb N}$ to be the function such that $\gamma_A(v_B) = 1+\sum_{v \in S \cap B}\gamma(v)$, and $\gamma_A(x)=\gamma(x)$ for every $x \in S \cap A$.
So $G_A-S_A$ does not contain $k-\bar{r} -1 = k-\sum_{v \in S_A}\gamma_A(v)$ edge-disjoint $H$-immersions.
Furthermore, every vertex in $S_A$ has degree in $G_A$ less than $\max\{g(k,\bar{r}), g(k,1+\bar{r})\} = g(k,1+\bar{r}) = g(k,\sum_{v\in S_A}\gamma_A(v))$ since $g$ is $H$-legal.
Hence the tuple $(g,G_A,k,S_A,\gamma_A)$ satisfies (i).
Since $k-\sum_{v \in S_A}\gamma_A(v) = k-1-\bar{r} <k-\bar{r}$, by (iii), $(g,G_A,k,S_A,\gamma_A)$ does not satisfy (ii).
So there exists $Z_A \subseteq E(G_A)$ with $\lvert Z_A \rvert \leq f(k,\sum_{v \in S_A}\gamma_A(v)) = f(k,\bar{r} +1)$ such that $G_A-Z_A$ does not contain an $H$-immersion.

Similarly, there exists $Z_B \subseteq E(G_B)$ with $\lvert Z_B \rvert \leq f(k,\bar{r} +1)$ such that $G_B-Z_B$ does not contain an $H$-immersion.
Note that every edge of $G_A$ incident with $v_B$ is an edge between $A$ and $B$.
So $Z_A$ is a subset of $E(G)$.
Similarly, $Z_B$ is a subset of $E(G)$.

Let $Z'$ be the set of edges of $G$ with one end in $A$ and one end in $B$.
Define $Z_0 = Z_A \cup Z_B \cup Z'$.
Note that $\lvert Z_0 \rvert \leq 2f(k, \bar{r} +1) + g(k,1+\bar{r})$.

Since $G[A]-Z_0$ is a subgraph of $G_A-Z_A$ and $G[B]-Z_0$ is a subgraph of $G_B-Z_B$, we know that $G[A]-Z_0$ and $G[B]-Z_0$ do not contain $H$-immersions.
Note that $Z_0$ contains all edges between $A$ and $B$.
Applying Claim 2 by taking $\theta=g(k,1+\bar{r})$ and $W=\emptyset$, we know that there exists $Z$ with $Z_0 \subseteq Z \subseteq E(G)$ and $\lvert Z \rvert \leq \lvert Z_0 \rvert + (2^h-2)f'_{H,g}(kg(k,\bar{r} +g(k,1+\bar{r}))) \leq f(k, \bar{r})$ such that $G-Z$ has no $H$-immersion.
Hence $(g,G,k,S,\gamma)$ does not satisfy (ii), a contradiction.
$\Box$

\smallskip

\noindent{\bf Claim 4:} For every edge-cut $[A,B]$ of $G$ of order less than $g(k,1+\bar{r})$, exactly one of $G[A]-S$ or $G[B]-S$ contains an $H$-immersion.

\noindent{\bf Proof of Claim 4:}
Suppose to the contrary that this claim does not hold.
So there exists an edge-cut $[A,B]$ of $G$ of order less than $g(k,1+\bar{r})$ such that $G[A]-S$ and $G[B]-S$ do not contain $H$-immersions by Claim 3.
Applying Claim 2 by taking $\theta=g(k,1+\bar{r})$, $W=S$ and $Z_0$ to be the set of the edges between $A$ and $B$, we obtain $Z \subseteq E(G)$ with $\lvert Z \rvert \leq g(k,1+\bar{r})+(2^h-2)f'_{H,g}(kg(k,\bar{r} +g(k,1+\bar{r})))$ such that $G-(S \cup Z)$ has no $H$-immersion.

Let $Z'$ be the union of $Z$ and the set of edges incident with vertices in $S$.
Since $\lvert S \rvert \leq \bar{r}  \leq k-1$ and every vertex in $S$ has degree less than $g(k,\bar{r})$, $\lvert Z' \rvert \leq \lvert Z \rvert + (k-1)(g(k,\bar{r})-1) \leq kg(k,1+\bar{r})+(2^h-2)f'_{H,g}(kg(k,\bar{r} +g(k,1+\bar{r}))) \leq f(k,\bar{r})$.
Since $H$ has no isolated vertices, $G-Z'$ has no $H$-immersion, contradicting (ii).
$\Box$

\smallskip

Define $\E$ to be the collection of edge-cuts of $G$ such that $[A,B] \in \E$ if and only if $[A,B]$ has order less than $g(k,1+\bar{r})$ and $G[B]-S$ contains an $H$-immersion.
Note that Claim 4 implies that $G[A]-S$ does not contain an $H$-immersion for every $[A,B] \in \E$.

\noindent{\bf Claim 5:} $\E$ is an edge-tangle in $G$ of order $g(k,1+\bar{r})$.

\noindent{\bf Proof of Claim 5:} 
Claim 4 implies that $\E$ satisfies (E1).

Suppose that $\E$ does not satisfy (E2).
So there exist edge-cuts $[A_1,B_1],[A_2,B_2],[A_3,B_3] \in \E$ with $B_1 \cap B_2 \cap B_3 = \emptyset$.
Hence $\{A_1,B_1 \cap A_2, B_1 \cap B_2 \cap A_3\}$ is a partition of $V(G)$.
Let $[C_1,D_1]=[A_1 \cup (B_1 \cap A_2), B_1 \cap B_2 \cap A_3]$.
Note that $G[C_1]-(A_1 \cup S) \subseteq G[B_1 \cap A_2]-S \subseteq G[A_2]-S$.
Since $[A_2,B_2] \in \E$, $G[C_1]-(A_1 \cup S)$ does not contain an $H$-immersion.
Since $[A_3,B_3] \in \E$, $G[D_1]-(A_1 \cup S) \subseteq G[A_3]-(A_1 \cup S)$ does not contain an $H$-immersion.
Note that every edge between $C_1,D_1$ is either between $A_1,B_1$ or between $A_2,B_2$.
So the order of $[C_1,D_1]$ is less than $2g(k,1+\bar{r})$.
Applying Claim 2 by taking $\theta=2g(k,1+\bar{r})$, $W=A_1 \cup S$, $[A,B]=[C_1,D_1]$ and $Z_0$ to be the set of all edges between $C_1$ and $D_1$, there exists $Z_1^* \subseteq E(G)$ with $\lvert Z_1^* \rvert \leq 2g(k,1+\bar{r}) + (2^h-2)f'_{H,g}(kg(k,\bar{r} +2g(k,1+\bar{r})))$ such that $G-(A_1 \cup S \cup Z_1^*)$ has no $H$-immersion.
Hence $G[B_1]-(S \cup Z_1^*) = G-(A_1 \cup S \cup Z_1^*)$ does not contain an $H$-immersion.
Since $[A_1,B_1] \in \E$, $G[A_1]-(S \cup Z_1^*)$ does not contain an $H$-immersion.
Applying Claim 2 by taking $\theta=g(k,1+\bar{r})$, $W=S$, $[A,B]=[A_1,B_1]$ and $Z_0$ to be the union of $Z_1^*$ and the set of all edges between $A_1,B_1$, there exists $Z_2^* \subseteq E(G)$ with $\lvert Z_2^* \rvert \leq (\lvert Z_1^* \rvert+g(k,1+\bar{r})) + (2^h-2)f'_{H,g}(kg(k,\bar{r} +g(k,1+\bar{r}))) \leq 3g(k,1+\bar{r}) + (2^h-2)f'_{H,g}(kg(k,\bar{r} +2g(k,1+\bar{r})))+(2^h-2)f'_{H,g}(kg(k,\bar{r} +g(k,1+\bar{r}))) \leq 3g(k,1+\bar{r}) + 2(2^h-2)f'_{H,g}(kg(k,\bar{r} +2g(k,1+\bar{r})))$ such that $G-(S \cup Z_2^*)$ does not contain an $H$-immersion.
Let $Z_3^*$ be the union of $Z_2^*$ and the set of all edges of $G$ incident with $S$.
Since $H$ has no isolated vertices, $G-Z_3^*$ does not contain an $H$-immersion.
Note that $\lvert Z_3^* \rvert \leq \lvert Z_2^* \rvert + \lvert S \rvert (g(k,\bar{r})-1) \leq (\bar{r} +3)g(k,1+\bar{r}) + 2(2^h-2)f'_{H,g}(kg(k,\bar{r} +2g(k,1+\bar{r}))) \leq f(k,\bar{r})$.
It contradicts (ii).
So $\E$ satisfies (E2).

Finally, suppose that there exists $[A,B] \in \E$ such that there are less than $g(k,1+\bar{r})$ edges incident with $B$, then $G[B]-(E(G[B]) \cup S)$ has no $H$-immersion.
Since $[A,B] \in \E$, $G[A]-(E(G[B]) \cup S)=G[A]-S$ has no $H$-immersion.
Applying Claim 2 by taking $\theta=g(k,1+\bar{r})$, $W=S$, and $Z_0$ to be the union of $E(G[B])$ and the set of edges between $A,B$, we know there exists $Z \subseteq E(G)$ with $\lvert Z \rvert \leq 2g(k,1+\bar{r}) + (2^h-2)f'_{H,g}(kg(k,\bar{r}+g(k,1+\bar{r})))$ such that $G-(Z \cup S)$ has no $H$-immersion.
Let $Z^*$ be the union of $Z$ and the set of all edges of $G$ incident with $S$.
Then $G-Z^*$ does not contain an $H$-immersion.
But $\lvert Z^* \rvert \leq 2g(k,1+\bar{r}) + (2^h-2)f'_{H,g}(kg(k,\bar{r}+g(k,1+\bar{r}))) + \bar{r}(g(k,\bar{r})-1) \leq (\bar{r} + 2)g(k,1+\bar{r}) + (2^h-2)f'_{H,g}(kg(k,\bar{r}+g(k,1+\bar{r}))) \leq f(k,\bar{r})$, contradicting (ii).
Hence $\E$ satisfies (E3).
$\Box$

\smallskip

Let $T$ be the tree and $\P=\{X_t: t \in V(T)\}$ the partition of $V(G)$ satisfying Lemma \ref{decomp into 4-edge-conn}.
For each $t \in (T)$, we call $X_t$ the bag at $t$.
For each edge $e \in E(T)$, there exists an edge-cut $[A_e,B_e]$ of $G$ such that each $A_e$ and $B_e$ is the union of the bags of the vertices in a component of $T-e$.
So $[A_e,B_e]$ has order at most three and the edges between $A_e$ and $B_e$ are the parallel edges with the same ends.
Since $\E$ is an edge-tangle of order greater than three, $[A_e,B_e] \in \E$ or $[B_e,A_e] \in \E$ but not both.
If $[A_e,B_e] \in \E$, then we direct $e$ such that $B_e$ contains the bag of the head of $e$; otherwise, we direct $e$ in the opposite direction.
Hence, we obtain an orientation of $T$ and there exists a vertex $t^*$ of $T$ of out-degree zero.

\noindent{\bf Claim 6:} There exist a set $R$ of loops of $G[X_{t^*}]$ with $\lvert R \rvert \leq (k-1)hd$ and a set $U \subseteq E(T)$ with $\lvert U \rvert \leq (k-1)hd$ such that every edge in $U$ is incident with $t^*$, and for every $H$-immersion $\Pi=(\pi_V,\pi_E)$ in $G$, one of the following holds.
	\begin{itemize}
		\item $\Pi(H)$ intersects $S$.
		\item $\Pi(H)$ contains a non-loop edge of $G[X_{t^*}]$.
		\item $\Pi(H)$ contains an edge in $R$.
		\item $V(\Pi(H)) \cap X_{t^*} \neq \emptyset$ and there exists $e \in U$ such that $V(\Pi(H)) \cap A_e \neq \emptyset$.
	\end{itemize}

\noindent{\bf Proof of Claim 6:}
Recall that for any edge $e$ of $T$, the edges between $[A_e,B_e]$ are parallel edges.
So for every $H$-immersion $\Pi=(\pi_V,\pi_E)$ in $G$ and every edge $x$ of $H$, since $\pi_E(x)$ is a path or a cycle, there are at most two edges $e$ of $T$ incident with $t^*$ such that $V(\pi_E(x)) \cap A_e \neq \emptyset$.
Therefore, for every $H$-immersion $\Pi=(\pi_V,\pi_E)$ in $G$, there exists a set $W_\Pi$ of edges of $T$ incident with $t^*$ with $\lvert W_\Pi \rvert \leq 2\lvert E(H) \rvert \leq hd$ such that $\Pi(H) \cap G[\bigcap_{e \in W_\Pi}B_e]-X_{t^*}=\emptyset$, and for every $e \in W_\Pi$, $V(\Pi(H)) \cap A_e \neq \emptyset$. 
In addition, for every $H$-immersion $\Pi=(\pi_V,\pi_E)$ in $G$, $\Pi(H)$ contains a loop $e'$ of $G$ only if there exists a loop $e$ of $H$ such that $\pi_E(e)=e'$.
So for every $H$-immersion $\Pi=(\pi_V,\pi_E)$ in $G$, there exists a set $R_\Pi$ of loops of $G$ incident with $X_{t^*}$ with $\lvert R_\Pi \rvert \leq \lvert E(H) \rvert \leq hd$ such that $R_\Pi$ consists of the loops of $G$ incident with $X_{t^*}$ contained in $\Pi(H)$.

Let ${\mathcal C}$ be a maximal collection of $H$-immersions in $G-S$ such that
	\begin{itemize}
		\item for every member $\Pi$ of ${\mathcal C}$, $\Pi(H)$ does not contain any non-loop edge of $G[X_{t^*}]$, and 
		\item for distinct members $\Pi_1,\Pi_2$ of ${\mathcal C}$, $R_{\Pi_1} \cap R_{\Pi_2} = \emptyset$ and $W_{\Pi_1} \cap W_{\Pi_2}=\emptyset$.
	\end{itemize}
Note that members of ${\mathcal C}$ are pairwise edge-disjoint $H$-immersions in $G-S$.
So $\lvert {\mathcal C} \rvert < k-\bar{r}  \leq k$.

Define $R=\bigcup_{\Pi \in {\mathcal C}}R_\Pi$ and define $U=\bigcup_{\Pi \in {\mathcal C}}W_\Pi$.
Hence $\lvert R \rvert \leq (k-1)hd$ and $\lvert U \rvert \leq (k-1)hd$.

Let $\Pi$ be an $H$-immersion in $G$.
We may assume that $\Pi$ does not satisfy the first two conclusions of this claim, for otherwise we are done.
So $\Pi$ is an $H$-immersion in $G-S$.

Suppose that $V(\Pi(H)) \cap X_{t^*} = \emptyset$.
Since $\Pi(H) \cap G[\bigcap_{e \in W_\Pi}B_e]-X_{t^*}=\emptyset$, $\Pi(H) \subseteq G[\bigcup_{e \in W_\Pi}A_e]$.
Hence $\Pi$ is an $H$-immersion in $G[\bigcup_{e \in W_\Pi}A_e]-S$.
So $[\bigcup_{e \in W_\Pi}A_e, \bigcap_{e \in W_\Pi}B_e] \not \in \E$ by the definition of $\E$.
Since $\lvert W_\Pi \rvert \leq hd$, the order of $[\bigcup_{e \in W_\Pi}A_e, \bigcap_{e \in W_\Pi}B_e]$ is at most $3hd$ which is less than the order of $\E$.
Since $[A_e,B_e] \in \E$ for each $e \in W_\Pi$, $[\bigcup_{e \in W_\Pi}A_e, \bigcap_{e \in W_\Pi}B_e] \in \E$ by Claim 5 and Lemma \ref{easy edge-tangle}, a contradiction.

Hence $V(\Pi(H)) \cap X_{t^*} \neq \emptyset$.
Suppose $R_\Pi=\emptyset$ and $W_\Pi = \emptyset$.
Since $W_\Pi = \emptyset$, $\Pi(H) \subseteq G[X_{t^*}]$.
Since $R_\Pi = \emptyset$ and $\Pi(H) \subseteq G[X_{t^*}]$, $\Pi(H)$ does not contain any loop of $G$.
Since $H$ has no isolated vertex, $\Pi(H)$ contains a non-loop edge of $G[X_{t^*}]$, so $\Pi$ satisfies the second conclusion of this claim, a contraction.

Hence either $R_\Pi \neq \emptyset$ or $W_\Pi \neq \emptyset$.
If $\Pi \in {\mathcal C}$ and $R_\Pi \neq \emptyset$, then $\Pi$ satisfies the third conclusion of this claim.
If $\Pi \in {\mathcal C}$ and $W_\Pi \neq \emptyset$, then $\Pi$ satisfies the fourth conclusion of this claim.
So we may assume that $\Pi \not \in {\mathcal C}$.
By the maximality of ${\mathcal C}$, there exists $\Pi' \in {\mathcal C}$ such that either $R_\Pi \cap R_{\Pi'} \neq \emptyset$ or $W_\Pi \cap W_{\Pi'} \neq \emptyset$.
If $R_\Pi \cap R_{\Pi'} \neq \emptyset$, then $\Pi(H)$ contains an edge in $R$, so $\Pi$ satisfies the third conclusion of this claim.
So we may assume that $W_\Pi \cap W_{\Pi'} \neq \emptyset$.
Then there exists $e \in W_{\Pi'} \subseteq U$ such that $V(\Pi(H)) \cap A_e \neq \emptyset$.
So $\Pi$ satisfies the fourth conclusion of this claim.
$\Box$

\smallskip

\noindent{\bf Claim 7:} For every $[A,B] \in \E$ of order less than $g(k,1+\bar{r})-3(k-1)hd$, $G[B \cap X_{t^*}]$ contains at least $f(k,\bar{r})-kg(k, 1+\bar{r})-(k-1)hd-(2^h-2)f'_{H,g}(kg(k,\bar{r} +g(k,1+\bar{r})))$ non-loop edges.

\noindent{\bf Proof of Claim 7:}
Let $m=f(k,\bar{r})-kg(k, 1+\bar{r})-(k-1)hd-(2^h-2)f'_{H,g}(kg(k,\bar{r} +g(k,1+\bar{r})))$.
Suppose to the contrary that there exists $[A,B] \in \E$ of order less than $g(k,1+\bar{r})-3(k-1)hd$ such that $G[B \cap X_{t^*}]$ contains less than $m$ non-loop edges.

Let $U$ be the set of edges of $T$ incident with $t^*$ and $R$ the set of loops of $G[X_{t^*}]$ mentioned in Claim 6.
Note that $\lvert U \rvert \leq (k-1)hd$ and $\lvert R \rvert \leq (k-1)hd$.
Let $A' = A \cup \bigcup_{e \in U} A_e$ and $B' = B \cap \bigcap_{e \in U} B_e$.
Note that the order of $[A',B']$ is at most $\lvert [A,B] \rvert + 3(k-1)hd < g(k,1+\bar{r})$.
By Lemma \ref{easy edge-tangle}, $[A', B'] \in \E$.
By Claim 6, for every $H$-immersion $\Pi=(\pi_V,\pi_E)$ in $G-S$, one of the following holds.
	\begin{itemize}
		\item $\Pi(H)$ contains a non-loop edge of $G[X_{t^*}]$ or a loop in $R$. 
		\item $V(\Pi(H)) \cap A' \neq \emptyset$.
	\end{itemize}
Let $Z_0$ be the set consisting of all non-loop edges in $G[B \cap X_{t^*}]$ and all edges of $G$ between $A'$ and $B'$.
Let $Z=Z_0 \cup R$.
Then $G[B']-(S \cup Z)$ has no $H$-immersion.
In addition, by the definition of $\E$, $G[A']-S$ has no $H$-immersion.
Applying Claim 2 by taking $(\theta,W,[A,B],Z_0)=(g(k,1+\bar{r}),S,[A',B'],Z)$, there exists $Z' \subseteq E(G)$ with $\lvert Z' \rvert \leq \lvert Z \rvert + (2^h-2)f'_{H,g}(kg(k,\bar{r} +g(k,1+\bar{r}))) \leq (m + g(k,1+\bar{r}) + (k-1)hd) + (2^h-2)f'_{H,g}(kg(k,\bar{r} +g(k,1+\bar{r}))) \leq f(k,\bar{r}) - (k-1)g(k,\bar{r})$ such that $G-(S \cup Z')$ has no $H$-immersion.

Let $Z_S$ be the set of edges incident with vertices $S$.
So $\lvert Z_S \rvert \leq (k-1)g(k,\bar{r})$.
Let $Z^* = Z' \cup Z_S$.
Then $\lvert Z^* \rvert \leq f(k,\bar{r})$ and $G-Z^*$ has no $H$-immersion, contradicting (ii).
$\Box$

\smallskip

Since $f(k,\bar{r})-kg(k, 1+\bar{r})-(k-1)hd-(2^h-2)f'_{H,g}(kg(k,\bar{r} +g(k,1+\bar{r}))) \geq 2$, Claim 7 implies that $X_{t^*}$ contains at least two vertices and hence $G[X_{t^*}]$ does not have an edge-cut of order three.

For every vertex $v$ in $X_{t^*}$, define $Q_v=\{u \in V(G)-X_{t^*}:$ every path in $G$ from $u$ to $X_{t^*}$ contains $v\}$.
Note that $Q_v$ is empty if $N_G(v) \subseteq X_{t^*}$.
Furthermore, since $G$ is connected, by Lemma \ref{decomp into 4-edge-conn}, for every $u \in V(G)-X_{t^*}$, there exists a unique $v \in X_{t^*}$ such that $u \in Q_v$.

Define $\E'$ to be the set of edge-cuts $[A',B']$ of $G[X_{t^*}]$ of order less than $g(k,1+\bar{r})-3(k-1)hd$ such that $[A',B'] \in \E'$ if and only if $[A' \cup \bigcup_{v \in A'} Q_v, B' \cup \bigcup_{v \in B'} Q_v] \in \E$.

\noindent{\bf Claim 8:} $\E'$ is an edge-tangle of order $g(k,1+\bar{r})-3(k-1)hd$ in $G[X_{t^*}]$.

\noindent{\bf Proof of Claim 8:}
For every edge-cut $[A',B']$ of $G[X_{t^*}]$, the order of $[A' \cup \bigcup_{v \in A'} Q_v, B' \cup \bigcup_{v \in B'} Q_v]$ equals the order of $[A',B']$.
Since $\E$ is an edge-tangle, either $[A' \cup \bigcup_{v \in A'} Q_v, B' \cup \bigcup_{v \in B'} Q_v] \in \E$ or $[B' \cup \bigcup_{v \in B'} Q_v,A' \cup \bigcup_{v \in A'} Q_v] \in \E$.
So $\E'$ satisfies (E1).

If there exist $[A_i',B_i'] \in \E'$ for $i\in [3]$ such that $A_1' \cup A_2' \cup A_3' = X_{t^*}$, then $[A_i' \cup \bigcup_{v \in A_i'}Q_v, B_i' \cup \bigcup_{v \in B_i'}Q_v] \in \E$ for $i \in [3]$, but $A_1' \cup A_2' \cup A_3' \cup \bigcup_{v \in A_1' \cup A_2' \cup A_3'} Q_v = V(G)$, a contradiction.
So $\E'$ satisfies (E2).

In addition, for every edge-cut $[A',B'] \in \E'$, the number of edges of $G[X_{t^*}]$ incident with $B'$ is at least $\lvert E(G[B']) \rvert = \lvert E(G[(B' \cup \bigcup_{v \in B'} Q_v) \cap X_{t^*}]) \rvert \geq f(k,\bar{r})-kg(k, 1+\bar{r})-(k-1)hd-(2^h-2)f'_{H,g}(kg(k,\bar{r} +g(k,1+\bar{r}))) \geq g(k,1+\bar{r})-3(k-1)hd$ by Claim 7.
So $\E'$ satisfies (E3).
$\Box$

\smallskip

By Claim 8, $\E'$ is an edge-tangle of order $g(k, 1+\bar{r})-3(k-1)hd \geq w_{\ref{edge-tangle in 4-edge-connected control edge-minor}}(w_{\ref{isolating immersion}}(H,k+ \bar{r} \cdot g(k,\bar{r})),\theta_{\ref{isolating immersion}}(H,k+ \bar{r} \cdot g(k,\bar{r})))$ in $G[X_{t^*}]$, where the last inequality follows from the assumption that $g$ is $H$-legal.

Define $\E_k$ and $\E'_k$ to be the subsets of $\E$ and $\E'$ consisting of edge-cuts of order less than $\theta_{\ref{isolating immersion}}(H,k+ \bar{r} \cdot g(k,\bar{r}))$, respectively.
So $\E_k$ and $\E'_k$ are edge-tangles of order $\theta_{\ref{isolating immersion}}(H,k+ \bar{r} \cdot g(k,\bar{r}))$ in $G$ and $G[X_{t^*}]$, respectively.
Let $w_k = w_{\ref{isolating immersion}}(H,k+ \bar{r} \cdot g(k,\bar{r}))$.

Since $G[X_{t^*}]$ does not have an edge-cut of order three, by Theorem \ref{edge-tangle in 4-edge-connected control edge-minor}, $\E'_k$ controls a $K_{w_k}$-thorns $\alpha$ in $G[X_{t^*}]$.
Since $\alpha$ is in $G[X_{t^*}] \subseteq G$, $\E_k$ controls $\alpha$.

Since $G-S$ does not contain $k-\bar{r}$ edge-disjoint $H$-immersions and every vertex in $S$ has degree in $G$ at most $g(k,\bar{r})-1$, $G$ does not contain $k-\bar{r}  + \lvert S \rvert (g(k,\bar{r})-1) \leq k-\bar{r}  + \bar{r} \cdot g(k,\bar{r}) \leq k+ \bar{r} \cdot g(k,\bar{r})$ edge-disjoint $H$-immersions.
Since $\E_k$ is an edge-tangle in $G$ of order $\theta_{\ref{isolating immersion}}(H,k+ \bar{r} \cdot g(k,\bar{r}))$ controlling a $K_{w_k}$-thorns, by Lemma \ref{isolating immersion}, there exist $Z^* \subseteq E(G)$ with $\lvert Z^* \rvert \leq \xi_k$ and $[A,B] \in \E_k-Z^* \subseteq \E-Z^*$ of order zero such that $G[B]-Z^*$ has no $H$-immersion, where $\xi_k = \xi_{\ref{isolating immersion}}(H,k+ \bar{r} \cdot g(k,\bar{r}))$.
Note that $\xi_k < \theta_{\ref{isolating immersion}}(H,k+ \bar{r} \cdot g(k,\bar{r})) \leq g(k,1+\bar{r})$ by Lemma \ref{isolating immersion} and the assumption that $g$ is $H$-legal.

In addition, $G[A]-S$ does not contain an $H$-immersion since $[A,B] \in \E$.
So every $H$-immersion in $G[A]$ intersects an edge in $Z_S$, where $Z_S$ is the set of edges of $G$ incident with $S$.
Hence $G[A]-(Z^* \cup Z_S)$ and $G[B]-(Z^* \cup Z_S)$ do not contain $H$-immersions.
Since $[A,B] \in \E$, $[A,B]$ is an edge-cut of $G$ of order less than $g(k,1+\bar{r})$.
By Claim 2, there exists $Z^{**}$ with $Z^{**} \subseteq E(G)$ and $\lvert Z^{**} \rvert \leq (\lvert Z^* \rvert + \lvert Z_S \rvert + g(k,1+\bar{r})) + (2^h-2)f'_{H,g}(kg(k,\bar{r} +g(k,1+\bar{r}))) \leq (2g(k,1+\bar{r}) + (k-1)g(k,\bar{r})) + (2^h-2)f'_H(kg(k,\bar{r} +g(k,1+\bar{r}))) \leq f(k,\bar{r})$ such that $G-Z^{**}$ has no $H$-immersion, contradicting (ii).
This completes the proof.
\end{pf}

\bigskip

Now we drop the requirement of having no isolated vertices from Lemma \ref{EP not isolated}.
Lemma \ref{EP immersion nearly 4-edge-connected} proves Theorem \ref{main intro} for the case that $G$ has only one component, as every 4-edge-connected graph is nearly 3-cut-free.

\begin{lemma} \label{EP immersion nearly 4-edge-connected}
For every graph $H$, there exists a function $f: {\mathbb N} \rightarrow {\mathbb N}$ such that for every nearly 3-cut free graph $G$ and every positive integer $k$, either $G$ contains $k$ edge-disjoint $H$-immersions, or there exists $Z \subseteq E(G)$ with $\lvert Z \rvert \leq f(k)$ such that $G-Z$ contains no $H$-immersion.
\end{lemma}

\begin{pf}
Let $H$ be a graph.
Note that when $H$ has no edge, every graph $G$ contains at least one $H$-immersion if and only if $G$ contains at least $\lvert V(H) \rvert$ vertices if and only if $G$ contains arbitrarily many edge-disjoint $H$-immersions.
So this theorem holds if $H$ has no edge.
	
Hence we may assume that $H$ contains at least one edge.
Let $H'$ be the graph obtained from $H$ by deleting all isolated vertices.
So $H'$ is a graph with at least one edge and with no isolated vertices.
For every positive integer $k$, define $f(k)$ to be the number $f_{\ref{EP not isolated}}(k,0)$, where $f_{\ref{EP not isolated}}$ is the function mentioned in Theorem \ref{EP not isolated} by taking $H$ to be $H'$ and further taking $g=g^*$.
We apply Theorem \ref{EP not isolated} by further taking $S=\emptyset$ and $\gamma$ to be a function with empty domain, we know that for every nearly 3-cut free graph $G$ and positive integer $k$, either $G$ contains $k$ edge-disjoint $H'$-immersions, or there exists $Z \subseteq E(G)$ with $\lvert Z \rvert \leq f(k)$ such that $G-Z$ does not contain an $H'$-immersion.

We shall prove that $f$ is a function satisfying the conclusion of this lemma.

Let $G$ be a nearly 3-cut free graph.
If $\lvert V(G) \rvert < \lvert V(H) \rvert$, then clearly $G$ does not contain an $H$-immersion, and we are done by choosing $Z=\emptyset$.
So we may assume that $\lvert V(G) \rvert \geq \lvert V(H) \rvert$.
Hence, for every $W \subseteq E(G)$ and every $H'$-immersion $(\pi'_V,\pi'_E)$ of $G-W$, we can extend $\pi'_V$ to an injection $\pi_V$ with domain $V(H)$ by further mapping isolated vertices of $H$ to some vertices of $G-\pi'_V(V(H'))$ such that $(\pi_V,\pi'_E)$ is an $H$-immersion in $G-W$ with $E(\pi'_E(E(H')))=E(\pi'_E(E(H)))$.
Therefore, for every $W \subseteq E(G)$ and every integer $k$, $G-W$ contains $k$ edge-disjoint $H$-immersions if and only if $G-W$ contains $k$ edge-disjoint $H'$-immersions.

Now let $k$ be a positive integer.
If $G$ does not contain $k$ edge-disjoint $H$-immersions, then $G$ does not contain $k$ edge-disjoint $H'$-immersions, so there exists $Z \subseteq E(G)$ with $\lvert Z \rvert \leq f_{\ref{EP not isolated}}(k,0)=f(k)$ such that $G-Z$ has no $H'$-immersion by Theorem \ref{EP not isolated}.
But it implies that $G-Z$ has no $H$-immersions.
This proves the lemma.
\end{pf}

\bigskip

Theorem \ref{main intro} is an immediate corollary of the following theorem.

\begin{theorem} \label{EP comp near 4-edge-conn}
For every graph $H$, there exists a function $f: {\mathbb N} \rightarrow {\mathbb N}$ such that for every graph $G$ whose every component is nearly 3-cut free and for every positive integer $k$, either $G$ contains $k$ edge-disjoint $H$-immersions, or there exists $Z \subseteq E(G)$ with $\lvert Z \rvert \leq f(k)$ such that $G-Z$ contains no $H$-immersion.
\end{theorem}

\begin{pf}
Let $H$ be a graph, and let $c$ be the number of components of $H$.
We define the following.
	\begin{itemize}
		\item For every graph $R$, define $f_R$ to be the function $f$ mentioned in Lemma \ref{EP immersion nearly 4-edge-connected} by taking $H=R$.
		\item For every $i\in [c]$, let $\F_i$ be the set of graphs consisting of $i$ components of $H$.
		\item For any positive integers $m \geq 2$ and $n$, let $f_1(n) = (n-1) \cdot \max\{f_R(n): R \in \F_1\}$ and let $f_m(n)=m^{nm}f_{m-1}(n) + (nm-1) \cdot \max\{f_R(n): R \in \F_m\}$.
		\item Define $f: {\mathbb N} \rightarrow {\mathbb N}$ to be the function such that $f(x)=f_c(x)$ for every $x \in {\mathbb N}$.
	\end{itemize}

We claim that for every $m \in [c]$, for every graph $W \in \F_m$, for every positive integer $k$ and for every graph $G$ whose every component is nearly 3-cut free, either $G$ contains $k$ edge-disjoint $W$-immersions, or there exists $Z \subseteq E(G)$ with $\lvert Z \rvert \leq f_m(k)$ such that $G-Z$ does not contain a $W$-immersion.
Note that this claim implies this theorem as $H \in \F_c$.
We shall prove this claim by induction on $m$.

Let $m \in [c]$, $W \in \F_m$, $k$ a positive integer and $G$ a graph whose every component is nearly 3-cut free.
We assume that $G$ does not contain $k$ edge-disjoint $W$-immersions.
It suffices to show that there exists $Z \subseteq E(G)$ with $\lvert Z \rvert \leq f_m(k)$ such that $G-Z$ has no $W$-immersion.

We first assume that $m=1$.
Let $G_1,G_2,...,G_p$ be the components of $G$ containing a $W$-immersion, and let $k_i$ be the maximum number of edge-disjoint $W$-immersions in $G_i$ for each $i \in [p]$.
If $\sum_{i=1}^p k_i \geq k$, then $G$ contains $k$ edge-disjoint $W$-immersions, a contradiction.
So $\sum_{i=1}^p k_i < k$.
In particular, $p<k$.
By Lemma \ref{EP immersion nearly 4-edge-connected}, for every $i \in [p]$, there exists $Z_i \subseteq E(G_i)$ with $\lvert Z_i \rvert \leq f_W(k_i+1)$ such that $G_i-Z_i$ has no $W$-immersion.
Since $m=1$, $W$ is connected, so $G-Z$ has no $W$-immersion, where $Z=\bigcup_{i=1}^pZ_i \subseteq E(G)$.
Note that $\lvert Z \rvert \leq \sum_{i=1}^p f_W(k_i+1) \leq (k-1)f_W(k) \leq f_1(k)$.
This proves the base case of the induction.

So we may assume that $m \geq 2$ and our claim holds for every smaller $m$.

Note that $W$ has $m$ components.
Let $W_1,W_2,...,W_m$ be the components of $W$.
For every $i \in [m]$, define $S_i$ to be the set of components of $G$ containing an $W_i$-immersion.
If $\lvert S_i \rvert \geq km$ for every $i \in [m]$, then $G$ contains $km$ components $G_1,...,G_{km}$ of $G$ such that $G_{(i-1)k+j} \in S_i$ for each $i \in [m]$ and each $j \in [k]$, so $G$ contains $k$ edge-disjoint $W$-immersions, a contradiction.
Therefore, there exists $t \in [m]$ such that $\lvert S_t \rvert <km$.

Define $L$ to be the disjoint union of the components of $G$ in $S_t$, and define $R=G-V(L)$.
Note that $R$ has no $W_t$-immersion by the definition of $S_t$.
Since $m \geq 2$, by the induction hypothesis, if $L$ does not contain $k$ edge-disjoint $W_t$-immersions, then there exists $Z_t \subseteq E(L)$ with $\lvert Z_t \rvert \leq f_1(k) \leq f_m(k)$ such that $L-Z_t$ has no $W_t$-immersion, so $G-Z_t$ has no $W_t$-immersion (since $W_t$ is connected) and hence has no $W$-immersion.

So we may assume that $L$ contains $k$ edge-disjoint $W_t$-immersions.
Hence $R$ does not contain $k$ edge-disjoint $(W-V(W_t))$-immersions, for otherwise $G$ contains $k$ edge-disjoint $W$-immersions.
Note that $W-V(W_t) \in \F_{m-1}$.
By the induction hypothesis, there exists $Z_R \subseteq E(R)$ with $\lvert Z_R \rvert \leq f_{m-1}(k)$ such that $R-Z_R$ has no $(W-V(W_t))$-immersion.
So $R-Z_R$ has no $W$-immersion.
In addition, for each component $C$ of $L$, $C$ is nearly 3-cut free and has no $k$ edge-disjoint $W$-immersions, so there exists $Z_C \subseteq E(C)$ with $\lvert Z_C \rvert \leq f_W(k)$ such that $C-Z_C$ has no $W$-immersion.

Define $Z_0 = Z_R \cup \bigcup_C Z_C$, where the second union is taken over all components $C$ of $L$.
Note that the number of components of $L$ equals $\lvert S_t \rvert \leq km-1$.
Therefore, $\lvert Z_0 \rvert \leq f_{m-1}(k)+(km-1)f_W(k)$, and $R-Z_0$ and $C-Z_0$ do not contain a $W$-immersion for every component $C$ of $L$.

Let $\ell$ be the number of components of $L$.
Define $Q_0=R$, and for every $i \in [\ell]$, define $Q_i$ to be the $i$-th component of $L$.
Note that $Q_i-Z_0$ has no $W$-immersion for every $i$ with $0 \leq i \leq \ell$.

We say that $(P_0,P_1,...,P_\ell)$ is a {\it $(\ell+1)$-partition} of $[m]$ if $P_0,P_1,...,P_\ell$ are pairwise disjoint (possibly empty) proper subsets of $[m]$ with $\bigcup_{i=0}^\ell P_i = [m]$.
Since $G$ has no $k$ edge-disjoint $W$-immersions, for every $(\ell+1)$-partition $\P=(P_0,...,P_\ell)$ of $[m]$, there exists $j$ with $0 \leq j \leq \ell$ such that $Q_j$ does not contain $k$ edge-disjoint $(\bigcup_{i \in P_j}W_i)$-immersions, so there exists $Z_\P \subseteq E(Q_j)$ with $\lvert Z_\P \rvert \leq f_{\lvert P_j \rvert}(k) \leq f_{m-1}(k)$ such that $Q_j-Z_\P$ has no $(\bigcup_{i \in P_j}W_i)$-immersions.

Define $Z^*$ to be the union of $Z_0$ and $Z_\P$ over all $(\ell+1)$-partitions $\P$ of $[m]$.
Since $\ell = \lvert S_t \rvert \leq km-1$ and there are at most $m^{\ell+1}-1$ different $(\ell+1)$-partitions of $[m]$, $\lvert Z^* \rvert \leq \lvert Z_0 \rvert + (m^{\ell+1}-1)f_{m-1}(k) \leq f_{m-1}(k)+(km-1)f_W(k)+(m^{km}-1)f_{m-1}(k) \leq f_m(k)$.

To prove the theorem, it suffices to prove that $G-Z^*$ has no $W$-immersion.
Suppose to the contrary that $G-Z^*$ contains a $W$-immersion.
Since $Q_i-Z^*$ has no $W$-immersion for every $0 \leq i \leq \ell$, there exists a $(\ell+1)$-partition $\P=(P_0,P_1,...P_\ell)$ of $[m]$ such that for every $j$ with $0 \leq j \leq \ell$, $Q_j-Z^*$ contains a $(\bigcup_{i \in P_j}W_i)$-immersion.
However, it contradicts the definition of $Z_\P$.
This completes the proof.
\end{pf}

\bigskip

We remark that Kakimura and Kawarabayashi \cite{kk1} proved that for every integer $t$, there exists a function $f$ such that for every 3-minimal-cut free graph $G$ and integer $k$, either $G$ contains $k$ edge-disjoint $K_t$-immersions, or there exists $Z \subseteq E(G)$ with $\lvert Z \rvert \leq f(k)$ such that $G-Z$ has no $K_t$-immersion, where a graph is {\it 3-minimal-cut free} if it is connected and it cannot be made disconnected by deleting exactly three edges while it remains connected by deleting at most two of those three edges.
This result is a simple corollary of Theorem \ref{EP comp near 4-edge-conn} when $t \geq 3$ (and the case $t \leq 2$ is easy).
Let $G$ be a 3-minimal-cut free graph, and let $G'$ be the graph obtained from $G$ by deleting all cut-edges and loops and then deleting all resulting isolated vertices.
Note that every component of $G'$ is 2-edge-connected and 3-minimal-cut free.
If a 2-edge-connected and 3-minimal-cut free graph has an edge-cut $[A,B]$ of order three, then one can delete at most two edges between $A$ and $B$ to make the graph disconnected, but it implies that some edge in $[A,B]$ is a cut-edge of the original graph, contradicting that it is 2-edge-connected.
Hence every component of $G'$ does not contain an edge-cut of order three and hence is nearly 3-cut free.
In addition, when $t \geq 3$, the optimal solutions for packing and covering $K_t$-immersions in $G$ are the same as the optimal solutions for packing and covering $K_t$-immersions in $G'$.
Hence Theorem \ref{EP comp near 4-edge-conn} implies the result in \cite{kk1}.

Now we prove Theorem \ref{half-integral intro}.
The following is a restatement.

\begin{theorem} 
For every loopless graph $H$, there exists a function $f: {\mathbb N} \rightarrow {\mathbb N}$ such that for every positive integer $k$ and every graph $G$, either there exists $k$ $H$-half-integral immersions $(\pi^{(1)}_V,\pi^{(1)}_E), ..., (\pi^{(k)}_V,\pi^{(k)}_E)\}$ in $G$ such that for each edge $e$ of $G$, there exist at most two distinct pairs $(i,e')$ with $1 \leq i \leq k$ and $e' \in E(H)$ such that $e \in \pi^{(i)}_E(e')$, or there exists $Z \subseteq E(G)$ with $\lvert Z \rvert \leq f(k)$ such that $G-Z$ contains no $H$-half-integral immersion.
\end{theorem}

\begin{pf}
Define $f$ to be the function $f$ mentioned in Theorem \ref{EP comp near 4-edge-conn} by taking $H=H$.

Let $k$ be a positive integer.
Let $G$ be a graph, and let $G'$ be the graph obtained from $G$ by duplicating each edge.
Note that every edge-cut of $G'$ has even order.
If $[A,B]$ is an edge-cut of a component of $G'$ of order between one and three, then it has order two and the two edges between $A$ and $B$ are parallel edges with the same ends.
So every component of $G'$ is nearly 3-cut free.
By Theorem \ref{EP comp near 4-edge-conn}, either $G'$ contains $k$ edge-disjoint $H$-immersions, or there exists $Z' \subseteq E(G')$ with $\lvert Z' \rvert \leq f(k)$ such that $G'-Z'$ does not contain an $H$-immersion.

Note that since $H$ is loopless, for every $H$-immersion $(\pi_V,\pi_E)$ in $G'$ and $e \in E(H)$, $\pi_E(e)$ is a path in $G'$, so there exists no $e' \in E(G)$ such that $\pi_E(e)$ contains $e'$ and its copy in $G'$. 
So for every $H$-immersion $(\pi'_V,\pi'_E)$ in $G'$, there exists an $H$-half-integral immersion $(\pi_V,\pi_E)$ such that $\pi_V(v)=\pi_V'(v)$ for every $v \in V(H)$, and for every $e \in E(H)$, $\pi_E(e)$ consists of the edges $z$ of $G$ such that some copy of $z$ belongs to $E(\pi_E'(e))$.
Similarly, if $G'$ contains $k$ edge-disjoint $H$-immersions, then $G$ contains $k$ $H$-half-integral immersions $(\pi^{(1)}_V,\pi^{(1)}_E), ..., (\pi^{(k)}_V,\pi^{(k)}_E)$ such that for each edge $e$ of $G$, there exist at most two distinct pairs $(i,e')$ with $1 \leq i \leq k$ and $e' \in E(H)$ such that $e \in \pi^{(i)}_E(e')$, so we are done.

So we may assume that there exists $Z' \subseteq E(G')$ with $\lvert Z' \lvert \leq f(k)$ such that $G'-Z'$ has no $H$-immersion.
Then $G-Z$ has no $H$-half-integral immersion, where $Z \subseteq E(G)$ is the set consisting of the edges of $G$ having a copy in $Z'$.
Note that $\lvert Z \rvert \leq \lvert Z' \rvert$.
This completes the proof.
\end{pf}

\bigskip

\noindent{\bf Acknowledgements.}
Theorem \ref{excluding immersion} is based on part of an unpublished manuscript of the author in 2013.
The author thanks Professor Robin Thomas and Professor Paul Wollan for some inspiring discussion that led to that manuscript.
Lemma \ref{colorful vertex spider} was proved by the author when he worked on \cite{lpw}, and its proof was included in \cite{lpw} when the first version of this paper was submitted.
The author thanks Professor Luke Postle and Professor Paul Wollan for agreeing for moving the proof of Lemma \ref{colorful vertex spider} into this paper to meet a request of anonymous reviewers of this paper.
The author thanks Professor Reinhard Diestel for pointing out \cite{dhl, do2,do} and thanks anonymous reviewers for careful reading and suggestions.

\end{document}